\documentclass[leqno,12pt]{amsart} 
\usepackage{amssymb}
\theoremstyle{plain} 
\newtheorem{theorem}{\indent\sc Theorem}[section]
\newtheorem{lemma}[theorem]{\indent\sc Lemma}
\newtheorem{corollary}[theorem]{\indent\sc Corollary}
\newtheorem{proposition}[theorem]{\indent\sc Proposition}

\theoremstyle{definition} 

\newtheorem{remark}[theorem]{\indent\sc Remark}
\newtheorem{example}[theorem]{\indent\sc Example}

\begin{document}

\font\eightrm=cmr8
\font\eightit=cmti8
\font\eighttt=cmtt8
\def\bbR{\text{\bf R}}
\def\rtr{\text{\bf R}\hskip-.7pt^3}
\def\bbRP{\text{\bf R}\text{\rm P}}
\def\bbC{\text{\bf C}}
\def\bbZ{\text{\bf Z}}
\def\hyp{\hskip.5pt\vbox
{\hbox{\vrule width3ptheight0.5ptdepth0pt}\vskip2.2pt}\hskip.5pt}
\def\er{r\hskip.3pt}
\def\df{d\hskip-.8ptf}
\def\tb{T\hskip-.3pt\bs}
\def\tcb{{T_{\hskip-.6pt\cy}\bs}}
\def\tctb{{T_{\hskip-.6pt\cy(t)}\bs}}
\def\tab{{T\hskip.2pt^*\hskip-2.3pt\bs}}
\def\tacb{{T_\cy^*\hskip-1.9pt\bs}}
\def\txm{{T\hskip-.4pt_x\hskip-.9ptM}}
\def\tm{{T\hskip-.3ptM}}
\def\tn{{T\hskip-.3ptN}}
\def\tam{{T^*\!M}}
\def\fc{F\hskip-2pt_\cy}
\def\lo{\mathcal{B}}
\def\op{\mathcal{L}}
\def\bop{\mathcal{F}}
\def\vt{\mathcal{P}}
\def\vtx{\vt_{\nh x}}
\def\h{\mathcal{V}}
\def\jm{\mathcal{I}}
\def\ke{\mathcal{K}}
\def\bu{\mathcal{T}}
\def\xz{\mathcal{X}}
\def\yz{\mathcal{Y}}
\def\zz{\mathcal{Z}}
\def\q{q}
\def\bq{\hat q}
\def\p{p}
\def\w{\vt^\perp}
\def\x{v}
\def\y{y}
\def\vp{\vt^\perp}
\def\vd{\vt\hh'}
\def\vdx{\vd{}\hskip-4.5pt_x}
\def\bz{b\hh}
\def\cy{{y}}
\def\rwo{\,\hs\text{\rm rank}\hskip2.7ptW\hskip-2.7pt=\hskip-1.2pt1}
\def\rwho{\,\hs\text{\rm rank}\hskip2.2ptW^h\hskip-2.2pt=\hskip-1pt1}
\def\rw{\,\hs\text{\rm rank}\hskip2.4ptW\hskip-1.5pt}
\def\im{\varPhi}
\def\js{J}
\def\ism{H}
\def\fe{F}
\def\fy{f}
\def\dfc{dF\hskip-2.3pt_\cy\hskip.4pt}
\def\dfct{dF\hskip-2.3pt_\cy(t)\hskip.4pt}
\def\dic{d\im\hskip-1.4pt_\cy\hskip.4pt}
\def\vl{\Lambda}
\def\hy{\mathcal{V}}
\def\mv{V}
\def\bs{\varSigma}
\def\hs{\hskip.7pt}
\def\hh{\hskip.4pt}
\def\nh{\hskip-.7pt}
\def\nnh{\hskip-1pt}
\def\hrz{^{\hskip.5pt\text{\rm hrz}}}
\def\vrt{^{\hskip.2pt\text{\rm vrt}}}
\def\th{\varTheta}
\def\vg{\varGamma}
\def\my{\mu}
\def\ny{\nu}
\def\gy{\lambda}
\def\gm{\gamma}
\def\sj{\sigma}
\def\lg{\langle}
\def\rg{\rangle}
\def\lr{\lg\,,\rg}
\def\vs{vector space}
\def\rvs{real vector space}
\def\vf{vector field}
\def\tf{tensor field}
\def\tvn{the vertical distribution}
\def\dn{distribution}
\def\pt{point}
\def\tc{tor\-sion\-free connection}
\def\ea{equi\-af\-fine}
\def\rt{Ricci tensor}
\def\sy{symmetric}
\def\pde{partial differential equation}
\def\pf{projectively flat}
\def\pfs{projectively flat surface}
\def\pfc{projectively flat connection}
\def\pftc{projectively flat tor\-sion\-free connection}
\def\su{surface}
\def\sco{simply connected}
\def\psr{pseu\-\hbox{do\hs-}Riem\-ann\-i\-an}
\def\inv{-in\-var\-i\-ant}
\def\trinv{trans\-la\-tion\inv}
\def\diml{-di\-men\-sion\-al}
\def\prl{-par\-al\-lel}
\def\skc{skew-sym\-met\-ric}
\def\sky{skew-sym\-me\-try}
\def\Sky{Skew-sym\-me\-try}
\def\dbly{-dif\-fer\-en\-ti\-a\-bly}
\def\cs{conformally symmetric}
\def\cf{conformally flat}
\def\ls{locally symmetric}
\def\ecs{essentially conformally symmetric}
\def\rr{Ric\-ci-re\-cur\-rent}
\def\kf{Killing field}
\def\om{\omega}
\def\vol{\varOmega}
\def\ve{\varepsilon}
\def\kx{\kappa}
\def\mf{manifold}
\def\mfd{-man\-i\-fold}
\def\bmf{base manifold}
\def\bd{bundle}
\def\tbd{tangent bundle}
\def\ctb{cotangent bundle}
\def\bp{bundle projection}
\def\prc{pseu\-\hbox{do\hs-}Riem\-ann\-i\-an metric}
\def\prd{pseu\-\hbox{do\hs-}Riem\-ann\-i\-an manifold}
\def\Prd{pseu\-\hbox{do\hs-}Riem\-ann\-i\-an manifold}
\def\npd{null parallel distribution}
\def\pj{-pro\-ject\-a\-ble}
\def\pd{-pro\-ject\-ed}
\def\lcc{Le\-vi-Ci\-vi\-ta connection}
\def\vb{vector bundle}
\def\vbm{vec\-tor-bun\-dle morphism}
\def\kerd{\text{\rm Ker}\hskip2.7ptd}
\def\ri{{\rho}}
\def\rih{\widetilde\rho}
\def\deh{\widetilde{\text{\rm D}}}
\def\bsh{\widetilde\bs}
\def\heh{\widetilde h}
\def\scal{\text{\rm s}}
\def\naf{\nabla\!\fy}
\def\bna{\hs\overline{\nh\nabla\nh}\hs}
\def\br{\hskip3pt\overline{\hskip-3ptR\nh}\hs}
\def\ts{total space}
\def\pmb{\pi}

\renewcommand{\thepart}{\Roman{part}}
\title[Conformally symmetric manifolds]{Projectively flat surfaces,\\
null parallel distributions, and\\    
conformally symmetric manifolds}      
\author[A. Derdzinski]{Andrzej Derdzinski} 

\author[W. Roter]{Witold Roter}

%
\address{Department of Mathematics \hskip41pt 
Institute of Mathematics and Computer Science \endgraf
The Ohio State University \hskip49.2pt 
Wroc\l aw University of Technology \endgraf
Columbus, OH 43210 \hskip82.9pt 
Wy\-brze\-\.ze Wys\-pia\'n\-skiego 27, 50-370 Wroc\l aw \endgraf
USA \hskip165.8pt Poland \endgraf
{\eightit E-mail}\ {\eightit address}\hs{\eightrm:}\hskip4pt{\eighttt 
andrzej@math.ohio-state.edu} \hskip10.5pt 
{\eightit E-mail}\ {\eightit address}\hs{\eightrm:}\hskip4pt{\eighttt 
roter@im.pwr.wroc.pl} 
}

\begin{abstract}
We determine the local structure of all \prd s $\,(M,g)$ in dimensions 
$\,n\ge4\,$ whose Weyl conformal tensor $\,W\,$ is parallel and has rank 
$\,1\,$ when treated as an operator acting on exterior $\,2$-forms at each 
\pt. If one fixes three discrete parameters: the dimension $\,n\ge4$, the 
metric signature 
\hbox{$-$\hskip1pt$-$\hskip2pt$\dots$\hskip0pt$+$\hskip1pt$+$\hs,} and a sign 
factor $\,\ve=\pm\hs1\,$ accounting for sem\-i\-def\-i\-nite\-ness of 
$\,W\nh$, then the lo\-cal-isom\-e\-try types of our metrics $\,g\,$ 
correspond bijectively to equivalence classes of \su s $\,\bs\,$ with 
\ea\ \pftc s; the latter equivalence relation is provided by unimodular affine 
local diffeomorphisms. The \su\ $\,\bs\,$ arises, locally, as the leaf 
space of a co\-di\-men\-sion-two parallel \dn\ on $\,M$, naturally associated 
with $\,g$. We exhibit examples in which the leaves of the \dn\ form a 
fibration with the \ts\ $\,M\,$ and base $\,\bs$, for a closed \su\ $\,\bs\,$ 
of any prescribed diffeomorphic type.

Our result also completes a local classification of \prc s with parallel Weyl 
tensor that are neither \cf\ nor \ls: for those among such metrics which are 
not \rr, $\rwo$, and so they belong to the class mentioned above; on the other 
hand, the \rr\ ones have already been classified by the second author.
\end{abstract}

\maketitle

\voffset=-22pt\hoffset=7pt

\section*{Introduction}
The main result of the present paper, Theorem~\ref{class}, describes the local 
structure of \prc s whose Weyl conformal tensor $\,W\hs$ is parallel and, as 
an operator acting on exterior $\,2$-forms, has rank $\,1\,$ at each \pt. 
Combined with a theorem of the second author \cite{roter}, our description 
leads to a local classification of all \ecs\ \prd s. Here are some details.

A \prd\ $\,(M,g)\,$ with $\,\dim M\ge4\,$ is said to be {\it conformally 
symmetric\/} \cite{chaki-gupta} if its Weyl tensor $\,W\hs$ is parallel. 
Obvious examples arise when $\,(M,g)\,$ is \cf\ or \ls; \cs\ \mf s which are 
neither have usually been referred to as {\it essentially conformally 
symmetric\hh}. For a sample of recent results on \cs\ \mf s, see 
\cite{sharma}, \cite{rong}, \cite{deszcz}, \cite{deszcz-hotlos}.

An \ecs\ \prc\ is always indefinite \cite[Theorem~2]{derdzinski-roter-77}; if 
it is Lo\-rentz\-i\-an, it must also be \rr\ 
\cite[Corollary 1 on p.\ 21]{derdzinski-roter-78}. Known examples of \ecs\ 
indefinite metrics include both \rr\ \cite{roter} and non-\rr\ ones 
\cite{derdzinski-77}, in every dimension $\,n\ge4$.

Here we call a \prd\ {\it Ric\-ci-re\-cur\-rent\/} if, for every tangent \vf\ 
$\,v$, the Ricci tensor $\,\ri\,$ and the covariant derivative 
$\,\nabla_{\!v}\ri\,$ are linearly dependent at each \pt, or, equivalently, 
$\,\nabla\nnh\ri=\xi\otimes\ri\,$ on the open set where $\,\ri\ne0$, for some 
$\,1$-form $\,\xi$.

In \cite{roter} the second author completely determined the local structure of 
\rr\ \ecs\ metrics at \pt s where $\,\ri\,$ and $\,\nabla\nnh\ri\,$ are 
nonzero: in suitable coordinates, such metrics have the canonical form of 
\cite[formula~(34)]{roter}.

Our Theorem~\ref{class}, mentioned above, leads to a classification of the 
remaining \ecs\ \prc s. The reason is that the Weyl tensor has rank $\,1\,$ 
for every such metric which is not \rr\ 
\cite[Theorem 4 on p.\ 17]{derdzinski-roter-80}. Unlike the result of 
\cite{roter}, Theorem~\ref{class} requires no gen\-er\-al-po\-si\-tion 
hypothesis.

The text consists of three parts, each devoted to one of the three topics 
named in the title. Part I presents some known classification theorems about 
\pf\ connections on \su s, followed by results on solvability of a specific 
\pde\ (Sections \ref{adeq} -- \ref{prft}), both of which are needed as a 
reference for Part III.

Part II deals with \npd s $\,\vt\,$ on \prd s such that the \lcc\ $\,\nabla\,$ 
is, in a natural sense, $\,\w\nh$\pj\ onto a \tc\ $\,\hs\text{\rm D}\,$ on the 
(local) leaf space $\,\bs\,$ of $\,\w\nnh$. We observe that 
$\,\w\nh$-pro\-ject\-a\-bil\-i\-ty of $\,\hs\nabla\,$ is characterized by a 
simple curvature condition and, therefore, it holds for the $\,2$\diml\ \npd\ 
$\,\vt\,$ present on every \cs\ \mf\ with $\rwo$. This is one of the steps 
needed for the argument in Part III.

In Part III we establish our main result (Theorem~\ref{class}), first showing 
that, locally, in dimensions $\,n\ge4$, a \cs\ metric $\,g\,$ with $\rwo\,$ is 
a warped product in which the totally geodesic factor is four\diml, \cs\ and 
has $\rwo$, while the umbilical factor is flat. The problem is thus reduced to 
the case $\,n=4$. When $\,n=4$, we prove that $\,g\,$ is conformal to a metric 
of a kind first classified by Ruse \cite{ruse}; therefore, $\,g\,$ itself is 
one of Patterson and Walker's Riemann extension metrics 
\cite{patterson-walker}. We also establish, for every $\,n\ge4$,  projective 
flatness of the connection $\,\hs\text{\rm D}\,$ arising, as in the last 
paragraph, on the surface $\,\bs\,$ which is, locally, the leaf space of 
$\,\w\nnh$. The lo\-cal-isom\-e\-try type of $\,g\,$ then turns out to be 
parametrized by the dimension $\,n\hh$, the metric signature 
\hbox{$-$\hskip1pt$-$\hskip1pt$\dots$\hskip0pt$+$\hskip1pt$+$\hs,} a factor 
$\,\ve=\pm\hs1\,$ describing sem\-i\-def\-i\-nite\-ness of $\,W\nh$, and the 
(local) \ea\ equivalence class of $\,\hs\text{\rm D}\hh$. See 
Section \ref{comm}.

Finally, in Section \ref{cotq}, we describe examples showing that any 
prescribed closed \su\ $\,\bs\hs$ can be realized as the {\it global\/} leaf 
space of $\,\hs\w\nnh$ for some non-\rr\ \ecs\ \mf\ $\,(M,g)\,$ of any given 
dimension $\,n\ge4$. More precisely, the leaves of $\,\w$ then are the fibres 
of a \bd\ with the \ts\ $\,M\,$ and base $\,\bs$.

The authors wish to thank Zbigniew Olszak, Barbara Opozda and Udo Simon for 
valuable comments.

\setcounter{section}{0}
\section{Preliminaries}\label{prel}
Throughout this paper, all manifolds, \bd s, their sections and sub\bd s, 
as well as connections and mappings, including \bd\ morphisms, are assumed to 
be of class $\,C^\infty\nnh$. A manifold is by definition connected; a \bd\ 
morphism may operate only between two \bd s with the same \bmf, and acts by 
identity on the base.

`\nnh A \bd' always means `a $\,C^\infty$ locally trivial \bd' and the same 
symbol is used both for a given \bd\ and for its \ts. The \bd\ projection 
onto the \bmf\ is denoted by $\,\pmb$, and $\,\hs\kerd\pmb\,$ stands for \tvn.

By a {\it differential\/ $\,k$-form valued in a vector bundle\/} 
$\,\xz\,$ over a \mf\ $\,\bs\,$ we mean, as usual \cite[p.\ 24]{besse}, any 
\vbm\ $\,(\tb)^{\wedge k}\nnh\to\xz$. This includes the case of 
ordinary (real-val\-ued) forms, where $\,\xz\,$ is the product \bd\ 
$\,\bs\times\bbR\hs$, the sections of which are functions $\,\bs\to\bbR\hs$.

The symbols $\,\nabla\,$ and $\,\hs\text{\rm D}\hs\,$ will be used for various 
connections in \vb s. Our sign convention about the curvature tensor 
$\,R=R^\nabla$ of a connection $\,\nabla\,$ in a \vb\ $\,\xz\,$ over a 
\mf\ $\,\bs\,$ is
\begin{equation}
R(u,v)\psi\hskip7pt=\hskip7pt\nabla_{\!v}\nabla_{\!u}\psi\,
-\,\nabla_{\!u}\nabla_{\!v}\psi\,+\,\nabla_{[u,v]}\psi\hs,
\label{cur}
\end{equation}
for sections $\,\psi\,$ of $\,\xz\,$ and \vf s $\,u,v\,$ tangent to 
$\,\bs$. Thus,
\begin{equation}
R\,\,\mathrm{\ is\ a\ 2\hyp{}form\ on\ }\,\bs\,\mathrm{\ valued\ in\ }\,\hs
\text{\rm Hom}\hs(\xz,\xz)\hs.
\label{rtf}
\end{equation}
Here $\,\hs\text{\rm Hom}\hs(\mathcal{X},\mathcal{Y})$, for real \vb s 
$\,\mathcal{X},\mathcal{Y}\,$ over a manifold $\,\bs$, is the \vb\ over 
$\,\bs$, the sections of which are \vbm s $\,\mathcal{X}\to\mathcal{Y}$. For 
instance, 
$\,\mathcal{X}\hh^*\nh=\,\text{\rm Hom}\hs(\mathcal{X},\bs\times\bbR)\,$ is 
the dual of $\,\mathcal{X}$. By (\ref{cur}), for connections in the \tbd,
\begin{equation}
R_{jkl}{}^m\,=\,\partial_k\vg_{\!jl}^m\,-\,\,\partial_j\vg_{\!kl}^m
\,+\,\,\vg_{\!ks}^m\vg_{\!jl}^s\,
-\,\,\vg_{\!js}^m\vg_{\!kl}^s\hs,
\label{rlm}
\end{equation}
$\partial_j$ and $\,\vg_{\!jk}^l$ being the partial derivatives and the 
connection components.

The \rt\ $\,\ri=\ri^{\hs\text{\rm D}}$ of a connection 
$\,\hs\text{\rm D}\hs\,$ on a \mf\ $\,\bs\,$ (that is, in the \tbd\ $\,\tb$) 
is given by $\,\ri\hh(u,w)=\,\text{\rm Trace}\,[v\mapsto R(u,v)w\hh]$, where 
$\,R=R^{\hs\text{\rm D}}$ and $\,u,v,w$ are vectors tangent to $\,\bs\,$ at 
any \pt. It is well known that, if $\,\hs\text{\rm D}\,$ is 
tor\-sion\-free and $\,\dim\bs=\er$, the following four conditions are 
equivalent:
\begin{enumerate}
  \def\theenumi{{\rm\alph{enumi}}}
\item the connection induced by $\hs\text{\rm D}\hs$ in the line \bd\ 
$(\tb)^{\wedge\er}\nh$ is flat,
\item a nonzero $\,\hs\text{\rm D}\hs$\prl\ differential $\,r$-form exists on 
every \sco\ open subset of $\,\bs$, 
\item the operator $\,R^{\hs\text{\rm D}}(u,v):\tcb\to\tcb\,$ given by 
$\,w\mapsto R^{\hs\text{\rm D}}(u,v)w\,$ is traceless for every 
$\,\cy\in\bs\,$ and all $\,u,v\in\tcb$,
\item $\ri^{\hs\text{\rm D}}$ is symmetric.
\end{enumerate}
In fact, the curvature tensor of the connection that $\,\hs\text{\rm D}\hs\,$ 
induces in $\,(\tb)^{\wedge\er}$, viewed as a real-val\-ued $\,2$-form (by 
(\ref{rtf}) with the identification 
$\,\hs\text{\rm Hom}\hs((\tb)^{\wedge\er}\nnh,(\tb)^{\wedge\er})\nh
=\nh\bs\times\bbR$), sends $\,u,v\in\tcb\,$ to 
$\,\hs\text{\rm Trace}\,[R^{\hs\text{\rm D}}(u,v)]$, while, by the first 
Bianchi identity, $\,\hs\text{\rm Trace}\,[R^{\hs\text{\rm D}}(u,v)]
=\ri^{\hs\text{\rm D}}\nh(u,v)-\ri^{\hs\text{\rm D}}\nh(v,u)$. Thus, (a) is 
equivalent both to (c) and to (d). Finally, the connections that 
$\,\hs\text{\rm D}\hs\,$ induces in $\,(\tb)^{\wedge\er}$ and 
$\,(\tab)^{\wedge\er}=[(\tb)^{\wedge\er}]{}^*$ are each other's duals, so that 
(a) holds if and only if (b) does.

A fixed connection $\,\nabla\,$ in a \vb\ 
$\,\xz\,$ over a \mf\ $\,\bs$ gives rise to the operator of {\it 
exterior derivative\/} $\,d^{\hs\nabla}$ acting on $\,\xz$-val\-ued 
differential forms on $\,\bs\,$ which, for the standard flat connection in the 
product \bd\ $\,\bs\times\bbR\hs$, is the ordinary exterior derivative $\,d\,$ 
for real-val\-ued forms \cite[p.\ 24]{besse}. Explicitly, for an 
$\,\xz$-val\-ued $\,1$-form $\,\varPsi$,
\begin{equation}
(d^{\hs\nabla}\nnh\varPsi)(u,v)\,=\,\nabla_{\!u}[\varPsi(v)]\,
-\,\nabla_{\!v}[\varPsi(u)]\,-\,\varPsi([u,v])\hs.
\label{dnp}
\end{equation}
We will establish the local solvability of various systems of linear \pde s 
(homogeneous or not) on a \sco\ \mf\ $\,\bs\,$ by introducing a flat 
connection $\,\nabla\,$ in a suitable \vb\ $\,\xz\,$ over $\,\bs\,$ and 
showing that
\begin{enumerate}
  \def\theenumi{{\rm\roman{enumi}}}
\item our (homogeneous) system amounts to requiring an unknown section of 
$\,\xz\,$ to be $\,\nabla$\prl, or,
\item the (nonhomogeneous) system can be rewritten so as to be imposed on 
$\,\xz$-val\-ued differential forms, and then its solvability is equivalent to 
$\,d^{\hs\nabla}\hskip-2.4pt$-exact\-ness of the right-hand side (which, 
locally, is the same as its easily-verified 
$\,d^{\hs\nabla}\hskip-2.4pt$-closed\-ness).
\end{enumerate}
Note that the Poincar\'e lemma can be applied, in (ii), since flatness of 
$\,\nabla\,$ and simple connectivity of $\,\bs\,$ allow us to treat 
$\,\xz$-val\-ued differential forms as forms taking values in a fixed 
fi\-nite\diml\ \vs, in such a way that $\,d^{\hs\nabla}$ becomes the ordinary 
$\,d\hh$.
\begin{remark}\label{twice}Twice-co\-var\-i\-ant \tf s $\,\tau\,$ on a \mf\ 
$\,\bs\,$ will also be treated as $\,\tab$-val\-ued $\,1$-forms. Namely, we 
choose $\,\tau(w,\,\cdot\,)$ to be the section of $\,\tab\,$ to which the 
$\,1$-form $\,\tau\,$ sends a \vf\ $\,w\,$ on $\,\bs$. Here $\,\tau\,$ is not 
assumed symmetric; the use of $\,\tau(\,\cdot\,,w)$ instead of 
$\,\tau(w,\,\cdot\,)\,$ would amount to replacing $\,\tau\,$ by its {\it 
transpose\/} $\,\tau^*\nnh$, which is the twice-co\-var\-i\-ant \tf\ assigning 
$\,\tau(w\hs'\nnh,w)\,$ to \vf s $\,w,w\hs'\nnh$. For instance, if 
$\,\tau=\hs\text{\rm D}\hh\xi\,$ for a connection $\,\hs\text{\rm D}\hs\,$ and 
a $\,1$-form $\,\xi\,$ on $\,\bs$, we have 
$\,\tau(w,\,\cdot\,)=\hs\text{\rm D}_w\xi$, and, in local coordinates, 
the components of $\,\tau=\hs\text{\rm D}\hh\xi$ and 
$\,\gy=(\hh\text{\rm D}\hs\xi)^*$ are $\,\tau_{jk}=\xi_{\hh k,j}$ and 
$\,\gy_{jk}=\xi_{\hh j,k}$.
\end{remark}
\begin{remark}\label{codaz}Let $\,\hs\text{\rm D}\hs\,$ be a \tc\ on a \mf\ 
$\,\bs$. A twice-co\-var\-i\-ant \sy\ \tf\ $\,\tau\,$ on $\,\bs\,$ is said to 
satisfy the {\it Codazzi equation\/} if $\,d^{\hs\text{\rm D}}\nh\tau=\nh0\,$ 
for $\,\tau\,$ treated as a $\,\tab$-val\-ued $\,1$-form (Remark~\ref{twice}), 
which in coordinates reads $\tau_{jl,k}\hskip-2.5pt=\nh\tau_{kl,j}$, cf.\ 
(\ref{dnp}). Suppose now that, in addition, the \rt\ $\,\ri^{\hs\text{\rm D}}$ 
is \sy. In view of the second Bianchi identity, the curvature tensor 
$\,R^{\hs\text{\rm D}}$ has zero divergence 
($\text{\rm div}{}^{\hs\text{\rm D}}\nnh R^{\hs\text{\rm D}}\hskip-2.5pt
=0\hh$) if and only if the Codazzi equation 
$\,d^{\hs\text{\rm D}}\nh\ri^{\hs\text{\rm D}}\hskip-2.5pt=\nh0$ holds for 
$\,\ri^{\hs\text{\rm D}}\nnh$. (In coordinates, this means that the 
condition $\,R_{jkl}{}^s{}_{\nh,\hh s}=0\,$ is equivalent to 
$\,R_{jl,k}\hskip-2.5pt=\nh R_{kl,j}$.) Riemannian \mf s with 
$\,\hs\text{\rm div}{}^\nabla\hskip-1.8ptR^\nabla\nnh=0\,$ for the \lcc\ 
$\,\nabla\,$ are said to have {\it harmonic curvature\/} \cite{besse}.
\end{remark}
\begin{remark}\label{onetw}We always treat $\,2$-forms on a \mf\ $\,\bs\,$ 
valued in ordinary $\,1$-forms (that is, $\,\tab$-val\-ued) as $\,1$-forms on 
$\,\bs\,$ valued in ordinary $\,2$-forms (i.e., 
$\,(\tab)^{\wedge2}\nnh$-val\-ued), via the obvious identification.
\end{remark}
\begin{remark}\label{ddtau}Aside from differentials of functions, we will 
apply the ex\-te\-ri\-or de\-riv\-a\-tive operators 
$\,d\hh,\hs d^{\hs\text{\rm D}}$ and $\,d^{\hs\nabla}\nnh$ {\it only to 
$1$-forms\/} valued in various \vb s. For instance, given a 
twice-co\-var\-i\-ant \tf\ $\,\tau\,$ and a tor\-sion\-free connection 
$\,\hs\text{\rm D}\hs\,$ on a \mf\ $\,\bs$, the exterior derivative 
$\,d^{\hs\text{\rm D}}\nh\tau\,$ (of $\,\tau\,$ viewed as a $\,\tab$-val\-ued 
$\,1$-form, cf.\ Remark~\ref{twice}), is itself treated, according to 
Remark~\ref{onetw}, as a $\,1$-form valued in $\,2$-forms. If, in addition, 
$\,\dim\bs=2\,$ and $\,\ri^{\hs\text{\rm D}}$ is \sy, we use this last 
interpretation to define the exterior derivative 
$\,dd^{\hs\text{\rm D}}\nh\tau$, writing $\,d\,$ instead of 
$\,d^{\hs\text{\rm D}}$ to reflect the fact that the existence, locally in 
$\,\bs$, of a $\,\hs\text{\rm D}\hs$\prl\ area form (see (a) -- (d) in 
Section \ref{prel}), allows us to regard $\,(\tab)^{\wedge2}\nnh$-val\-ued 
forms, locally, as real-val\-ued.

The $\,(\tab)^{\wedge2}\nnh$-val\-ued $\,2$-form 
$\,A=dd^{\hs\text{\rm D}}\nh\tau\,$ assigns to \vf s $\,u,v\,$ the 
real-val\-ued $\,2$-form $\,A(u,v)\,$ with the lo\-cal-co\-or\-di\-nate 
expression $\,[A(u,v)]_{lm}=u^jv^k\hh A_{jklm}$, where $\,A_{jklm}
=\tau_{mk,lj}\nh-\nh\tau_{lk,mj}\nh-\nh\tau_{mj,lk}\nh+\nh\tau_{lj,mk}$.
\end{remark}

\section{Exterior products}\label{expr}
If $\,\hs\text{\rm D}\hs\,$ is a \tc\ on a \su\ and the \rt\ 
$\,\ri^{\hs\text{\rm D}}$ is \sy, the curvature tensor $\,R^{\hs\text{\rm D}}$ 
is given by
\begin{equation}
R^{\hs\text{\rm D}}=\,\hs\ri^{\hs\text{\rm D}}\nnh\wedge\hs\text{\rm Id}\hs,
\hskip7pt\mathrm{that\ is,\ in\ coordinates,}\hskip6pt
R_{jkl}{}^m\hs
=\,R_{jl}\delta_k^{\hs m}-\hs R_{kl}\delta_j^{\hs m}\nnh.\hskip-5pt
\label{rri}
\end{equation}
Namely, for $\,R=R^{\hs\text{\rm D}}$ and $\,\ri=\ri^{\hs\text{\rm D}}$ 
treated as algebraic objects, at any fixed \pt, $\,\ri\,$ is the Ricci 
contraction of $\,R=\ri\wedge\text{\rm Id}\hs$, and so the operator 
$\,\ri\mapsto R=\ri\wedge\text{\rm Id}\hs$, being injective, is also 
surjective. (Both $\,R,\ri\,$ range over $\,3$\diml\ spaces, cf.\ (c) in 
Section \ref{prel}.)

Aside from formula (\ref{rri}), we will use the symbol $\,\wedge\,$ only in 
two situations, which involve a real-val\-ued $\,1$-form $\,\xi\,$ and 
twice-co\-var\-i\-ant \tf s $\,\tau,\gy\,$ on any given \mf\ $\,\bs$. 
Namely, $\,\xi\wedge\gy\,$ stands for the $\,1$-form valued in $\,2$-forms, 
sending a \vf\ $\,w\,$ to the $\,2$-form that assigns to \vf s $\,u,v\,$ the 
function $\,\xi(u)\hh\gy(v,w)-\xi(v)\hh\gy(u,w)$, while $\,\tau\wedge\gy\,$ is 
a $\,2$-form valued in $\,2$-forms, and it sends \vf s $\,w,w\hs'$ to the 
$\,2$-form $\,\tau(w,\,\cdot\,)\wedge\gy(w\hs'\nnh,\,\cdot\,)-
\tau(w\hs'\nnh,\,\cdot\,)\wedge\gy(w,\,\cdot\,)\,$ (which in turn associates 
with $\,u,v\,$ the function $\,\tau(w,u)\hh\gy(w\hs'\nnh,v)
-\tau(w,v)\hh\gy(w\hs'\nnh,u)-\tau(w\hs'\nnh,u)\hh\gy(w,v)
+\tau(w\hs'\nnh,v)\hh\gy(w,u)$). Note that 
$\,\xi,\tau\,$ and $\,\gy\,$ may be viewed as differential forms on $\,\bs\,$ 
valued in differential forms (Remark~\ref{twice}), and then $\,\wedge\,$ 
becomes the usual exterior product (with the identification described in 
Remark~\ref{onetw}). For twice-co\-var\-i\-ant \tf s $\,\tau,\gy\,$ on a {\it 
surface\/} $\,\bs$,
\begin{equation}
\gy\wedge\tau\hs=\,dd^{\hs\text{\rm D}}\nh\gy\,=\,0\hskip7pt\mathrm{\ if\ }
\,\tau\,\mathrm{\ is\ symmetric\ and\ }\,\gy\,\mathrm{\ is\ 
skew\hyp{}symmetric.}
\label{tgz}
\end{equation}
In fact, in local coordinates $\,A=\tau\wedge\gy\,$ is given by 
$\,A_{jklm}\nh=\nh\tau_{jl}\gy_{km}\nh-\nh\tau_{jm}\gy_{kl}\nh
-\nh\tau_{kl}\gy_{jm}\nh+\nh\tau_{km}\gy_{jl}$, with the only essential 
component $\,A_{1212}\nh=0$. On the other hand, locally, $\,\gy=f\alpha$, 
where $\,\alpha\,$ is a fixed $\,\hs\text{\rm D}\hs$\prl\ area form, and so 
$\,d^{\hs\text{\rm D}}\nh\gy=\df\nnh\wedge\alpha=-\hs\df\nnh\otimes\alpha\,$ 
(equality of $\,1$-forms valued in $\,2$-forms; the sum 
$\,\df\nnh\otimes\alpha+\df\nnh\wedge\alpha\,$ must vanish, as it is a 
real-val\-ued differential $\,3$-form, while $\,\dim\bs=2$). Thus, 
$\,dd^{\hs\text{\rm D}}\nh\gy=-\hs d\df\nnh\otimes\alpha=0$.

For a \tc\ $\,\hs\text{\rm D}\hs\,$ and a $\,1$-form $\,\xi\,$ on a \mf\ 
$\,\bs$,
\begin{equation}
\text{\rm(i)}\hskip9ptd^{\hs\text{\rm D}}\text{\rm D}\hs\xi\hs
=\hs\xi\hs R^{\text{\rm D}}\nh,\hskip44pt\text{\rm(ii)}\hskip9pt
d^{\hs\text{\rm D}}\nnh(\hh\text{\rm D}\hs\xi)^*\nh
=\hs\text{\rm D}\hs d\hh\xi\hs+\hs\xi\hs R^{\text{\rm D}}\nh.
\label{rid}
\end{equation}
Here (i) is the {\it Ricci identity\/} (in coordinates, 
$\,\xi_{\hh j,kl}-\xi_{\hh j,lk}=R_{lkj}{}^s\hs\xi_s$), 
$\,R^{\text{\rm D}}$ denotes, as usual, the curvature tensor of 
$\,\hs\text{\rm D}\hh$, and $\,\xi\hs R^{\text{\rm D}}$ is the $\,2$-form 
valued in $\,1$-forms that assigns to \vf s $\,w,w\hs'$ on $\,\bs\,$ the 
the composite \vbm\ $\,\tb\to\tb\to\bs\times\bbR\,$ in which 
$\,R^{\text{\rm D}}(w,w\hs'\hh)\,$ is followed by $\,\xi\,$ (cf.\ 
(\ref{rtf})). To establish (ii), note that the coordinate version of (i) gives 
$\,\xi_{\hh k,lj}-\xi_{\hh j,lk}=(\xi_{\hh k,j}-\xi_{\hh j,k}){}_{,\hh l}
+R_{jlk}{}^s\hs\xi_s-R_{klj}{}^s\hs\xi_s$, the last two terms of which may be 
replaced by $\,R_{jkl}{}^s\hs\xi_s$ in view of the first Bianchi identity.

If, in addition, $\,\dim\bs=2\,$ and the \rt\ $\,\ri^{\hs\text{\rm D}}$ is \sy,
\begin{equation}
\text{\rm(i)}\hskip9pt\xi\hs R^{\text{\rm D}}\nh
=-\hs\xi\wedge\ri^{\hs\text{\rm D}}\nh,\hskip48pt
\text{\rm(ii)}\hskip9ptd^{\hs\text{\rm D}}\text{\rm D}\hs\xi\hs
=-\hs\xi\wedge\ri^{\hs\text{\rm D}}\nh,
\label{ris}
\end{equation}
by (\ref{rri}) and (\ref{rid}\hh-i). The coordinate form of (\ref{ris}\hh-ii) 
is $\,\,\xi_{\hh j,kl}-\xi_{\hh j,lk}
=\xi_{\hh k}\hs R_{lj}-\xi_{\hh l}\hs R_{kj}$.

The Ricci identity for twice-co\-var\-i\-ant \tf s $\,\gy\,$ (sections of 
$\,(\tam)^{\otimes2}$) reads
\begin{equation}
d^{\hs\nabla}\nh\nabla\gy\hs=\hs\gy\hs\cdot R\hs,
\hskip7pt\mathrm{that\ is,}\hskip6pt\gy_{\hh lm,kj}-\gy_{\hh lm,jk}
=R_{jkl}{}^s\hs\gy_{sm}+R_{jkm}{}^s\hs\gy_{\hh ls}\hs.
\label{rit}
\end{equation}
Here $\,M\,$ is any \mf\ with a fixed \tc\ (this time denoted by $\,\nabla$) 
and $\,R=R^\nabla\nnh$, while the covariant derivative $\,\nabla\gy$, treated 
as a $\,(\tam)^{\otimes2}\nnh$-val\-ued $\,1$-form on $\,M$, associates with a 
\vf\ $\,w\,$ the section $\,\nabla_{\!w}\gy\,$ of $\,(\tam)^{\otimes2}\nnh$. 
Finally, $\,\gy\hs\cdot R\,$ denotes the $\,(\tam)^{\otimes2}\nnh$-val\-ued 
$\,2$-form, sending \vf s $\,w,w\hs'$ to the section of $\,(\tam)^{\otimes2}$ 
which assigns to \vf s $\,u,v\,$ the function 
$\,\gy(R(w,w\hs'\hh)u,v)+\gy(u,R(w,w\hs'\hh)v)$.

Given a twice-co\-var\-i\-ant \tf\ $\,\tau\,$ and a \vf\ $\,u\,$ on a \prd\ 
$\,(M,g)$, we define $\,\tau u\,$ to be the \vf\ on $\,M\,$ which is the image 
of $\,u\,$ under the \vbm\ $\,\tm\to\tm\,$ obtained from $\,\tau\,$ by index 
raising; thus,
\begin{equation}\label{tau}
\text{\rm$g(\tau u,v)=\tau(u,v)\,$ for\ all\ vector\ fields\ $\,u\,$\ and\ 
$\,v$.}
\end{equation}
In a \prd\ $\,(M,g)\,$ of any dimension $\,n$, we denote by $\,\nabla\,$ its 
\lcc, and by $\,R,\ri\,$ (rather than $\,R^\nabla\nnh,\ri^\nabla$) its 
curvature and \rt s. The same symbol $\,R\,$ is used for the four-times 
covariant curvature tensor (a $\,2$-form valued in $\,2$-forms), with 
$\,R(u,v,w,w'\hs)=g(R(u,v)w,w'\hs)$. If $\,n\ge4$, the Weyl conformal tensor 
of $\,(M,g)\,$ is defined by $\,W\nnh=R-(n-2)^{-1}\hs g\wedge\hh\sj$, where 
$\,\sj=\ri-\hs(2n-2)^{-1}\,\text{\rm s}\hs g\,$ is the {\it Schouten 
tensor}, with $\,\hs\text{\rm s}\hs=\hs\text{\rm Trace}_{\hh g}\hh\ri$ 
standing for the scalar curvature, and $\,\wedge\,$ as above.

For \vf s $\,u,v$, a differential $\,2$-form $\,\om\hh$, the \lcc s 
$\,\nabla,\widetilde\nabla\,$ and \rt s $\,\ri,\widetilde\ri\,$ of conformally 
related metrics $\,g\,$ and $\,\widetilde g=\fy^{-2}g=e^{2\phi}g\,$ on a 
manifold $\,M$, with functions $\,\fy>0\,$ and $\,\phi=-\hs\log\fy$, their 
$\,g$-gra\-di\-ents $\,\nabla\nh\phi,\naf$, and the $\,g$-La\-plac\-i\-an 
$\,\Delta f=g^{jk}\fy_{,\hh jk}$ of $\,\fy$, we have (cf.\ 
\cite[formulae~(16.8), (16.9), (16.13) on pp.\ 528--529]{dillen-verstraelen})
\begin{enumerate}
  \def\theenumi{{\rm\alph{enumi}}}
\item $\widetilde\nabla_{\!u}v\,=\,\nabla_{\!u}v\,+\,g(u,\nabla\nh\phi)\hs v\,
+\,g(v,\nabla\nh\phi)\hs u\,-\,g(u,v)\nabla\nh\phi$\hs,
\item $\widetilde\nabla_{\!u}\hs\om\,=\,\nabla_{\!u}\hs\om\,
-\,2(d_u\phi)\hs\om+\om(u,\,\cdot\,)\wedge\hs d\phi\,
+\,g(u,\,\cdot\,)\wedge\hs\om\hh(\nabla\nh\phi,\,\cdot\,)$\hs,
\item $\widetilde\ri\,=\hs\ri\,+(n-2)\fy^{-1}\nabla\df
+\left[\fy^{-1}\Delta\fy-(n-1)\fy^{-2}g(\naf,\naf)\right]g$, where 
$\,n=\dim M$,
\end{enumerate}
$d_u$ in (b) being the directional derivative, so that 
$\,d_u\phi=g(u,\nabla\nh\phi)$.

\newpage
\part{PROJECTIVELY FLAT SURFACES}\label{pfs}
Except for Sections~\ref{adeq} -- \ref{prft}, the material in Part I is known, 
and consists of classification results about \pf\ equi\-af\-fine \tc s on 
\su s. A self-con\-tain\-ed presentation of those results is provided here for 
the reader's convenience: such a connection serves as the single 
non-dis\-crete parameter in our classification of \cs\ \mf s with $\rwo\,$ 
(see Section \ref{tlst}). For more on \pf\ \su s, 
see \cite{pinkall-schwenk-schellschmidt-simon} and Simon's article on affine 
differential geometry \cite[pp.\ 905--961]{dillen-verstraelen}.

\section{Projective flatness in dimension $\,2$}\label{pfid}
A connection $\,\hs\text{\rm D}\hh\,$ on a \mf\ $\,\bs\,$ is called {\it 
projectively flat\/} if $\,\bs$ is covered by coordinate systems in which the 
geodesics of $\,\hs\text{\rm D}\,$ appear as (re-pa\-ram\-e\-trized) 
straight-line segments.

We begin with a well-known lemma, going back to Weyl \cite[p.~100]{weyl}:
\begin{lemma}\label{smgeo}Two \tc s\/ $\,\hs\text{\rm D}\hs$ and\/ 
$\,\hs\deh\hs$ on a \mf\/ $\,\bs\,$ have the same 
re-pa\-ram\-e\-trized geodesics if and only if\/ 
$\,B=\hs2\hs\xi\odot\nh\text{\rm Id}\hs\,$ for the \tf\/ 
$\,B=\hs\deh\hs\hs-\hs\text{\rm D}$ and a real-val\-ued\/ $\,1$-form\/ 
$\,\xi\,$ on $\,\bs$, or, in coordinates,  
$\,B_{jk}^l=\xi_j\delta_k^l\nh+\xi_k\delta_j^l$. Then also
\begin{enumerate}
  \def\theenumi{{\rm\roman{enumi}}}
\item the \rt s $\,\ri^{\hs\text{\rm D}}$ of\/ $\,\hs\text{\rm D}\hh\,$ and\/ 
$\,\rih\,$ of\/ $\,\hs\deh\hh\,$ are related by 
$\,\rih=\ri^{\hs\text{\rm D}}\nnh+(\er-1)\hskip1pt\xi\otimes\xi
-\hs\er\hs\text{\rm D}\hs\xi+(\hh\text{\rm D}\hs\xi)^*$ for $\,\er=\dim\bs$, 
with $\,(\hh\text{\rm D}\hs\xi)^*$ as in Remark\/ $\ref{twice}$,
\item $d\hs\xi=0\,$ if both\/ $\,\hs\text{\rm D}\hh\,$ and\/ $\,\hs\deh\hs$ 
have \sy\ \rt s.
\end{enumerate}
\end{lemma}
\begin{proof}In view of the lo\-cal-co\-or\-di\-nate form of the geodesic 
equation, the two connections have the same geodesics if and only if 
$\,B_vv=c\hs v\,$ for every $\,\cy\in\bs\,$ and $\,v\in\tcb$, with some 
$\,c\in\bbR\,$ depending on $\,v$. The latter condition gives 
$\,B=\hs2\hs\xi\odot\nh\text{\rm Id}\hs\,$ (by polarization in a fixed basis 
of $\,\tcb$, for any $\,\cy\in\bs$). Next, if 
$\,\hs\deh\hs-\hs\text{\rm D}\,
=\hs2\hs\xi\odot\nh\text{\rm Id}\hs$, the curvature tensors 
$\,R^{\hs\text{\rm D}}$ of $\,\hs\text{\rm D}\hs\,$ and $\,\widetilde R\,$ of 
$\,\hs\deh\hs\,$ are related by 
$\,\widetilde R=R^{\text{\rm D}}\nh-\hs d\hs\xi\otimes\text{\rm Id}\,-
\hs\text{\rm D}\hs\xi\wedge\nh\text{\rm Id}\,
+\hs\xi\otimes(\xi\wedge\nh\text{\rm Id}\hh)\,$ (that is, 
$\,\widetilde R_{jkl}{}^m=R_{jkl}{}^m+(\xi_{\hh j,k}-\xi_{\hh k,j})\delta_l^m
+\xi_{\hh l,k}\delta_j^m-\xi_{\hh l,j}\delta_k^m
+\xi_{\hh l}(\xi_{\hh j}\delta_k^m-\xi_{\hh k}\delta_j^m)$), as one easily 
verifies using (\ref{rlm}) in coordinates in which the components of 
$\,\hs\text{\rm D}\hs\,$ at the given \pt\ vanish. This implies (i). Now (ii) 
follows: if $\,\ri^{\hs\text{\rm D}}$ and $\,\widetilde\ri\,$ are \sy, so 
must be $\,\hs\text{\rm D}\hs\xi$.
\end{proof}
\begin{remark}\label{ddfrd}For a \tc\ $\,\hs\text{\rm D}\hs\,$ on a \su\ 
$\,\bs\,$ such that the \rt\ $\,\ri^{\hs\text{\rm D}}$ is \sy, a function 
$\,\fy:\bs\to(0,\infty)$, and the $\,1$-form $\,\xi=-\hs d\hh\log\fy$, the 
condition $\,\hs\text{\rm D}\hh\df=-\fy\nnh\ri^{\hs\text{\rm D}}$ is 
equivalent to flatness of the connection 
$\,\hs\deh\hs=\hs\text{\rm D}\hs+2\hs\xi\odot\nh\text{\rm Id}\hs$. In fact, by 
Lemma~\ref{smgeo}, $\,\rih=0\,$ if and only if 
$\,\hs\text{\rm D}\hh\df=-\fy\nnh\ri^{\hs\text{\rm D}}\nnh$, while $\,\rih\,$ 
determines $\,\widetilde R$, as in (\ref{rri}).
\end{remark}
The next result is the $\,2$\diml\ case of a theorem of Weyl 
\cite[p.~105]{weyl}:
\begin{theorem}\label{prjfl}A \tc\/ $\,\hs\text{\rm D}\hs\,$ on a \su\ 
$\,\bs\,$ with a \sy\ \rt\ $\,\ri^{\hs\text{\rm D}}\nh$ is \pf\ if and only 
if\/ $\,\hs\ri^{\hs\text{\rm D}}\nh$ satisfies the Codazzi equation 
$\,d^{\hs\text{\rm D}}\hskip-2.2pt\ri^{\hs\text{\rm D}}\hskip-2.5pt=0$.
\end{theorem}
\begin{proof}In the \vb\ $\,\tab\oplus(\bs\nh\times\bbR)$, whose sections 
$\,(\xi,\fy)\,$ consist of a $\,1$-form $\,\xi\,$ and a function $\,\fy\,$ on 
$\,\bs$, we define a connection $\,\nabla\,$ by 
$\,\nabla(\xi,\fy)=(\hh\text{\rm D}\hh\xi
+\fy\nnh\ri^{\hs\text{\rm D}}\nnh,\hs\df-\xi)$, that is, 
$\,\nabla_{\!u}(\xi,\fy)=(\hh\text{\rm D}_u\hh\xi
+\fy\nnh\ri^{\hs\text{\rm D}}\nh(u,\,\cdot\,),\hs d_uf-\xi(u))\,$ for any \vf\ 
$\,u$. From (\ref{cur}) with $\,\psi=(\xi,\fy)\,$ we get 
$\,R^\nabla\nnh(u,v)\psi=(\xi\hh'(u,v),0)$, where 
$\,\xi\hh'=-\,d^{\hs\text{\rm D}}\text{\rm D}\hs\xi
-\xi\wedge\ri^{\hs\text{\rm D}}\nh
-\fy d^{\hs\text{\rm D}}\hskip-2.2pt\ri^{\hs\text{\rm D}}\nnh$, with 
$\,d^{\hs\text{\rm D}}\text{\rm D}\hs\xi\,$ as in (\ref{rid}\hh-i). Thus, by 
(\ref{ris}\hh-ii), 
$\,d^{\hs\text{\rm D}}\hskip-2.2pt\ri^{\hs\text{\rm D}}\hskip-2.5pt=\nh0\,$ if 
and only if $\,\nabla\,$ is flat.

If $\,\nabla\hs$ is flat, we may choose, locally, a $\,\nabla$\prl\ section 
$\,(\xi,\fy)\,$ of $\,\tab\oplus(\bs\nh\times\bbR)$ with $\,\fy>0$. Now 
$\,\hs\text{\rm D}\hs\,$ is \pf\ by Remark~\ref{ddfrd} and Lemma~\ref{smgeo}, 
as $\,\hs\text{\rm D}\hh\df=-\fy\nnh\ri^{\hs\text{\rm D}}\nnh$.

Conversely, let $\,\hs\text{\rm D}\hs\,$ be \pf. By Lemma~\ref{smgeo}, we can 
find, locally, a positive function $\,\fy\,$ such that the \tc\ 
$\,\hs\deh\hs=\hs\text{\rm D}\hs+2\hs\xi\odot\nh\text{\rm Id}\hs$, with 
$\,\xi=-\hs d\hh\log\fy$, is flat. Remark~\ref{ddfrd} then gives 
$\,\hs\text{\rm D}\hh\df=-\fy\nnh\ri^{\hs\text{\rm D}}$ and, applying 
$\,d^{\hs\text{\rm D}}$ to both sides (treated as $\,\tab$-val\-ued 
$\,1$-forms) we get $\,d^{\hs\text{\rm D}}\text{\rm D}\hh\df
=-\hs\df\wedge\ri^{\hs\text{\rm D}}\nh
-\fy d^{\hs\text{\rm D}}\hskip-2.2pt\ri^{\hs\text{\rm D}}\nnh$, while 
$\,d^{\hs\text{\rm D}}\text{\rm D}\hh\df=-\hs\df\wedge\ri^{\hs\text{\rm D}}$ 
in view of (\ref{ris}\hh-ii) with $\,\xi\,$ replaced by $\,\df$. Thus, 
$\,d^{\hs\text{\rm D}}\hskip-2.2pt\ri^{\hs\text{\rm D}}\hskip-2.5pt=0$.
\end{proof}
Suppose that $\,\bs\,$ is an $\,\er$\diml\ \mf, $\,\hy\,$ is a \rvs\ of 
dimension $\,\er+1$, and $\,\im:\bs\to\hy\,$ is an immersion transverse to all 
lines in $\,\hy\hs$ passing through $\,0$. By the {\it cen\-tro\-af\-fine 
connection\/} corresponding to $\,\im\,$ we mean the connection 
$\,\hs\text{\rm D}\hs\,$ on $\,\bs\,$ obtained by first pulling back the 
standard flat connection of $\,\hy\hs$ to the pull\-back \bd\ 
$\,\im^*T\hy=\bs\times\hy\nnh$, and then projecting it onto the summand 
$\,\tb\,$ in $\,\im^*T\hy=\tb\oplus\mathcal{N}\nh$. Here $\,\mathcal{N}\,$ is 
the {\it normal bundle\/} over $\,\bs\,$ whose fibre $\,\mathcal{N}_\cy$, for 
any $\,\cy\in\bs$, is the line in $\,\hy=T_{\im(\cy)}\hy$ spanned by 
$\,\im(\cy)$. See \cite[pp.\ 927,\ 934]{dillen-verstraelen}.
\begin{remark}\label{ctraf}Let $\,\hs\text{\rm D}\hs\,$ be the 
cen\-tro\-af\-fine connection on $\,\bs$, for $\,\bs,\hy,\er\,$ and $\,\im\,$ 
as above.
\begin{enumerate}
  \def\theenumi{{\rm\alph{enumi}}}
\item We may assume, locally, that $\,\bs\,$ is a sub\mf\ of $\,\hy\hs$ and 
$\,\im\,$ is the inclusion mapping. The geodesics of $\,\hs\text{\rm D}\hs$ 
then have the form $\,\bs\cap\varPi$, for planes $\,\varPi\,$ through $\,0\,$ 
in $\,\hy$, since for any parametrization of such a curve $\,\bs\cap\varPi$, 
the acceleration lies in $\,\varPi$, and so its component tangent to $\,\bs\,$ 
is tangent to $\,\bs\cap\varPi$. Cf.\ \cite[p.\ 927]{dillen-verstraelen}.
\item The connection $\,\hs\text{\rm D}\,$ is tor\-sion\-free and \pf.
\item The \rt\ $\,\ri^{\hs\text{\rm D}}$ is \sy\ and 
$\,\ri^{\hs\text{\rm D}}\nh=(1-\er)\hs\bz$, where $\,\bz\,$ is the second 
fundamental form of the immersion $\,\im$, obtained in the usual fashion from 
the trivialization of the normal \bd\ provided by the radial (identity) \vf.
\end{enumerate}
In fact, projective flatness of $\,\hs\text{\rm D}\,$ is immediate from 
(a): central projections of $\,\bs\,$ into hyperplanes not containing $\,0\,$ 
send geodesics of $\,\hs\text{\rm D}\,$ into lines. To verify the other 
claims in (b) -- (c) we may use the traditional notation in which the 
inclusion mapping $\,\bs\to\hy$, here serving also as a trivializing section 
of the normal \bd, is represented by the symbol $\,\mathbf{r}\hs$ and treated 
as a $\,\hy$-val\-ued function on $\,\bs$, while its partial 
differentiations of all orders, relative to fixed coordinates in $\,\bs$, are 
represented by successive subscripts; for instance, $\,\mathbf{r}_j$ are the 
coordinate vector fields. In the expansion
\begin{equation}
\mathbf{r}_{jk}\,=\,\vg_{\!jk}^s\hs\mathbf{r}_s\,+\,\bz_{jk}\hs\mathbf{r}
\label{rjk}
\end{equation}
the coefficients $\,\vg_{\!jk}^s$ and $\,\bz_{jk}$ are the component 
functions of $\,\hs\text{\rm D}\,$ and, respectively, of the second 
fundamental form $\,\bz\,$ mentioned in (c). Thus, $\,\hs\text{\rm D}\hs\,$ is 
tor\-sion\-free and $\,\bz\,$ is \sy, since  
$\,\mathbf{r}_{jk}\nh=\mathbf{r}_{kj}$. Also, 
differentiating (\ref{rjk}), we now get, from (\ref{rlm}), 
$\,R_{jkl}{}^m\hs=\,\bz_{kl}\hs\delta_j^{\hs m}
-\hs\bz_{jl}\hs\delta_k^{\hs m}$ (and $\,\bz_{kl,j}\nh=\bz_{jl,k}$), as 
$\,\mathbf{r}_{jlk}\nh=\mathbf{r}_{klj}$. Hence 
$\,\ri^{\hs\text{\rm D}}\nnh=(1-\er)\hs\bz$, as required.
\end{remark}

\section{Flatness of an associated connection}\label{foaa}
For a fixed \tc\ $\,\hs\text{\rm D}\,$ on a \mf\ $\,\bs\,$ and a $\,1$-form 
$\,\xi\,$ on $\,\bs$, let
\begin{equation}
\lo\xi\,=\,\text{\rm D}\hs\xi\hs+\hs(\hh\text{\rm D}\hs\xi)^*
\hskip9pt\mathrm{or,\ in\ coordinates,}\hskip7pt
(\lo\xi)_{jk}\nh=\hs\xi_{\hh k,j}\nh+\xi_{\hh j,k}\hs.
\label{lxi}
\end{equation}
This defines a first-or\-der linear differential operator $\,\lo\,$ sending 
$\,1$-forms on $\,\bs\,$ to twice-co\-var\-i\-ant \sy\ \tf s. We use the 
symbol $\,\hs\text{\rm Ker}\,\lo\,$ for its kernel:
\begin{equation}
\text{\rm Ker}\,\lo\,\,=
\,\,\{\xi:\xi\text{\rm\ is\ a\ $\,1$-form\ on\ $\,\bs\,$\ and\ 
$\,\lo\xi=0$}\}\hs. 
\label{krb}
\end{equation}
The second covariant derivative $\,\hs\text{\rm D}\text{\rm D}\hs\xi\,$ of a 
$\,1$-form $\,\xi\,$ on $\,\bs\,$ can be expressed in terms of the tensor 
$\,\tau=\lo\xi$, via the identity 
$\,(\hh\text{\rm D}_w\text{\rm D}\hs\xi\hh)(u,v)
=\xi\hs R^{\text{\rm D}}\nh(v,u)w+(\hs\partial\tau)(v,u,w)$, valid for all 
\vf s $\,u,v,w\,$ on $\,\bs$, where 
$\,2\hs(\hs\partial\tau)(v,u,w)=(\hh\text{\rm D}_u\tau)(v,w)
+(\hh\text{\rm D}_w\tau)(v,u)-(\hh\text{\rm D}_v\tau)(u,w)$. (In coordinates, 
this reads $\,\xi_{\hh j,kl}=R_{jkl}{}^s\hs\xi_s
+(\tau_{jl,k}+\tau_{jk,l}-\tau_{kl,j})/2$, with 
$\,\tau_{jk}\nh=\xi_{\hh k,j}\nh+\xi_{\hh j,k}\hs$, and easily follows since 
the Ricci identity (\ref{rid}\hh-i), in its coordinate form 
$\,\xi_{\hh j,kl}-\xi_{\hh j,lk}=R_{lkj}{}^s\hs\xi_s$, gives 
$\,\tau_{jl,k}+\tau_{jk,l}-\tau_{kl,j}=2\hs\xi_{\hh j,kl}+R_{kjl}{}^s\hs\xi_s
+R_{klj}{}^s\hs\xi_s+R_{ljk}{}^s\hs\xi_s$, while 
$\,R_{klj}{}^s\hs\xi_s+R_{ljk}{}^s\hs\xi_s\nh=R_{kjl}{}^s\hs\xi_s$ by the 
first Bianchi identity.) In particular,
\begin{equation}
\text{\rm D}\text{\rm D}\hs\xi\,
=\,-\hs\xi\hs R^{\text{\rm D}}\hskip28pt\mathrm{whenever}\hskip9pt
\lo\xi=0\hs.
\label{ddr}
\end{equation}
The equality $\,\hs\text{\rm D}\text{\rm D}\hs\xi=-\hs\xi\hs R^{\text{\rm D}}$ 
means that $\,(\hh\text{\rm D}_w\text{\rm D}\hs\xi\hh)(u,v)
=\xi\hs R^{\text{\rm D}}\nh(v,u)w\,$ for all \vf s $\,u,v,w\,$ on $\,\bs$. 
The $\,2$-form $\,\hs\text{\rm D}\hs\xi\,$ is viewed here as a $\,0$-form 
valued in $\,2$-forms, which makes $\,\hs\text{\rm D}\text{\rm D}\hs\xi\,$ a 
$\,1$-form valued in $\,2$-forms, while $\,\xi\hs R^{\text{\rm D}}\nnh$, 
defined in the lines following (\ref{rid}), is now treated as a $\,1$-form 
valued in $\,2$-forms, in agreement with Remark~\ref{onetw}.

In the following lemma we solve the equation $\,\lo\xi=0$, under 
additional assumptions on $\,\bs\,$ and $\,\hs\text{\rm D}\hh$, using the 
approach of (i) in Section \ref{prel}.
\begin{lemma}\label{fltco}Given a \pftc\/ $\,\hs\text{\rm D}\hs$ on a \su\ 
$\,\bs\,$ such that the \rt\ $\,\ri^{\hs\text{\rm D}}$ is \sy, let us define a 
connection\/ $\,\nabla\hs$ in the \vb\ 
$\,\xz=(\tab)^{\wedge2}\nnh\oplus\tab\,$ by $\,\nabla(\th,\xi)
=(\hh\text{\rm D}\hh\th\hh-\xi\wedge\ri^{\hs\text{\rm D}}\nnh,
\,\text{\rm D}\hh\xi-\hs\th\hh)$. Then $\,\nabla\hs$ is flat, and the 
assignment\/ $\,\xi\mapsto(\th,\xi)\,$ given by 
$\,\th=\hs\text{\rm D}\hs\xi\,$ is a linear isomorphism of the space 
$\,\hs\text{\rm Ker}\,\lo\,$ in\/ {\rm(\ref{krb})} onto the space of all\/ 
$\,\nabla$\prl\ sections $\,(\th,\xi)\,$ of\/ $\,\xz$.
\end{lemma}
The pair $\,(\th,\xi)\,$ in Lemma~\ref{fltco} is a section of 
$\,\xz$, with the components that are a real-val\-ued $\,2$-form $\,\th\,$ and 
a real-val\-ued$\,1$-form $\,\xi\,$ on $\,\bs$. The meaning of the equality 
$\,\nabla(\th,\xi)
=(\hh\text{\rm D}\hh\th\hh-\xi\wedge\ri^{\hs\text{\rm D}}\nnh,
\,\text{\rm D}\hh\xi-\hs\th\hh)\,$ is 
$\,\nabla_{\!u}(\th,\xi)=(\hh\text{\rm D}_{\hh u}\th\hh
-\xi\wedge\ri^{\hs\text{\rm D}}\nh(u,\,\cdot\,),
\,\text{\rm D}_{\hh u}\xi-\hs\th\hh(u,\,\cdot\,))\,$ for every \vf\ $\,u\,$ 
tangent to $\,\bs$. We precede the proof of Lemma~\ref{fltco} with a remark.
\begin{remark}\label{detby}On a \mf\ $\,\bs\,$ with a fixed \tc\ 
$\,\hs\text{\rm D}\hh$, a $\,1$-form $\,\xi$ such that 
$\,\lo\xi=0\,$  is uniquely determined by its value $\,\xi_\cy$ and 
covariant derivative $\,(\hh\text{\rm D}\hs\xi\hh)_\cy$ at any given \pt\ 
$\,\cy\in\bs$. In fact, by (\ref{ddr}), the pair $\,(\th,\xi)$, with 
$\,\th=\hs\text{\rm D}\hs\xi$, then is a $\,\nabla$\prl\ section of the \vb\ 
$\,\xz=(\tab)^{\wedge2}\nnh\oplus\tab$, for the connection $\,\nabla\,$ in 
$\,\xz$ given by $\,\nabla(\th,\xi)
=(\hh\text{\rm D}\hh\th\hh+\xi\hs R^{\text{\rm D}}\nnh,
\,\text{\rm D}\hh\xi-\hs\th\hh)$. (Notation as above; 
$\,\xi\hs R^{\text{\rm D}}\nnh$, appearing in (\ref{rid}), is treated here as 
a $\,1$-form valued in $\,2$-forms, cf.\ Remark~\ref{onetw}.)
\end{remark}
\begin{proof}[Proof of Lemma\/ $\ref{fltco}$]From (\ref{cur}) with 
$\,\psi=(\th,\xi)\,$ and (\ref{dnp}), 
$\,R^\nabla\nnh(u,v)\psi=(\th\hh'\nnh,\xi\hh'\hh)$, where 
$\,\th\hh'\nh=\hat R(u,v)\hs\th
+\xi\wedge[\hh(d^{\hs\text{\rm D}}\nh\ri^{\hs\text{\rm D}}\hh)(u,v)]
+(\th\wedge\ri^{\hs\text{\rm D}}\hh)(u,v)\,$ and 
$\,\xi\hh'=-\hs(d^{\hs\text{\rm D}}\text{\rm D}\hs\xi
+\xi\wedge\ri^{\hs\text{\rm D}})\hh(u,v)$, with 
$\,d^{\hs\text{\rm D}}\text{\rm D}\hs\xi\,$ as in (\ref{rid}) and $\,\hat R\,$ 
denoting the curvature tensor of the connection that $\,\hs\text{\rm D}\hs\,$ 
induces in $\,(\tab)^{\wedge2}\nnh$. Flatness of $\,\nabla\,$ now follows 
since $\,\hat R=0\,$ (see (a) -- (d) in Section \ref{prel}), 
$\,d^{\hs\text{\rm D}}\nh\ri^{\hs\text{\rm D}}\hskip-2.5pt=0\,$ by 
Theorem~\ref{prjfl}, $\,\th\wedge\ri^{\hs\text{\rm D}}\nnh=0\,$ due to 
(\ref{tgz}), and $\,\xi\hh'=0\,$ in view of (\ref{ris}\hh-ii).

The assignment $\,\xi\mapsto(\th,\xi)\,$ is obviously injective. That it 
maps some subspace of $\,\hs\text{\rm Ker}\,\lo$ {\it onto} the space of 
$\,\nabla$\prl\ sections of $\,\xz\,$ is also clear: if $\,\nabla(\th,\xi)=0$, 
then $\,\hs\text{\rm D}\hh\xi=\hs\th\hh$, and so 
$\,\lo\xi=\hs\th\hh+\hs\th\hh^*=0$. Finally, 
$\,\nabla(\hs\text{\rm D}\hs\xi,\xi)=0\,$ whenever $\,\lo\xi=0$, in 
view of Remark~\ref{detby} and (\ref{ris}\hh-i).
\end{proof}
\begin{lemma}\label{immer}Suppose that\/ $\,\hs\text{\rm D}\,$ is a \pftc\ on 
a \sco\ \su\ $\,\bs$, the \rt\ $\,\ri^{\hs\text{\rm D}}$ is \sy, and\/ 
$\,\alpha\,$ is a fixed $\,\hs\text{\rm D}\hs$\prl\ area form on $\,\bs$, cf.\ 
{\rm(b)} in Section $\ref{prel}$. Let\/ $\,\hs\text{\rm Ker}\,\lo\,$ be the 
space given by {\rm(\ref{krb})}. Then
\begin{enumerate}
  \def\theenumi{{\rm\roman{enumi}}}
\item $\dim\hs\text{\rm Ker}\,\lo=3$,
\item a mapping $\,\fe:\bs\to\hs\text{\rm Ker}\,\lo\,$ can be defined by letting $\,\fe(\cy)$, for 
$\,\cy\in\bs$, be the unique $\,\xi\in\hs\text{\rm Ker}\,\lo\,$ with $\,\xi_\cy=0\,$ and\/ 
$\,(\hh\text{\rm D}\hs\xi\hh)_\cy=\alpha_\cy$,
\item the mapping $\,\fe:\bs\to\hs\text{\rm Ker}\,\lo\,$ defined in {\rm(ii)} is an immersion, and
\item its differential\/ $\,\dfc:\tcb\to\hs\text{\rm Ker}\,\lo\,$ at any 
$\,\cy\in\bs\,$ sends $\,v\in\tcb$ to the element $\,\eta=\dfc v\,$ of\/ 
$\,\hs\text{\rm Ker}\,\lo\,$ characterized by 
$\,\eta_\cy=-\hs\alpha_\cy(v,\,\cdot\,)\,$ and\/ 
$\,(\hh\text{\rm D}\hs\eta\hh)_\cy=0$.
\end{enumerate}
\end{lemma}
\begin{proof}Assertions (i) and (ii) are immediate from Lemma~\ref{fltco}. Now 
let $\,t\mapsto\cy=\cy(t)$ be a curve in $\,\bs$, and let 
$\,\xi=\xi(t)\in\hs\text{\rm Ker}\,\lo\,$ equal $\,\fe(\cy(t))$. With 
$\,v=v(t)\,$ and $\,\eta=\eta(t)$ standing for $\,\dot\cy\,$ and 
$\,d\hh\xi/dt$, we thus have $\,\eta=\dfc v\,$ for each $\,t$. Differentiating 
the equalities $\,\xi_\cy=0\,$ and 
$\,(\hh\text{\rm D}\hs\xi\hh)_\cy=\alpha_\cy$ covariantly along the curve, and 
noting that $\,\hs\text{\rm D}\hs\alpha=0$, we get 
$\,0=\eta_\cy+\hs\text{\rm D}_v\hh\xi=\eta_\cy+\alpha_\cy(v,\,\cdot\,)\,$ 
and $\,(\hh\text{\rm D}\hs\eta\hh)_\cy
+(\hh\text{\rm D}_v\text{\rm D}\hs\xi\hh)_\cy=0$. However, 
$\,(\hh\text{\rm D}_v\text{\rm D}\hs\xi\hh)_\cy=0\,$ by (\ref{ddr}), since 
$\,\xi_\cy=0$. This proves (iv). Consequently, (iii) follows, as 
$\,\eta_\cy\ne0$, by (iv), whenever $\,\cy\in\bs$ and $\,\eta=\dfc v\in\hs\text{\rm Ker}\,\lo\,$ 
for $\,v\in\tcb\smallsetminus\{0\}$.
\end{proof}
\begin{remark}\label{afdif}For $\,\bs,\hs\text{\rm D}\hh,\alpha\,$ and 
$\,\lo\,$ as in Lemma~\ref{immer} and another such quadruple 
$\,\bs\hs'\nnh,\hs\text{\rm D}\hs'\hs,\alpha\hs'\nnh,\lo\hs'\nnh$, let 
$\,H:\bs\to\bs\hs'$ be an affine diffeomorphism (a diffeomorphism sending 
$\,\hs\text{\rm D}\hs\,$ onto $\,\hs\text{\rm D}\hs'$). Under the 
push-for\-ward linear isomorphism 
$\,\hs\text{\rm Ker}\,\lo\to\hs\text{\rm Ker}\,\lo\hs'\nnh$, induced by $\,H$, 
the immersion $\,\fe$ of Lemma~\ref{immer} obviously corresponds to the 
composite in which the analogous immersion 
$\,\fe\hs'\nnh:\bs\hs'\nnh\to\hs\text{\rm Ker}\,\lo\hs'$ is followed by the 
isomorphism $\,\hs\text{\rm Ker}\,\lo\hs'\nnh\to\hs\text{\rm Ker}\,\lo\hs'$ of 
multiplication by the constant $\,c\,$ such that 
$\,H^*\nh\alpha\hs'\nh=c^{-1}\alpha$.
\end{remark}
The following lemma will be needed in Section \ref{scoc}.
\begin{lemma}\label{parvf}Given a \pftc\/ $\,\hs\text{\rm D}\hs$ on a \sco\ 
\su\/ $\,\bs\,$ such that the \rt\/ $\,\ri^{\hs\text{\rm D}}$ is \sy, let\/ 
$\,\fe:\bs\to\hs\text{\rm Ker}\,\lo\,$ be as in\/ Lemma\/ $\ref{immer}$. For 
any \vf\/ $\,v\,$ on\/ $\,\bs$, the following two conditions are 
equivalent\/{\rm:}
\begin{enumerate}
  \def\theenumi{{\rm\alph{enumi}}}
\item $v\,$ is\/ $\,\hs\text{\rm D}\hs$\prl\/{\rm;}
\item the function $\,\cy\mapsto\dfc v_\cy$, valued in 
$\,\hs\text{\rm Ker}\,\lo$, is constant on\/ $\,\bs$.
\end{enumerate}
\end{lemma}
\begin{proof}If $\,v\,$ is $\,\hs\text{\rm D}\hs$\prl, so is the $\,1$-form 
$\,\eta=-\hs\alpha(v,\,\cdot\,)$, where $\,\alpha\,$ is a fixed 
$\,\hs\text{\rm D}\hs$\prl\ area form on $\,\bs$. The same $\,\eta\,$ thus 
satisfies the conditions listed in Lemma~\ref{immer}(iv) at all \pt s 
$\,\cy\in\bs\,$ simultaneously, so that $\,\dfc v_\cy=\eta\,$ for all 
$\,\cy\in\bs$, and (b) follows. Conversely, if the function 
$\,\cy\mapsto\dfc v_\cy$ is constant, and equal to 
$\,\eta\in\hs\text{\rm Ker}\,\lo$, then, by Lemma~\ref{immer}(iv), $\,\eta\,$ 
is $\,\hs\text{\rm D}\hs$\prl\ and $\,\eta=-\hs\alpha(v,\,\cdot\,)$. Hence 
$\,v\,$ is $\,\hs\text{\rm D}\hs$\prl\ as well.
\end{proof}
For a cen\-tro\-af\-fine connection $\,\hs\text{\rm D}\hs\,$ on a \su\ 
$\,\bs$, defined as in Section \ref{pfid} for some $\,\hy\,$ and $\,\im\,$ 
with $\,\er=2$, we have an explicit description of the space 
$\,\hs\text{\rm Ker}\,\lo\,$ given by (\ref{krb}). Namely, any fixed volume 
form $\,\vol\in(\hy^*)^{\wedge3}\smallsetminus\{0\}\,$ gives rise to an 
isomorphism $\,\hy\to\text{\rm Ker}\,\lo\,$ sending $\,\mathbf{w}\in\hy$ to 
the $\,1$-form $\,\xi\,$ on $\,\bs\,$ with 
$\,\xi_\cy(v)=\vol(\mathbf{w},\im(\cy),\dic v)$ for $\,\cy\in\bs\,$ and 
$\,v\in\tcb$. That $\,\xi\in\text{\rm Ker}\,\lo\,$ is easily seen in 
the notation of (\ref{rjk}). Specifically, $\,\xi\,$ has the components 
$\,\xi_{\hh j}=\xi(\mathbf{r}_j)
=\vol(\mathbf{w},\mathbf{r},\mathbf{r}_j)$, and so, by (\ref{rjk}), 
$\,\xi_{\hh j,k}=\hs\partial_k\hh\xi_{\hh j}-\vg_{\!kj}^s\hs\xi_{\hh s}
=\vol(\mathbf{w},\mathbf{r}_k,\mathbf{r}_j)$, which is \skc\ in $\,j,k$. 
The assignment $\,\mathbf{w}\mapsto\xi\,$ is injective: if 
$\,\mathbf{w}\in\hy\,$ is nonzero, the transversality assumption about 
$\,\im\,$ guarantees that $\,\mathbf{w},\im(\cy)\,$ and $\,\dic v\,$ are 
linearly independent for some $\,\cy\in\bs\,$ and $\,v\in\tcb$. Since 
$\,\dim\,\text{\rm Ker}\,\lo\le3\,$ by Lemma~\ref{immer}, it 
follows that $\,\mathbf{w}\mapsto\xi\,$ is an isomorphism.

Furthermore, the immersion $\,\im:\bs\to\hy\,$ is {\it equi\-af\-fine\/} 
relative to $\,\hs\text{\rm D}\hh$, in the sense that, for some (or any) fixed 
$\,\vol\in(\hy^*)^{\wedge3}\smallsetminus\{0\}$, the formula 
$\,\alpha_\cy(u,v)=\vol(\im(\cy),\dic u,\dic v)$, for $\,\cy\in\bs\,$ and 
$\,v\in\tcb$, defines a $\,\hs\text{\rm D}\hs$\prl\ area form $\,\alpha\,$ on 
$\,\bs$. This is clear since 
$\,\alpha_{jk}=\vol(\mathbf{r},\mathbf{r}_j,\mathbf{r}_k)$, and so (\ref{rjk}) 
gives 
$\,\alpha_{jk,l}=\hs\partial_l\hh\alpha_{jk}-\vg_{\!lj}^s\hs\alpha_{sk}
-\vg_{\!lk}^s\hs\alpha_{js}=0$. (Note that $\,\vol\,$ is constant, and 
$\,\vol(\mathbf{r}_l,\mathbf{r}_j,\mathbf{r}_k)=0\,$ as 
$\,\mathbf{r}_l,\mathbf{r}_j,\mathbf{r}_k$ are tangent to $\,\bs$.)
\begin{remark}\label{cafed}Let us fix $\,\vol$, define $\,\alpha$, and 
identify $\,\hy\,$ with $\,\hs\text{\rm Ker}\,\lo\,$ as above. Then 
$\,\fe=\im$, for the immersion $\,\fe\,$ of Lemma~\ref{immer}. This is clear 
from the relations 
$\,\alpha_{jk}=\vol(\mathbf{r},\mathbf{r}_j,\mathbf{r}_k)=\xi_{\hh k,j}$, cf.\ 
Remark~\ref{twice}.
\end{remark}
\begin{remark}\label{cptim}If an immersion $\,\fe\,$ of a closed \su\ 
$\,\bs\,$ in a real $\,3$-space $\,\hy\hh\,$ is transverse to all lines 
through $\,0$, then $\,\bs\,$ is diffeomorphic to $\,S^2$ and $\,\fe\,$ is an 
embedding. In fact, let $\,S=\{\eta\in\hy:|\eta|=1\}\,$ for a fixed Euclidean 
norm $\,|\hskip3pt|\,$ in $\,\hy$. The locally diffeomorphic mapping 
$\,\bs\ni\cy\mapsto\fe(\cy)/|\fe(\cy)|\in S\,$ must be a covering, and hence a 
diffeomorphism.
\end{remark}

\section{Local classification}\label{locl}
Since every immersion is, locally, an embedding, the following theorem 
provides a complete local classification of \pftc s on \su s with \sy\ \rt s, 
which is a special case of a result of Kurita \cite{kurita}.
\begin{theorem}\label{lccls}There exists a natural bijective correspondence 
between the equivalence classes of
\begin{enumerate}
  \def\theenumi{{\rm\alph{enumi}}}
\item[(a)] pairs $\,(\bs,\hs\text{\rm D}\hh)\,$ formed by a \sco\ \su\ 
$\,\bs\,$ and a \pftc\ $\,\hs\text{\rm D}\hs$ on $\,\bs\,$ such that the 
\rt\ of\/ $\,\hs\text{\rm D}\,$ is \sy\ and the immersion 
$\,\fe:\bs\to\hs\text{\rm Ker}\,\lo\,$ defined in\/ Lemma\/ $\ref{immer}$ is an embedding,
\end{enumerate}
and the equivalence classes of
\begin{enumerate}
  \def\theenumi{{\rm\alph{enumi}}}
\item[(b)] \sco\ \su s $\,S\hs$ embedded in a fixed\/ $\,3$\diml\ \rvs\ 
$\,\hy\,$ and transverse to all lines in\/ $\,\hy\hs$ containing $\,0$.
\end{enumerate}
The equivalence relation in question consists in being congruent under a 
specific class of transformations\/{\rm:} affine diffeomorphisms of \su s with 
connections for {\rm(a)}, linear isomorphisms\/ $\,\hy\to\hy$ for {\rm(b)}.

Explicitly, the bijective correspodence assigns to the equivalence class of a 
pair $\,(\bs,\hs\text{\rm D}\hh)$ the equivalence class of the \su\ 
$\,S=\hs\ism(F(\bs))\subset\hy$, where $\,\ism\,$ is any linear isomorphism 
$\,\hs\text{\rm Ker}\,\lo\to\hy$. The inverse assignment sends a \su\ 
$\,S\subset\hy\,$ to $\,(\bs,\hs\text{\rm D}\hh)$, where $\,\bs=S\,$ and\/ 
$\,\hs\text{\rm D}\hs\,$ is the cen\-tro\-af\-fine connection of\/ $\,S$, 
described in Section $\ref{pfid}$.
\end{theorem}
\begin{proof}That, for a pair $\,(\bs,\hs\text{\rm D}\hh)\,$ as in (a), 
$\,S=\hs\ism(F(\bs))\,$ has the properties named in (b), is obvious from 
Lemma~\ref{immer}. Conversely, for $\,S\,$ and $\,\hy\,$ as in (b), the 
conditions listed in (a) are satisfied by the cen\-tro\-af\-fine connection 
$\,\hs\text{\rm D}\hs\,$ on $\,\bs=S$. (See Remarks~\ref{ctraf} 
and~\ref{cafed}, for the inclusion mapping $\,\im:\bs\to\hy$.) Both 
assignments are well defined on equivalence classes (cf.\ Remark~\ref{afdif}). 
It now remains to be shown that the two mappings between sets of equivalence 
classes are each other's inverses.

First, if $\,S,\hy\,$ are as in (b), $\,S\,$ coincides, by Remark~\ref{cafed}, 
with the image $\,\fe(\bs)\,$ of the immersion $\,\fe$ in Lemma~\ref{immer} 
for the cen\-tro\-af\-fine connection $\,\hs\text{\rm D}\hs\,$ on $\,\bs=S$.

Conversely, for $\,(\bs,\hs\text{\rm D}\hh)\,$ as in (a), the diffeomorphism 
$\,\fe:\bs\to\fe(\bs)$ sends $\,\hs\text{\rm D}\hs\,$ onto the 
cen\-tro\-af\-fine connection on the \su\ $\,\fe(\bs)$. In fact, suppose that 
a \vf\ $\,t\mapsto w=w(t)\in T_{\cy(t)}\bs\,$ is tangent to $\,\bs\,$ along a 
curve $\,t\mapsto\cy=\cy(t)\,$ in $\,\bs$, while the symbols 
$\,\dot\cy=\dot\cy(t)\,$ and $\,\eta=\eta(t)\,$ denote $\,d\cy/dt\,$ and 
$\,\dfc\dot\cy=\dot\xi\,$ (i.e., $\,d\hh\xi/dt$), and $\,\zeta=\zeta(t)$ 
stands for $\,\dfc w$. Differentiating the equality 
$\,\zeta_\cy=-\hs\alpha_\cy(w,\,\cdot\,)\,$ covariantly along the curve, and 
noting that $\,(\hh\text{\rm D}\hs\zeta\hh)_\cy=0\,$ (see 
Lemma~\ref{immer}(iv)), we obtain 
$\,\dot\zeta_\cy=-\hs\alpha_\cy(\hh\text{\rm D}_{\dot\cy}w,\,\cdot\,)$. Thus, 
again by Lemma~\ref{immer}(iv), $\,\dot\zeta=d\hh\zeta/dt\,$ differs from the 
$\,\dfc$-im\-age of $\,\hs\text{\rm D}_{\dot\cy}w$, at $\,\cy=\cy(t)$, by an 
element of $\,\hs\text{\rm Ker}\,\lo\,$ which is a $\,1$-form vanishing at 
$\,\cy$, and hence equal to a scalar times the normal vector $\,\fe(\cy(t))$. 
Thus, $\,\dfc$ sends $\,\hs\text{\rm D}_{\dot\cy}w\,$ to the tangent component 
of $\,\dot\zeta_\cy$.
\end{proof}

\section{Special classes of connections}\label{scoc}
In this section we use Theorem~\ref{lccls} to classify \pftc s 
$\,\hs\text{\rm D}\hs\,$ on \su s $\,\bs$, with \sy\ \rt s 
$\,\ri^{\hs\text{\rm D}}\nnh$, which satisfy further 
restrictive conditions. The conditions in question are local symmetry 
($\text{\rm D}R^{\text{\rm D}}\nnh=0$, equivalent, by (\ref{rri}), to 
$\,\hs\text{\rm D}\ri^{\hs\text{\rm D}}\nnh=0$) and 
(Ric\-ci-)\-re\-cur\-rence. Both results are well known, cf.\ 
\cite{opozda-91}, \cite{nomizu-opozda}.
\begin{theorem}\label{lcsym}Let\/ $\,\hs\text{\rm D}\hs$ be a \pftc\ on a 
\sco\ \su\ $\,\bs$, for which the \rt\/ $\,\ri^{\hs\text{\rm D}}$ is \sy\ and 
the immersion $\,\fe\,$ defined in\/ Lemma\/ $\ref{immer}$ is an embedding 
of\/ $\,\bs\,$ in the $\,3$\diml\ \vs\ $\,\hs\text{\rm Ker}\,\lo$.
\begin{enumerate}
  \def\theenumi{{\rm\roman{enumi}}}
\item Suppose, in addition, that\/ $\,\ri^{\hs\text{\rm D}}$ is\/ 
$\,\hs\text{\rm D}\hs$\prl. Then there exists a \sy\ bilinear form\/ $\,\lr\,$ 
in $\,\hs\text{\rm Ker}\,\lo\,$ such that\/ $\,\lg\fe(\cy),\fe(\cy)\rg=1\,$ 
for every\/ $\,\cy\in\bs\,$ and\/ $\,\lr\,$ has the algebraic type of the 
direct sum of a pos\-i\-tive-def\-i\-nite form in dimension $\,1\,$ and\/ 
$\,\ri^{\hs\text{\rm D}}\nnh$, at any \pt\ of\/ $\,\bs$. Thus, the image 
$\,\fe(\bs)\,$ is a relatively open subset of an algebraic \su\ $\,S\,$ in 
$\,\hs\text{\rm Ker}\,\lo$, and\/ $\,S\,$ itself may be described as follows.
\begin{enumerate}
  \def\theenumi{{\rm\alph{enumi}}}
\item[(a)] If\/ $\,\ri^{\hs\text{\rm D}}\nnh=0$, i.e., 
$\,\hs\text{\rm D}\hs\,$ is flat, $\,S\,$ is a plane not containing $\,0$.
\item[(b)] If\/ $\,\ri^{\hs\text{\rm D}}$ is\/ $\,\hs\text{\rm D}\hs$\prl, of 
rank $\,1$, and positive sem\-i\-def\-i\-nite, $\,S\,$ is an elliptic cylinder 
and its center axis contains $\,0$.
\item[(c)] If\/ $\,\ri^{\hs\text{\rm D}}$ is\/ $\,\hs\text{\rm D}\hs$\prl, 
of rank $\,1$, and negative sem\-i\-def\-i\-nite, $\,S\,$ is a hyperbolic 
cylinder whose center axis contains $\,0$.
\item[(d)] If\/ $\,\ri^{\hs\text{\rm D}}$ is\/ $\,\hs\text{\rm D}\hs$\prl, 
nondegenerate, and positive definite or negative definite or, respectively, 
indefinite, $\,S\,$ is an ellipsoid or a two-sheet\-ed hyperboloid or, 
respectively, a one-sheet\-ed hyperboloid, centered at\/ $\,0$.
\end{enumerate}
\item Conversely, if\/ $\,\fe(\bs)\,$ is contained in an algebraic \su\ 
$\,S\,$ in $\,\hs\text{\rm Ker}\,\lo\,$ with the properties listed in\/ 
{\rm(a)}, {\rm(b)}, {\rm(c)} or\/ {\rm(d)}, then $\,\ri^{\hs\text{\rm D}}$ 
is\/ $\,\hs\text{\rm D}\hs$\prl\ and has the algebraic type named in\/ 
{\rm(a)} -- {\rm(d)}.
\end{enumerate}
\end{theorem}
\begin{proof}Using a $\,\hs\text{\rm D}\hs$\prl\ area form $\,\alpha\,$ on 
$\,\bs\,$ (cf.\ (a) -- (d) in Section \ref{prel}), we define $\,\lr\,$ to be 
the \sy\ bilinear mapping assigning to $\,\xi,\eta\in\hs\text{\rm Ker}\,\lo\,$ 
the function $\,\lg\xi,\eta\rg$ on $\,\bs\,$ given by 
$\,\lg\xi,\eta\rg=\phi\psi+\ri^{\hs\text{\rm D}}\nh(v,w)$, with the functions 
$\,\phi,\psi\,$ and \vf s $\,v,w$ such that 
$\,\text{\rm D}\hs\xi=\phi\hh\alpha$, $\,\text{\rm D}\hs\eta=\psi\hh\alpha$, 
$\,\xi=\alpha(v,\,\cdot\,)\,$ and $\,\eta=\alpha(w,\,\cdot\,)$.

If $\,\hs\text{\rm D}\ri^{\hs\text{\rm D}}\nnh=0$, then, for any pair 
$\,\xi,\eta\in\hs\text{\rm Ker}\,\lo$, the function $\,\lg\xi,\eta\rg\,$ is constant. To 
see this, let us fix an arbitrary \vf\ $\,u\,$ on $\,\bs$, and assume that 
neither $\,\xi\,$ nor $\,\eta\,$ is identically zero. For $\,v,w,\phi,\psi\,$ 
determined by $\,\xi\,$ and $\,\eta\,$ as above, 
$\,\alpha(\hh\text{\rm D}_uv,\,\cdot\,)
=\hs\text{\rm D}_u\xi=\phi\hh\alpha(u,\,\cdot\,)$, so that 
$\,\text{\rm D}_uv=\phi\hh u$, and, similarly, 
$\,\text{\rm D}_uw=\psi\hh u$. On the other hand, by (\ref{ddr}) and 
\hbox{(\ref{ris}\hh-i),} 
$\,(d_u\phi)\hs\alpha=\hs\text{\rm D}_u\text{\rm D}\hs\xi
=\xi\wedge\ri^{\hs\text{\rm D}}\nh(u,\,\cdot\,)$, where $\,d_u$ is the 
directional derivative. Hence 
$\,(d_u\phi)\hs\xi=(d_u\phi)\hs\alpha(v,\,\cdot\,)
=-\hs\ri^{\hs\text{\rm D}}\nh(u,v)\hs\xi$, as $\,\xi(v)=\alpha(v,v)=0$. 
Noting that, by Remark~\ref{detby}, $\,\xi\ne0\,$ on a dense open subset of 
$\,\bs$, we now obtain $\,d_u\phi=-\hs\ri^{\hs\text{\rm D}}\nh(u,v)$, and, 
similarly, $\,d_u\psi=-\hs\ri^{\hs\text{\rm D}}\nh(u,w)$. Combining our 
formulae for $\,\text{\rm D}_uv$, $\text{\rm D}_uw$, $d_u\phi\,$ and 
$\,d_u\psi\,$ we see that $\,d_u\lg\xi,\eta\rg=0$, as required.

Assertions (i) and (a) -- (d) are now immediate. Conversely, under the 
assumption of (ii), let $\,\lr\,$ be a \sy\ bilinear form on $\,\hs\text{\rm Ker}\,\lo\,$ such 
that $\,\lg\xi,\xi\rg=1\,$ for every $\,\xi\in\fe(\bs)$. In the 
notation of (\ref{rjk}) we thus have 
$\,\lg\hh\mathbf{r},\hh\mathbf{r}\hh\rg=1$, and partial 
differentiation gives $\,\lg\hh\mathbf{r},\hh\mathbf{r}_j\rg=0$. 
Applying $\,\lg\hh\mathbf{r},\,\cdot\,\rg\,$ to (\ref{rjk}) we now get 
$\,\bz_{jk}=\lg\hh\mathbf{r},\hh\mathbf{r}_{jk}\rg$, that is, 
$\,\bz_{jk}=-\hs\lg\hh\mathbf{r}_j,\hh\mathbf{r}_k\rg$. Hence, again 
by (\ref{rjk}), $\,\bz_{jk,l}=\hs\partial_l\hh\bz_{jk}
-\vg_{\!lj}^s\hs\bz_{sk}-\vg_{\!lk}^s\hs\bz_{js}$ vanishes 
identically. As $\,\bz=-\hs\ri^{\hs\text{\rm D}}$ (see Remark~\ref{ctraf}(c)), 
this completes the proof.
\end{proof}
A \tf\ $\,B\,$ on a \mf\ endowed with a fixed connection 
$\,\hs\text{\rm D}\,$ is called {\it recurrent\/} if, for every tangent \vf\ 
$\,v$, the tensors $\,B\,$ and $\,\hs\text{\rm D}_vB\,$ are linearly dependent 
at every \pt. Equivalently, in the open set $\,\,U\,$ on which $\,B\ne0$, 
one has $\,\hs\text{\rm D}B=\zeta\otimes B$ for some $\,1$-form $\,\zeta$. The 
connection $\,\hs\text{\rm D}\hs\,$ itself is said to be {\it recurrent\hh}, 
or {\it Ric\-ci-re\-cur\-rent\hh}, if its curvature tensor (or, respectively, 
\rt) is recurrent. When this is the case for the \lcc\ of a \prd\ $\,(M,g)$, 
one refers to $\,(M,g)\,$ as a {\it re\-cur\-rent\/} or, respectively, {\it 
Ric\-ci-re\-cur\-rent\/} manifold. Finally, $\,(M,g)\,$ is called {\it 
conformally recurrent\/} if its Weyl conformal curvature tensor is recurrent.

For a \tc\ with a \sy\ \rt\ on a \su, being recurrent means, by (\ref{rri}), 
the same as being \rr. Such connections are called here `\rr' rather than 
`recurrent' (which will be convenient later, in Part III). 
\begin{theorem}\label{rcrec}Suppose that\/ $\,\hs\text{\rm D}\hs$ is a \pftc\ 
on a \sco\ \su\ $\,\bs\,$ and the \rt\/ $\,\ri^{\hs\text{\rm D}}$ is \sy, 
while\/ $\,\ri^{\hs\text{\rm D}}\nnh\ne0\,$ and\/ 
$\,\hs\text{\rm D}\hs\ri^{\hs\text{\rm D}}\nnh\ne0$ everywhere in\/ $\,\bs$. 
In addition, let the immersion $\,\fe:\bs\to\hs\text{\rm Ker}\,\lo\,$ defined 
in\/ {\rm Lemma~\ref{immer}} be an embedding. The following four conditions 
then are equivalent\/{\rm:}
\begin{enumerate}
  \def\theenumi{{\rm\alph{enumi}}}
\item $\text{\rm D}\hs\,$ is \rr,
\item there exists a nonzero $\,\hs\text{\rm D}\hs$\prl\ \vf\ on\/ 
$\,\bs$,
\item $\bs\,$ admits a nonzero $\,\hs\text{\rm D}\hs$\prl\ $\,1$-form,
\item the image \su\ $\,S=\fe(\bs)\,$ is a cylinder, that is, a union of 
mutually parallel line segments in\/ $\,\hs\text{\rm Ker}\,\lo$.
\end{enumerate}
\end{theorem}
\begin{proof}If $\,\hs\text{\rm D}\hs\,$ is \rr, and so 
$\,\hs\text{\rm D}\ri^{\hs\text{\rm D}}\nnh=\zeta\otimes\ri^{\hs\text{\rm D}}$ 
for some no\-where-van\-ish\-ing $\,1$-form $\,\zeta$, the Codazzi equation 
$\,d^{\hs\text{\rm D}}\nh\ri^{\hs\text{\rm D}}\hskip-2.5pt=\nh0\,$ gives 
$\,\zeta\wedge\ri^{\hs\text{\rm D}}\nnh=0$. Thus, $\,\zeta\,$ spans the image 
$\,\jm\,$ of the \vbm\ $\,\tb\to\tab\,$ sending each \vf\ $\,u\,$ to 
$\,\ri^{\hs\text{\rm D}}\nh(u,\,\cdot\,)$. (In fact, 
$\,\zeta\wedge\ri^{\hs\text{\rm D}}\nh(u,\,\cdot\,)=0$, and so 
$\,\ri^{\hs\text{\rm D}}\nh(u,\,\cdot\,)\,$ equals a function times 
$\,\zeta$.) Hence $\,\jm\,$ is a line sub\bd\ of $\,\tab\,$ and, as 
$\,\ri^{\hs\text{\rm D}}$ is recurrent, $\,\jm\,$ must be 
$\,\hs\text{\rm D}\hs$\prl\ (invariant under 
$\,\hs\text{\rm D}\hs$\prl\ transports), which in turn implies that the 
$\,1$-form $\,\zeta\,$ spanning $\,\jm\,$ is recurrent, i.e., 
$\,\hs\text{\rm D}\hs\zeta=\eta\otimes\zeta\,$ for some $\,1$-form $\,\eta$. 
The Ricci identity (\ref{ris}\hh-ii) gives 
$\,d^{\hs\text{\rm D}}\text{\rm D}\hs\zeta=0$, so that 
$\,0=d^{\hs\text{\rm D}}\text{\rm D}\hs\zeta
=d^{\hs\text{\rm D}}\nh(\eta\otimes\zeta)=(d\hskip.3pt\eta)\otimes\zeta$. 
Consequently, $\,d\hskip.3pt\eta=0\,$ and, as $\,\bs\,$ is \sco, 
$\,\eta=\df\,$ for some function $\,f:\bs\to\bbR\hs$. Since 
$\,\hs\text{\rm D}\hs\zeta=df\otimes\zeta$, the $\,1$-form 
$\,\xi=e^{-f}\zeta\,$ is $\,\hs\text{\rm D}\hs$\prl. Therefore, (a) implies 
(c).

Equivalence of (b) and (c) is clear since a fixed $\,\hs\text{\rm D}\hs$\prl\ 
area form $\,\alpha\,$ on $\,\bs\,$ leads to a correspondence 
$\,v\mapsto\xi=\alpha(v,\,\cdot\,)\,$ between \vf s $\,v\,$ and $\,1$-forms 
$\,\xi$. On the other hand, (b) and (d) are equivalent as a consequence of 
Lemma~\ref{parvf}. (Note that, if $\,\fe(\bs)\,$ is a cylinder, some nonzero 
vector in $\,\hs\text{\rm Ker}\,\lo\,$ is tangent to $\,\fe(\bs)\,$ at every 
\pt.)

Finally, (c) implies (a). In fact, if $\,\hs\text{\rm D}\hs\xi=0\,$ for a 
$\,1$-form $\,\xi\ne0$, the Ricci identity (\ref{ris}\hh-ii) yields 
$\,\xi\wedge\ri^{\hs\text{\rm D}}\nnh=0$, and so the line sub\bd\ $\,\jm\,$ of 
$\,\tab\,$ spanned by $\,\xi\,$ contains the image of the \vbm\ 
$\,u\mapsto\ri^{\hs\text{\rm D}}\nh(u,\,\cdot\,)$. (See the beginning of this 
proof.) As $\,\jm\,$ is $\,\hs\text{\rm D}\hs$\prl, $\,\ri^{\hs\text{\rm D}}$ 
must be recurrent.
\end{proof}
\begin{remark}\label{nscgc}If a \pftc\ $\,\hs\text{\rm D}\hs$ on a \sco\ \su\ 
$\,\bs\,$ is \rr, $\,\ri^{\hs\text{\rm D}}$ is \sy, and $\,\fe(\bs)\,$ 
contains no nontrivial line segment, then either $\,\ri^{\hs\text{\rm D}}$ or 
$\,-\hh\ri^{\hs\text{\rm D}}$ is a pos\-i\-tive-def\-i\-nite metric on 
$\,\bs\,$ with nonzero constant Gaussian curvature. In fact, 
$\,\ri^{\hs\text{\rm D}}$ must be parallel and definite, for if it were 
parallel but not definite, or non\prl, $\,\fe(U)\,$ would be a union of line 
segments, for some nonempty open set $\,\,U\subset\bs$. (See 
Theorem~\ref{lcsym}(i) and (d) in Theorem~\ref{rcrec}.) Thus, 
$\,\pm\hh\ri^{\hs\text{\rm D}}$ is a Riemannian metric with the \lcc\ 
$\,\hs\text{\rm D}\hs$, and $\,\hs\text{\rm D}\ri^{\hs\text{\rm D}}\nnh=0$.
\end{remark}

\section{Equi\-af\-fine connections}\label{eaco}
Let $\,\hs\text{\rm D}\hs\,$ be a \tc\ on a \su\ $\,\bs$. By a 
$\,\hs\text{\rm D}\hs$\prl\ {\it area element\/} we mean any nonzero 
$\,\hs\text{\rm D}\hs$\prl\ differential $\,2$-form $\,\pm\hs\alpha\,$ on 
$\,\bs\,$ defined, at each \pt, only up to a sign. (The sign $\,\pm\,$ 
indicates its double-val\-ued\-ness.) We will call $\,\hs\text{\rm D}\hs\,$ an 
{\it equi\-af\-fine\/} connection if a $\,\hs\text{\rm D}\hs$\prl\ area 
element exists on $\,\bs$.

By (a) -- (d) in Section \ref{prel}, $\,\hs\text{\rm D}\hs\,$ is {\it locally} 
\ea\ if and only if its \rt\ $\,\ri^{\hs\text{\rm D}}$ is \sy. In the 
case where $\,\bs\,$ is \sco, symmetry of $\,\ri^{\hs\text{\rm D}}$ implies 
(global) equi\-af\-fin\-i\-ty of $\,\hs\text{\rm D}\hs\,$ on $\,\bs$, and a 
single-val\-ued $\,\hs\text{\rm D}\hs$\prl\ area form on $\,\bs\,$ exists as 
well.
\begin{remark}\label{csunm}For an \ea\ \tc\ $\,\hs\text{\rm D}\hs\,$ on a 
{\it closed\/} \su\ $\,\bs$, every diffeomorphism $\,\bs\to\bs\,$ which is 
affine (that is, sends $\,\hs\text{\rm D}\hs\,$ onto itself) is also {\it 
unimodular\/} in the sense of preserving some (or any) 
$\,\hs\text{\rm D}\hs$\prl\ area element $\,\pm\hs\alpha\,$ on $\,\bs$. In 
fact, an affine diffeomorphism obviously sends $\,\pm\hs\alpha\,$ onto a 
constant multiple $\,\pm\hs c\hs\alpha$. However, since $\,\pm\hs\alpha\,$ 
constitutes a smooth positive measure density, invariance of the area 
$\,\int_\bs\alpha\,$ under diffeomorphisms implies that $\,c=\pm\hs1$.
\end{remark}
Given a \su\ $\,\bs\,$ with a \pftc\ $\,\hs\text{\rm D}\hs\,$ such that 
$\,\ri^{\hs\text{\rm D}}$ is \sy, let $\,\hs\widehat{\text{\rm D}}\hs\,$ be 
the pull\-back of $\,\hs\text{\rm D}\hs\,$ to the universal covering \su\ 
$\,\widehat\bs\,$ of $\,\bs$. While $\,\hs\widehat{\text{\rm D}}\hs\,$ is 
always \ea, it is clear that $\,\hs\text{\rm D}\hs\,$ is \ea\ on 
$\,\bs\,$ if and only if the affine diffeomorphisms 
$\,\widehat\bs\to\widehat\bs\,$ forming the deck transformation group 
$\,\pi_1\bs\,$ are all unimodular. On the other hand, every deck 
transformation in $\,\pi_1\bs\,$ gives rise to the corresponding 
push-for\-ward linear isomorphism 
$\,\hs\text{\rm Ker}\,\widehat{\lo}\to\hs
\text{\rm Ker}\,\widehat{\lo}\,$ of the $\,3$\diml\ space 
$\,\hs\text{\rm Ker}\,\widehat{\lo}\,$ defined as in (\ref{krb}) for 
$\,\hs\widehat{\text{\rm D}}\hs\,$ and $\,\widehat\bs\,$ rather than 
$\,\hs\text{\rm D}\hs\,$ and $\,\bs$. (Cf.\ Lemma~\ref{immer}(i).) In other 
words, $\,\pi_1\bs\,$ has a natural linear representation in 
$\,\hs\text{\rm Ker}\,\widehat{\lo}$, and the immersion 
$\,\widehat\fe:\widehat\bs\to\widehat{\lo}\,$ defined as in Lemma~\ref{immer}  
is obviously equivariant relative to it.
\begin{remark}\label{unmsl}For $\,\bs,\hs\text{\rm D}\hs,\widehat\bs\,$ and 
$\,\hs\widehat{\text{\rm D}}\hs\,$ as in the last paragraph, 
$\,\hs\text{\rm D}\hs\,$ is \ea\ on $\,\bs$ if and only if the representation 
of $\,\pi_1\bs\,$ in $\,\hs\text{\rm Ker}\,\widehat{\lo}\,$ just described 
consists of operators with determinant $\,\pm\hs1$. (In fact, according to 
Remark~\ref{cafed} and the paragraph preceding it, choosing a 
$\,\hs\widehat{\text{\rm D}}\hs$\prl\ area form in $\,\widehat\bs\,$ amounts 
to fixing a volume form in $\,\hs\text{\rm Ker}\,\widehat{\lo}$.)
\end{remark}
Let $\,\bs\,$ be a fixed closed \su. Any \pftc\ $\,\hs\text{\rm D}\hs\,$ on 
$\,\bs\,$ such that $\,\ri^{\hs\text{\rm D}}$ is \sy\ can be constructed as 
follows, in terms of the universal covering \su\ $\,\widehat\bs\,$ of 
$\,\bs\,$ and the action on $\,\widehat\bs\,$ of the deck transformation group 
$\,\pi_1\bs$. We begin by choosing linear representation of $\,\pi_1\bs\,$ in 
a real $\,3$-space $\,\hy$. In view of Theorem~\ref{lccls}, we now only 
need to prescribe an equivariant immersion 
$\,\widehat\im:\widehat\bs\to\hy$, transverse to all lines through $\,0$. 
To this end, we choose a fundamental domain for the action of $\,\pi_1\bs\,$ 
on $\,\widehat\bs$, in the form of a curvilinear polygon 
$\,Q\subset\widehat\bs\,$ on which the action of $\,\pi_1\bs\,$ realizes 
standard identifications between some pairs of edges. (See 
\cite[pp.\ 148--149]{spanier}.) Our $\,\widehat\im$, chosen arbitrarily on a 
small neighborhood of one vertex of $\,Q$, is then propagated to neighborhoods 
of other vertices, with the aid of linear transformations assigned to standard 
generators of $\,\pi_1\bs$. Next, $\,\widehat\im\,$ is extended from (smaller) 
neighborhoods of the vertices to narrow tubular neighborhoods of a half of the 
total number of edges, chosen so as to be either pairwise disjoint ($\bs\,$ 
non\-orient\-able), or so that each selected edge shares a vertex with exactly 
one other selected edge ($\bs\,$ orient\-able). Generators of $\,\pi_1\bs\,$ 
are used, again, to propagate the immersion $\,\widehat\im\,$ to neighborhoods 
of the remaining edges. Finally, we extend $\,\widehat\im\,$ to the interior 
of $\,Q$, leaving it unchanged near the boundary, and from there to all of 
$\,\widehat\bs$, in a manner uniquely determined by the requirement of 
equivariance.

Not all steps of this process are always possible, as the transversality 
requirement may preclude extensibility of the immersion. For instance, if 
$\,\bs\,$ is a closed \su\ and we try to immerse in $\,\hy\,$ not just 
$\,\widehat\bs$, but also $\,\bs\,$ itself (which amounts to choosing the 
trivial representation), the immersion $\,\widehat\im\,$ will not exist except 
when $\,\bs\,$ is diffeomorphic to the $\,2$-sphere (Remark~\ref{cptim}).
\begin{remark}\label{divfo}The divergence formula 
$\,\int_\bs(\text{\rm div}{}^{\hs\text{\rm D}}\nh w)\,\alpha=0\,$ for 
compactly supported $\,C^1$ vector fields $\,w\,$ on a \mf\ $\,\bs\,$ remains 
valid also when $\,\hs\text{\rm D}\hs\,$ is a just \tc\ on $\,\bs\,$ with a 
$\,\hs\text{\rm D}\hs$\prl\ volume element $\,\pm\hs\alpha\,$ (defined 
analogously as an area element in the case of surfaces). 

In fact, using a finite partition of unity, we may reduce the question to 
the case where the compact support of $\,w\,$ is contained in the domain of a 
coordinate system $\,\y^{\hs j}$. Since $\,\hs\text{\rm D}\hs\alpha=0$, we 
have $\,\vg_{\!jk}^k=\partial_j\log|\alpha_{1\ldots\er}|$, where 
$\,\er=\dim\bs$, and so 
$\,\int_\bs(\text{\rm div}{}^{\hs\text{\rm D}}\nh w)\,\alpha\,$ equals the 
Lebesgue integral of the Euclidean divergence 
$\,\partial_j(w^j|\alpha_{1\ldots\er}|)$, which vanishes since so do the 
integrals of the individual terms, for each fixed $\,j$.
\end{remark}

\section{The $\,2$-sphere and projective plane}\label{tspp}
For the $\,2$-sphere $\,S^2$ we have the following classification result, 
cf.\ \cite[p.\ 91]{pinkall-schwenk-schellschmidt-simon}.
\begin{corollary}\label{spher}The assignment described at the end of 
Theorem\/ $\ref{lccls}$ establishes a bijective correspondence between the 
equivalence classes of
\begin{enumerate}
  \def\theenumi{{\rm\roman{enumi}}}
\item[(i)] \pftc s $\,\hs\text{\rm D}\hh$ with \sy\ \rt s on the $\,2$-sphere 
$\,\bs$,
\end{enumerate}
and those of
\begin{enumerate}
  \def\theenumi{{\rm\roman{enumi}}}
\item[(ii)] $2$-spheres\/ $\,S\hh$ embedded in a fixed\/ $\,3$\diml\ \rvs\/ 
$\,\hy\hs$ and transverse to all lines in $\,\hy\hs$ passing through $\,0$.
\end{enumerate}
The equivalence relations are provided by\/{\rm:} affine diffeomorphisms 
between connections on $\,\bs\,$ in {\rm(i)}, linear isomorphisms\/ 
$\,\hy\to\hy\hs$ in {\rm(ii)}.
\end{corollary}
\begin{proof}This is immediate from Theorem~\ref{lccls} and Remark~\ref{cptim}.
\end{proof}
Connections $\,\hs\text{\rm D}\hh,\hs\deh\hs\,$ on \mf s $\,\bs,\bsh\,$ are 
called {\it projectively equivalent\/} if some diffeomorphism $\,\bs\to\bsh\,$ 
sends the geodesics of $\,\hs\text{\rm D}\hs\,$ onto (re-pa\-ram\-e\-trized) 
geodesics of $\,\hs\deh\hs$.

On the $\,2$-sphere $\,\bs\,$ there exists just one projective equivalence 
class of \pftc\ $\,\hs\text{\rm D}\hs\,$ with \sy\ \rt s (cf. 
\cite[Sect.\ 1.7]{pinkall-schwenk-schellschmidt-simon}). In fact, the 
diffeomorphism $\,\bs\to S\,$ described in Remark~\ref{cptim} for 
$\,\fe=\im\,$ and $\,\hy=\hs\text{\rm Ker}\,\lo\,$ (see Lemma~\ref{immer}) 
sends the geodesics of $\,\hs\text{\rm D}\hs\,$ (cf.\ Remark~\ref{ctraf}(a)) 
onto the great circles in $\,S$.
\begin{remark}\label{prpln}The assertion of Corollary~\ref{spher} remains true 
also when one replaces the phrase {\it$\,2$-sphere\/} in (i) with {\it real 
projective plane}, and adds to the end of (ii) the clause {\it as well as 
invariant under multiplication by $\,-\hh1\hh$}.

In fact, let $\,\varPsi:S^2\nh\to S^2$ be the fix\-ed-\pt\ free involution 
corresponding to the two-fold covering $\,S^2\nh\to\bbRP^2\nnh$, 
and let $\,\hs\text{\rm D}\hs\,$ be the pull\-back to $\,S^2$ of a given 
\pftc\ with a \sy\ \rt\ on $\,\bbRP^2\nnh$. Thus, 
$\,\varPsi^*\nh\alpha=-\hs\alpha\,$ for any fixed $\,\hs\text{\rm D}\hs$\prl\ 
area form $\,\alpha\,$ on $\,S^2\nnh$, since $\,\varPsi\,$ is an 
o\-ri\-en\-ta\-tion-re\-vers\-ing affine diffeomorphism (Remark~\ref{csunm}). 
From Remark~\ref{afdif} we now get $\,\fe(\varPsi(\cy))=-\hs\fe(\cy)\,$ for 
every $\,\cy\in S^2$ and the immersion $\,\fe\,$ of Lemma~\ref{immer}, 
while, by Remark~\ref{cptim}, $\,\fe\,$ is an embedding.
\end{remark}

\section{Projectively flat $\,2$-tori and Klein bottles}\label{pftt}
For any fixed real number $\,a\ne1$, the set 
$\,S=\{(z,|z|^a):z\in\bbC\smallsetminus\{0\}\}\,$ is a \su\ in the \rvs\ 
$\,\hy=\bbC\times\bbR\hs$, transverse to all lines containing $\,0$. Let 
$\,\hs\text{\rm D}\,$ denote the cen\-tro\-af\-fine connection on $\,S\,$ (see 
Section \ref{pfid}). The assignment $\,(z,|z|^a)\mapsto z\,$ is a 
diffeomorphism between $\,S\,$ and the multiplicative group 
$\,\bbC\smallsetminus\{0\}$, which turns $\,S\,$ into an Abelian Lie group. 
The group translations are affine diffeomorphisms $\,S\to S$, and, for 
$\,a=-\hs2$ only, they are all unimodular (Section \ref{eaco}): under our 
identification $\,\bbC\smallsetminus\{0\}\approx S$, the translation by 
$\,q\in\bbC\smallsetminus\{0\}\,$ corresponds to the restriction to $\,S\,$ of 
the linear automorphism $\,(z,t)\mapsto(qz,|q|^at)\,$ of 
$\,\hy=\bbC\times\bbR\hs$, which has the determinant $\,|q|^{a+2}$ (cf.\ 
Remarks~\ref{afdif} and~\ref{cafed}; note that automorphisms leaving $\,S\,$ 
invariant give rise to affine diffeomorphisms $\,S\to S$). Similarly, for any 
$\,q\in(0,\infty)$, the linear transformation 
$\,(z,t)\mapsto(q^{1/2}\hskip1pt\overline{\hskip-1ptz},q^{a/2}t)$ restricted 
to $\,S\,$ is an affine diffeomorphism, unimodular if $\,a=-\hs2$.
\begin{remark}\label{nosgm}Unless $\,a=0$, the above \su\ $\,S\,$ contains no 
nontrivial line segment. In fact, for a segment in 
$\,\bbC\smallsetminus\{0\}\,$ parametrized by $\,t\mapsto z=c\hs t+b$, where 
$\,b,c\,$ are nonzero complex numbers, $\,|z|^a$ cannot be a linear function 
of the parameter $\,t$, since by raising $\,|z|^a$ to the power $\,2/a\,$ one 
obtains the polynomial $\,|z|^2=|c\hs t+b\hs|^2$ of degree $\,2\,$ in $\,t$.
\end{remark}
The connections on the $\,2$-torus, discussed below, are among those described 
by Opozda \cite[Example~2.10]{opozda-04}. Our description is different, 
for reasons dictated by our applications.
\begin{example}\label{pftor}Any $\,q\in\bbC\,$ such that $\,0\ne|\hs q|\ne1\,$ 
obviously generates a multiplicative subgroup $\,G\,$ of 
$\,\bbC\smallsetminus\{0\}$, isomorphic to $\,\bbZ$, for which the quotient 
Lie group $\,(\bbC\smallsetminus\{0\})/G\,$ is a torus. For $\,a\ne1\,$ and 
$\,S\,$ as above, this fact, combined with the isomorphic identification 
$\,\bbC\smallsetminus\{0\}\approx S$, gives rise to a torus group 
$\,\bs=S/\bbZ$, where the action of $\,\bbZ\,$ on $\,S$ corresponds to 
$\,G\,$ acting on $\,\bbC\smallsetminus\{0\}$. The cen\-tro\-af\-fine 
connection $\,\hs\text{\rm D}\,$ on $\,S$, being \trinv, descends to a \trinv\ 
\pftc\ on $\,\bs$, also denoted by $\,\hs\text{\rm D}\hs$. Again, 
$\,\hs\text{\rm D}\hs\,$ is \ea\ on the torus $\,\bs\,$ if $\,a=-\hs2$. If, in 
addition, $\,q\in(1,\infty)$, the affine diffeomorphism 
$\,(z,t)\mapsto(q^{1/2}\hskip1pt\overline{\hskip-1ptz},q^{a/2}t)\,$ mentioned 
above, whose square is the generator of the $\,\bbZ\,$ action, commutes with 
the $\,\bbZ\,$ action on $\,S$, and so it descends a fix\-ed-\pt-free {\it 
affine involution\/} $\,\varPsi:\bs\to\bs$, while $\,\hs\text{\rm D}\,$ then 
gives rise to a \pftc, again denoted by $\,\hs\text{\rm D}\hs$, on the Klein 
bottle $\,\bs/\bbZ_2$, with the action of $\,\bbZ_2$ on $\,\bs\,$ generated by 
$\,\varPsi$. Finally, $\,\hs\text{\rm D}\,$ is not \rr\ unless $\,a=0$. 
Otherwise, Remarks~\ref{nosgm} and~\ref{nscgc} would contradict the 
Gauss-Bon\-net theorem.
\end{example}

\section{A differential equation}\label{adeq}
Let $\,\hs\text{\rm D}\hs\,$ be a \pftc\ on a surface $\,\bs\,$ such that the 
\rt\/ $\,\ri^{\hs\text{\rm D}}$ is \sy. We denote by $\,\op\,$ the 
sec\-ond-or\-der partial differential operator sending each 
twice-co\-var\-i\-ant \sy\ \tf\ $\,\tau\,$ on $\,\bs\,$ to a differential 
$\,2$-form $\,\op\tau$ on $\,\bs\,$ valued in $\,2$-forms, defined, in the 
notation of Remark~\ref{ddtau} and Section \ref{expr}, by
\begin{equation}
\op\tau\,\,=\,\,dd^{\hs\text{\rm D}}\nh\tau
\,+\,\tau\wedge\ri^{\hs\text{\rm D}}\nh.
\label{opf}
\end{equation}
The main result of this section, Theorem~\ref{slveq}(i), provides a simple 
topological condition sufficient for solvability of the linear equation 
$\,\op\tau=\nh A$, 
where the unknown $\,\tau\,$ is a twice-co\-var\-i\-ant \sy\ tensor on a 
\su\ $\,\bs\,$ with a connection $\,\hs\text{\rm D}\hs\,$ satisfying specific 
assumptions. Theorem~\ref{slveq}(ii) establishes the extent to which a 
solution $\,\tau\,$ is non\-u\-nique, by describing the solutions of the 
associated linear homogeneous equation. In local coordinates, the condition 
$\,\op\tau=A\,$ means that $\,A_{jklm}$ is equal to
\begin{equation}
\tau_{mk,lj}\nh-\nh\tau_{lk,mj}\nh-\nh\tau_{mj,lk}\nh+\nh\tau_{lj,mk}\nh
+\nh\tau_{mk}R_{lj}\nh-\nh\tau_{lk}R_{mj}\nh-\nh\tau_{mj}R_{lk}
+\nh\tau_{lj}R_{mk}\hs.
\label{lct}
\end{equation}
Here $\,A\,$ is a given four-times covariant tensor with 
$\,A_{kjlm}=-\hs A_{jklm}=A_{jkml}$. According to Theorem~\ref{slveq}(ii), 
when $\,\bs\,$ is \sco, (\ref{lct}) vanishes if and only if 
$\,\tau_{jk}\nh=\xi_{\hh k,j}\nh+\xi_{\hh j,k}$ for some $\,1$-form $\,\xi\,$ 
on $\,\bs$.

We begin with a lemma.
\begin{lemma}\label{clxct}Suppose that\/ $\,\hs\text{\rm D}\hs\,$ is a 
\pftc\ on a surface $\,\bs$, the \rt\/ $\,\ri^{\hs\text{\rm D}}$ of\/ 
$\,\hs\text{\rm D}\hs\,$ is \sy, $\,\lo\,$ and\/ $\,\op\,$ are the 
operators given by $(\ref{lxi})$ and\/ $(\ref{opf})$, while $\,\nabla$ denotes 
the flat connection in $\,\xz=(\tab)^{\wedge2}\nnh\oplus\tab$, defined in 
\hbox{Lemma\/ $\ref{fltco}$.}
\begin{enumerate}
  \def\theenumi{{\rm\alph{enumi}}}
\item For a fixed differential\/ $\,2$-form\/ $\,A\,$ on\/ $\,\bs\,$ valued in 
$\,2$-forms, the existence of a twice-co\-var\-i\-ant \sy\ \tf\/ $\,\tau\,$ 
on\/ $\,\bs\,$ with\/ $\,\op\tau=A\,$ is equivalent to 
$\,d^{\hs\nabla}\hskip-2.4pt$-exact\-ness of the 
$\,\xz$-val\-ued $\,2$-form\/ $\,(A,0)\,$ on $\,\bs$.
\item Any twice-co\-var\-i\-ant \sy\ \tf\/ $\,\tau\,$ on\/ $\,\bs\hs$ 
gives rise to the\/ $\,\xz$-val\-ued\/ $\,1$-form\/ 
$\,(d^{\hs\text{\rm D}}\nh\tau,\tau)$, and then
\begin{enumerate}
  \def\theenumi{{\rm\roman{enumi}}}
\item[i)] $\op\tau=0\,$ if and only if\/ 
$\,(d^{\hs\text{\rm D}}\nh\tau,\tau)\,$ is 
$\,d^{\hs\nabla}\hskip-2.4pt$-closed,
\item[ii)] $(d^{\hs\text{\rm D}}\nh\tau,\tau)\,$ is\/ 
$\,d^{\hs\nabla}\hskip-2.4pt$-exact if and only if\/ $\,\tau=\lo\xi\,$ 
for some real-val\-ued\/ $\,1$-form\/ $\,\xi$.
\end{enumerate}
\end{enumerate}
\end{lemma}
\begin{proof}By (\ref{dnp}), 
$\,d^{\hs\nabla}\hskip-2.4pt(\zeta,\gy)=(d\hh\zeta
+\gy\wedge\ri^{\hs\text{\rm D}}\nnh,\hh d^{\hs\text{\rm D}}\nh\gy-\zeta)\,$ 
whenever $\,\zeta\,$ is a $\,1$-form on $\,\bs\,$ valued in $\,2$-forms and 
$\,\gy\,$ is a twice-co\-var\-i\-ant tensor on $\,\bs\,$ viewed as a 
$\,\tab$-val\-ued $\,1$-form, cf.\ Remark~\ref{twice}, so that the pair 
$\,(\zeta,\gy)\,$ is an $\,\xz$-val\-ued $\,1$-form. (The term $\,-\hs\zeta\,$ 
in $\,d^{\hs\text{\rm D}}\nh\gy-\zeta$, treated here, in accordance with 
Remark~\ref{onetw}, as a $\,2$-form valued in $\,1$-forms, arises since 
summing $\,\zeta\hh$ cyclically over its arguments we get a real-val\-ued 
differential $\,3$-form on the \su\ $\,\bs$, that is, $\,0$.) Therefore, 
$\,d^{\hs\nabla}\hskip-2.4pt$-exact\-ness of $\,(A,0)\,$ in (a) means that 
$\,A=\op\gy\,$ for some twice-co\-var\-i\-ant tensor $\,\gy$, while, by 
(\ref{tgz}), $\,\op\gy\,$ remains unchanged when $\,\gy\,$ is replaced by its 
\sy\ part, which proves (a). Next, for $\,\tau\,$ as in (b), the above formula 
for $\,d^{\hs\nabla}\hskip-2.4pt(\zeta,\gy)$, applied to 
$\,\zeta=d^{\hs\text{\rm D}}\nh\tau\,$ and $\,\gy=\tau$, yields (i).

On the other hand, $\,d^{\hs\nabla}\hskip-2.4pt$-exact\-ness of 
$\,(d^{\hs\text{\rm D}}\nh\tau,\tau)\,$ means that 
$\,(d^{\hs\text{\rm D}}\nh\tau,\tau)=\nabla(\th,2\hh\xi)\,$ or, equivalently, 
$\,d^{\hs\text{\rm D}}\nh\tau=\,\text{\rm D}\hh\th\hh
-2\hh\xi\wedge\ri^{\hs\text{\rm D}}$ and 
$\,\tau=2\hs\text{\rm D}\hh\xi-\hs\th$, for some (real-val\-ued) 
$\,2$-form $\,\th\,$ and $\,1$-form $\,\xi$. Taking the transpose of the 
last equality (Remark~\ref{twice}) we get 
$\,\tau^*\nnh=2\hs(\hh\text{\rm D}\hh\xi)^*+\hs\th$, and so 
$\,\tau=(\tau+\tau^*)/2=\,\text{\rm D}\hs\xi\hs
+\hs(\hh\text{\rm D}\hs\xi)^*\nh=\lo\xi$.

Conversely, if $\,\tau=\hs\text{\rm D}\hs\xi\hs+\hs(\hh\text{\rm D}\hs\xi)^*$ 
for a $\,1$-form $\,\xi$, we get the $\,d^{\hs\nabla}\hskip-2.4pt$-exact\-ness 
relations $\,d^{\hs\text{\rm D}}\nnh\tau=\hs\text{\rm D}\hh\th\hh
-2\hh\xi\wedge\ri^{\hs\text{\rm D}}$ and 
$\,\tau=2\hs\text{\rm D}\hh\xi-\hs\th\,$ with the $\,2$-form 
$\,\th=\hs\text{\rm D}\hs\xi\hs-\hs(\hh\text{\rm D}\hs\xi)^*$ (that is, 
$\,\th=\hs d\hh\xi$). In fact, the second relation is obvious, and the first 
follows from (\ref{rid}) and (\ref{ris}).
\end{proof}
\begin{theorem}\label{slveq}Given a \pftc\/ $\,\hs\text{\rm D}\hs\,$ with a 
\sy\ \rt\ on a \sco\ \su\ $\,\bs$, let\/ $\,\lo,\op\,$ be as in\/ 
$(\ref{lxi})$ and\/ $(\ref{opf})$.
\begin{enumerate}
  \def\theenumi{{\rm\roman{enumi}}}
\item If\/ $\,\bs\,$ is noncompact and\/ $\,A\,$ is any differential\/ 
$\,2$-form on\/ $\,\bs\,$ valued in $\,2$-forms, then\/ 
$\,\op\tau=A\,$ for 
some twice-co\-var\-i\-ant \sy\ tensor\/ $\,\tau\hs$ on\/ $\,\bs$.
\item The twice-co\-var\-i\-ant \sy\ tensors\/ $\,\tau\hs$ with 
$\,\op\tau=0\,$ are precisely the tensors $\,\lo\xi$, for all 
real-val\-ued\/ $\,1$-forms\/ $\,\xi\hs$ on $\,\bs$.
\item A $\,1$-form $\,\xi\hs$ on $\,\bs\,$ is determined by the tensor\/ 
$\,\lo\xi\,$ uniquely up to adding an element of the $\,3$\diml\ \vs\ 
$\,\hs\text{\rm Ker}\,\lo$.
\end{enumerate}
\end{theorem}
\begin{proof}Assertion (i) is immediate from Lemma~\ref{clxct}(a): 
$\,\dim\bs=2$, so that the $\,(\tab)^{\wedge2}\nnh$-val\-ued $\,2$-form 
$\,(A,0)\,$ is $\,d^{\hs\nabla}\hskip-2.4pt$-closed, and hence 
$\,d^{\hs\nabla}\hskip-2.4pt$-exact, as $\,H^{2\nh}(\bs,\bbR)=\{0\}\,$ due to 
noncompactness of $\,\bs$. (Cf.\ the lines preceding Remark~\ref{twice}.) 
Similarly, (ii) follows from Lemma~\ref{clxct}(b), since $\,\bs\,$ is \sco, 
and so $\,d^{\hs\nabla}\hskip-2.4pt$-closed\-ness of 
$\,(d^{\hs\text{\rm D}}\nh\tau,\tau)\,$ is equivalent to its 
$\,d^{\hs\nabla}\hskip-2.4pt$-exact\-ness, while Lemma~\ref{immer}(i) 
gives (iii).
\end{proof}
\begin{remark}\label{landc}The differential operator $\,\op\,$ given by 
(\ref{opf}) (and (\ref{lct})), sending twice-co\-var\-i\-ant \sy\ \tf s to 
differential $\,2$-forms valued in $\,2$-forms, is well-defined on \mf s 
$\,\bs\,$ of any dimension $\,\er\ge2$, even though our rationale for writing 
$\,d\,$ instead of $\,d^{\hs\text{\rm D}}$ applies only when $\,\er=2$. The 
principal symbol of $\,-\hs\op/2\,$ equals that of the qua\-si-lin\-e\-ar 
operator sending a metric to its four-times covariant curvature tensor (cf.\ 
(\ref{rlm})). In some cases, the equality extends beyond the principal 
symbols; see (c) in Section \ref{riex}.
\end{remark}

\section{The case of closed surfaces}\label{tccs}
The construction in Section \ref{tlst} requires solvability of the equation
\begin{equation}
\op\tau\,\,=\,\,\ve\hs\alpha\otimes\alpha\hs,\hskip18pt\mathrm{with}\hskip9pt
\ve\in\{1,-1\}\hs,
\label{opt}
\end{equation}
where the unknown is a twice-co\-var\-i\-ant \sy\ \tf\ $\,\tau\,$ on a \su\ 
$\,\bs\,$ carrying an \ea\ \pftc\ $\,\hs\text{\rm D}\hs$, while 
$\,\pm\hs\alpha\,$ is a $\,\text{\rm D}\hs$\prl\ area element 
(Section \ref{eaco}), and $\,\op\,$ is given by (\ref{opf}). The value of 
$\,\ve\,$ is of no consequence for solvability of (\ref{opt}), since $\,\op\,$ 
is linear; in other words, solving (\ref{opt}) means, up to a factor, finding 
$\,\tau\,$ such that $\,\op\tau\,$ is parallel and nonzero. We have the 
following result.
\begin{theorem}\label{eqcls}On every closed \su\ $\,\bs$, {\rm(\ref{opt})} 
holds for some non-\rr, \ea, \pftc\ $\,\,\hs\text{\rm D}\hs$, some 
twice-co\-var\-i\-ant \sy\ \tf\/ $\,\tau$, and a $\,\text{\rm D}\hs$\prl\ area 
element\/ $\,\pm\hs\alpha$.
\end{theorem}
We precede the proof of Theorem~\ref{eqcls}, given in Section \ref{prft}, with 
three lemmas.

A fixed $\,\text{\rm D}\hs$\prl\ area element $\,\pm\hs\alpha\,$ on a \su\ 
$\,\bs\,$ with an \ea\ \tc\ $\,\hs\text{\rm D}\hs$, cf.\ 
Section \ref{eaco}, can be used to identify twice-co\-var\-i\-ant \sy\ tensors 
$\,\tau\,$ on $\,\bs\,$ (or, $\,2$-forms $\,A\,$ on $\,\bs\,$ valued in 
$\,2$-forms) with twice-con\-tra\-var\-i\-ant \sy\ tensors $\,T\,$ on 
$\,\bs\,$ (or, respectively, functions $\,\psi:\bs\to\bbR$). Specifically, 
such $\,\tau\,$ and $\,T\,$ are sections of the \vb s $\,[\tab]^{\odot2}$ and 
$\,[\tb]^{\odot2}\nnh$, while $\,\pm\hs\alpha\,$ gives rise to an isomorphism 
$\,\tb\to\tab\,$ acting on \vf s $\,v\,$ by $\,v\mapsto\alpha(v,\,\cdot\,)$, 
and hence defined at each \pt\ only up to a sign; however, the isomorphism 
$\,[\tb]^{\odot2}\nh\to[\tab]^{\odot2}$ induced by it, and sending $\,T\,$ to 
$\,\tau$, is single-val\-ued. In coordinates, 
$\,\tau_{jk}=\alpha_{jl}\hh\alpha_{km}T^{\hs lm}\nnh$. Similarly, our $\,A$ 
are sections of a line \bd\ in which $\,\alpha\otimes\alpha\,$ is a global 
single-val\-ued trivializing section, and so we can identify $\,A\,$ with 
$\,\psi\,$ such that $\,A=\hh\psi\hs\alpha\otimes\alpha$.
\begin{lemma}\label{opfid}Let\/ $\,\pm\hs\alpha\,$ be a 
$\,\text{\rm D}\hs$\prl\ area element on a \su\ $\,\bs\,$ with an \ea\ 
\pftc\ $\,\hs\text{\rm D}\hs\,$ having the \rt\ $\,\ri^{\hs\text{\rm D}}\nnh$. 
Under the above identifications, the operator $\,\op\,$ in {\rm(\ref{opf})} 
corresponds to an operator $\,\bop\,$ sending twice-con\-tra\-var\-i\-ant \sy\ 
\tf s $\,T\hs$ to functions $\,\bs\to\bbR\hs$. Explicitly, the operator 
$\,\bop\,$ is given by $\,\bop\hs T\,=\,\hs\text{\rm div}{}^{\hs\text{\rm D}}
(\text{\rm div}{}^{\hs\text{\rm D}}T)\,
+\,\lg\hh\ri^{\hs\text{\rm D}}\nnh,T\hh\rg$, where $\,\lr\,$ stands for the 
natural pairing between covariant and contravariant $\,2$-ten\-sors. In 
coordinates, $\,\bop\hs T=T^{jk}{}_{,\hh jk}+T^{jk}R_{jk}$.
\end{lemma}
\begin{proof}Let us fix a $\,\text{\rm D}\hs$\prl\ area element 
$\,\pm\hs\alpha\,$ on $\,\bs$, and let $\,\alpha^{-1}$ be the reciprocal 
bivector of $\,\alpha$, with the components $\,\alpha^{jk}$ such that 
$\,\alpha^{jl}\hh\alpha_{lk}=\delta^j_k$. Thus, $\,\alpha^{jk}\alpha_{jk}=-2$, 
and, whenever $\,A_{jklm}$ are the components of a $\,2$-form $\,A\,$ valued 
in $\,2$-forms, $\,4\psi=\alpha^{jk}\hh\alpha^{lm}A_{jklm}$ for the function 
$\,\psi\,$ such that $\,A=\psi\hs\alpha\otimes\alpha$. Now, if 
$\,\op\tau=\nh A$, the components $\,A_{jklm}$ are given by (\ref{lct}). 
Hence, for reasons of \sky, 
$\,4\psi=\alpha^{jk}\hh\alpha^{lm}(\tau_{mk,lj}\nh+\nh\tau_{mk}R_{lj})$, which 
equals $\,T^{jk}{}_{,\hh jk}+T^{jk}R_{jk}$ as $\,\alpha^{-1}$ is 
$\,\text{\rm D}\hs$\prl\ and $\,T^{jk}\nh
=\alpha^{jl}\hh\alpha^{km}\tau_{ml}$.
\end{proof}
\begin{remark}\label{fttor}In the notation of Lemma~\ref{opfid}, we have 
$\,\bop\hs T=1\,$ for the non-\rr\ \ea\ \trinv\ \pftc\ 
$\,\hs\text{\rm D}\hs\,$ on the $\,2$-torus $\,\bs\,$ described in 
Example~\ref{pftor} (with $\,a=-\hs2$) and a suitable \trinv\ 
twice-con\-tra\-var\-i\-ant \sy\ \tf\ $\,T\hs$ on $\,\bs$. In fact, it 
suffices to choose $\,T\hs$ so that the constant 
$\,\lg\hh\ri^{\hs\text{\rm D}}\nnh,T\hh\rg\,$ equals $\,1$, as the function 
$\,\hs\text{\rm div}{}^{\hs\text{\rm D}}
(\text{\rm div}{}^{\hs\text{\rm D}}T)$, being \trinv\ (that is, constant), 
must vanish by the divergence formula (Remark~\ref{divfo}).
\end{remark}
\begin{lemma}\label{trfru}Given a \pftc\ $\,\hs\text{\rm D}\hs\,$ on a \su\ 
$\,\bs$, a $\,\text{\rm D}\hs$\prl\ area element\/ $\,\pm\hs\alpha$, and a 
function $\,\fy:\bs\to(0,\infty)$, let\/ $\,\hs\deh\hs\,$ be the \pftc\ on 
$\,\bs\,$ with $\,\hs\deh\hs
=\hs\text{\rm D}\hs+2\hs\xi\odot\nh\text{\rm Id}\hs\,$ for 
$\,\xi=-\hs d\hh\log\fy$, cf.\ Lemma\/ $\ref{smgeo}$. The formula 
$\,\widetilde\alpha=\fy^{-3}\alpha\,$ then defines a 
$\,\hs\deh\hs$\prl\ area element\/ $\,\pm\hs\widetilde\alpha\,$ on 
$\,\bs$. Furthermore, if $\,\op\,$ and\/ $\,\bop\hs$ are the operators 
described in $(\ref{opf})$ and Lemma\/ $\ref{opfid}$, while 
$\,\widetilde\op,\widetilde\bop\,$ stand for their analogues corresponding 
to\/ $\,\hs\deh\hs\,$ and\/ $\,\widetilde\alpha$, then
\begin{enumerate}
  \def\theenumi{{\rm\alph{enumi}}}
\item $\widetilde{\op}(\fy^{-2}\tau)\,=\,\fy^{-2}\op\tau\,$ for any \sy\ 
covariant\/ $\,2$-ten\-sor\/ $\,\tau\,$ on\/ $\,\bs$, 
\item $\widetilde{\bop}(\fy^4T)\,=\,\fy^4\bop\hs T\,$ for any \sy\ 
contravariant\/ $\,2$-ten\-sor\/ $\,T$ on\/ $\,\bs$.
\end{enumerate}
\end{lemma}
\begin{proof}The $\,\hs\deh\hs$-divergence of any \vf\ 
$\,w\,$ (or, twice-con\-tra\-var\-i\-ant \sy\ \tf\ $\,T\hh$) clearly equals 
$\,\hs\text{\rm div}{}^{\hs\text{\rm D}}w+(\er+1)\hs\xi(w)\,$ (or, 
$\,\hs\text{\rm div}{}^{\hs\text{\rm D}}T\nh+(\er+3)T\xi$) whenever 
$\,\hs\text{\rm D}\hs\,$ and $\,\deh\hs\,$ are \tc s on 
a \mf\ $\,\bs\,$ of dimension $\,\er$ and $\,\hs\deh\hs
=\hs\text{\rm D}\hs+2\hs\xi\odot\nh\text{\rm Id}\hs$. Using 
Lemma~\ref{smgeo}(i), we now get $\,\widetilde{\bop}\hs T=\bop\hs T
+(\er^2\nh+5\er+2)\lg T,\xi\otimes\xi\rg+4\lg T,\hs\text{\rm D}\hs\xi\rg
+(2\er+4)\hs\xi(\text{\rm div}{}^{\hs\text{\rm D}}T)$. When $\,\er=2\,$ 
and $\,\xi=-\hs d\hh\log\fy$, this easily yields (b). As we then 
obviously have $\,\hs\deh\hs\widetilde\alpha=0$ 
for $\,\widetilde\alpha=\fy^{-3}\alpha$, assertion (a) follows.
\end{proof}
\begin{lemma}\label{surje}Let\/ $\,g\,$ a Riemannian metric of constant 
Gaussian curvature $\,K$ on a closed \su\ $\,\bs$, and let\/ $\,\bop\hs$ be 
the operator defined in Lemma\/ $\ref{opfid}$ with $\,\hs\text{\rm D}\hs$ 
replaced by the \lcc\ $\,\nabla$ of $\,g$. If\/ $\,K<0$, or $\,\bs\hs$ is 
nonorientable and\/ $\,K>0$, then\/ $\,\bop\,$ is a surjective operator from 
the space of\/ $\,C^\infty$ twice-con\-tra\-var\-i\-ant \sy\ \tf s on\/ 
$\,\bs\hs$ onto the space of\/ $\,C^\infty$ functions $\,\bs\to\bbR\hs$.
\end{lemma}
\begin{proof}Let $\,g^{-1}$ be the reciprocal metric, with the components 
$\,g^{jk}\nnh$. The operator $\,f\mapsto\bop(fg^{-1})\,$ clearly equals 
$\,\Delta+2K$, that is, sends any $\,C^\infty$ function $\,f\,$ to 
$\,\Delta f+2Kf$, where $\,\Delta\,$ is the Laplacian of $\,g$, acting by 
$\,\Delta f=g^{jk}f_{,\hh jk}$. Since $\,\Delta+2K\,$ is self-ad\-joint and 
elliptic, its surjectivity will be immediate once we establish its injectivity.

Injectivity of $\,\Delta+2K\,$ is clear when $\,K<0$, as $\,\Delta+2K\,$ then 
is a negative operator. Suppose now that $\,\bs\,$ is nonorientable and 
$\,K>0$. The two-fold covering of $\,\bs\,$ is the round sphere $\,S^2$ of 
curvature $\,K$, on which $\,2K\,$ is the lowest positive eigenvalue of 
$\,-\Delta$, and the corresponding eigenspace 
$\,\hs\text{\rm Ker}\,(\Delta+2K)\,$ consists of restrictions of linear 
functionals on $\,\rtr$ to $\,S^2$ treated as a sphere in $\,\rtr$ centered at 
$\,0$. Eigenfunctions thus cannot descend to the projective plane $\,\bs$, as 
they change sign under the antipodal involution, and injectivity of 
$\,\Delta+2K\,$ on $\,\bs\,$ follows. This completes the proof.
\end{proof}
Note that $\,\bop\,$ is not surjective in the remaining cases: if $\,K=0$, or 
$\,\bs\,$ is orientable and $\,K\,$ is a positive constant, the image of 
$\,\bop\,$ is the $\,L^2$-or\-thog\-o\-nal complement of the kernel of 
$\,\Delta+2K$. In fact, 
$\,\hs\text{\rm Ker}\,(\Delta+2K)=\hs\text{\rm Ker}\,\bop^{\hs*}$ for the 
adjoint $\,\bop^{\hs*}$ of $\,\bop$, given by 
$\,\bop^{\hs*}\psi=\nabla d\psi+K\psi g\,$ for any function $\,\psi$.

\section{Proof of Theorem~\ref{eqcls}}\label{prft}
To prove Theorem~\ref{eqcls}, it suffices, by Lemma~\ref{opfid}, to exhibit a 
connection $\,\hs\text{\rm D}\hs\,$ on $\,\bs$ with the stated properties 
and a twice-con\-tra\-var\-i\-ant \sy\ \tf\ $\,T\hh$ on $\,\bs\,$ such that 
$\,\bop\hs T=1$, i.e., $\,\hs\text{\rm div}{}^{\hs\text{\rm D}}
(\text{\rm div}{}^{\hs\text{\rm D}}T)
+\lg\hh\ri^{\hs\text{\rm D}}\nnh,T\hh\rg=1$.

When $\,\bs\,$ is diffeomorphic to the $\,2$-torus, $\,\hs\text{\rm D}\hs\,$ 
and $\,T\hs$ exist according to Remark~\ref{fttor}.

If our \su\ is diffeomorphic to the Klein bottle, we may denote it by 
$\,\bs/\bbZ_2$ (rather than $\,\bs$), and let $\,\bs,\hs\text{\rm D}\hs,T\,$ 
be as in Remark~\ref{fttor}, the $\,\bbZ_2$ action on the torus $\,\bs\,$ 
being generated by the affine involution $\,\varPsi\,$ of Example~\ref{pftor}. 
As $\,\bop\hs T=\bop\hs T\hh'\nnh=1$, where $\,T\hh'$ is the push-for\-ward of 
$\,T\,$ under $\,\varPsi$, the $\,\varPsi$\inv\ \tf\ $\,(T+T\hh'\hh)/2\,$ on 
$\,\bs\,$ descends to the required \tf\ on the Klein bottle $\,\bs/\bbZ_2$.

If $\,\bs\,$ is not diffeomorphic to the $\,2$-torus, or the Klein bottle, or 
the $\,2$-sphere, we may choose on $\,\bs\,$ a Riemannian metric $\,g\,$ 
of constant Gaussian curvature and denote by $\,\nabla\,$ its \lcc. The 
required connection $\,\hs\text{\rm D}\hs\,$ on $\,\bs\,$ can now be obtained 
from $\,\nabla$ by a projective modification. Namely, we set 
$\,\hs\text{\rm D}\hs=\nabla+\hs2\hs\xi\odot\nh\text{\rm Id}\hs\,$ (notation 
of Lemma~\ref{smgeo}), with $\,\xi=-\hs d\hh\log\fy\,$ for a suitable function 
$\,\fy:\bs\to(0,\infty)$. In view of Lemma~\ref{trfru}(a) with 
$\,\hs\deh\hs,\hs\text{\rm D}\hs\,$ replaced by 
$\,\hs\text{\rm D}\hs,\nabla$, the operator $\,\bop\,$ corresponding to 
$\,\hs\text{\rm D}\hs\,$ is surjective since so is, by Lemma~\ref{surje}, 
the analogue of $\,\bop\,$ corresponding to $\,\nabla$.

We still need to verify here that, for a suitable choice of $\,\fy$, the 
connection $\,\hs\text{\rm D}\hs\,$ will not be \rr. Although this might be 
justified by very general reasons (namely, Ric\-ci-\-re\-cur\-rence of 
$\,\hs\text{\rm D}\hs\,$ would amount to imposing on $\,\fy\,$ a system of 
\pde s), a direct geometric argument is also possible. Specifically, given a 
nonempty contractible open set $\,\,U\subset\bs\,$ for which the immersion 
$\,\fe:U\to\hs\text{\rm Ker}\,{\lo}\,$ of Lemma~\ref{immer} is an embedding, a 
compactly supported small deformation of the image $\,\fe(\bs)\,$ 
yields a \su\ $\,\widehat S\,$ in $\,\hs\text{\rm Ker}\,{\lo}$, again 
transverse to lines through $\,0$, and such that the radial projection 
$\,\pmb:\fe(U)\to\widehat S\,$ is a diffeomorphism; $\,\pmb\,$ sends 
$\,\fe$-im\-ages of geodesics in $\,\,U\,$ onto (re-pa\-ram\-e\-trized) 
geodesics of the cen\-tro\-af\-fine connection $\,\hs\deh\hs\,$ on 
$\,\widehat S$. (See Remark~\ref{ctraf}(a).) Let the connection 
$\,\hs\text{\rm D}\hs\,$ on $\,\bs\,$ now be the result of replacing 
$\,\nabla\nh$, \hskip.7ptjust on $\,\,U\nnh$, with the pull\-back of 
$\,\hs\deh\,$ under the composite $\,\pmb\circ\fe$. By Lemma~\ref{smgeo}, 
$\,\hs\text{\rm D}\hs=\nabla+2\hs\xi\odot\nh\text{\rm Id}\hs\,$ with 
$\,\xi=-\hs d\hh\log\fy\,$ for some function $\,\fy:\bs\to(0,\infty)\,$ 
equal to $\,1\,$ outside a compact subset of $\,\,U\nnh$. Choosing $\,S\,$ so 
that some nonempty open subset of $\,S\,$ is not contained in a quadric 
surface, and at the same time contains no nontrivial line segment, we now 
ensure that $\,\hs\text{\rm D}$ is not \rr. (Cf. Theorems~\ref{lcsym} 
and~\ref{rcrec}.)

Finally, the existence of the required $\,\hs\text{\rm D}\,$ and $\,T\hh$ on 
the $\,2$-sphere is immediate from their existence on the projective plane.

\newpage
\part{NULL PARALLEL DISTRIBUTIONS}\label{npd}
A general lo\-cal-co\-or\-di\-nate form of \prc s with \npd s was found by 
Walker \cite{walker}; for a co\-or\-di\-nate-free version, see 
\cite{derdzinski-roter-06}. In this part we discuss a certain class of \npd s 
without using Walker's theorem directly. However, our ultimate application (in 
Part III) of the results obtained here does lead to a special type of Walker 
coordinates; see formulae (\ref{ged}) in Section \ref{pfth}.

\section{Curvature conditions}\label{afcc}
Let $\,\vt\,$ be a \npd\ of dimension $\,\er\hs$ on an $\,n$\diml\ \prd\
$\,(M,g)$. Thus, the $\,g$-or\-thog\-o\-nal complement $\,\w$ is a parallel 
\dn\ of dimension $\,n-\er$. If the sign pattern of $\,g\,$ has $\,i_-$ 
minuses and $\,i_+$ pluses, then
\begin{equation}
\text{\rm(i)}\hskip9pt\er\,\le\,\,\text{\rm min}\hs(i_-,i_+)\hs,
\hskip24pt\text{\rm(ii)}\hskip9pt\vt\,\subset\,\w\nh,
\hskip24pt\text{\rm(iii)}\hskip9pt\er\,\le\,n/2\hs.\label{inq}
\end{equation}
In fact, $\,\vt\,$ is null, which gives (\ref{inq}\hh-ii) and hence
$\,\er\le n-\er$, that is, \hbox{(\ref{inq}\hh-iii).} Now (\ref{inq}\hh-i) 
follows: in a pseu\-\hbox{do\hs-}Euclid\-e\-an space with the sign pattern as 
above, $\,i_\pm$ is the maximum dimension of a subspace on which the inner 
product is positive/negative semidefinite.

Every \npd\ $\,\vt\,$ satisfies the curvature conditions
\begin{equation}
\text{\rm(a)}\hskip9ptR(\vt,\w\nnh,\,\cdot\,,\,\cdot\,)\,
=\,0\hs,\hskip11pt
\text{\rm(b)}\hskip9ptR(\vt,\vt,\,\cdot\,,\,\cdot\,)\,
=\,0\hs,\hskip11pt
\text{\rm(c)}
\hskip9ptR(\w\nnh,\w\nnh,\vt,\,\cdot\,)\,=\,0\hs,
\label{rzo}
\end{equation}
where (\ref{rzo}\hh-a) states that $\,R(v,u,w,w')=0\,$ for all \vf s 
$\,v,u,w,w'$ such that $\,v\,$ is a section of $\,\vt\,$ and $\,u\,$ is a 
section of $\,\w\nnh$, and similarly for (\ref{rzo}\hh-b), (\ref{rzo}\hh-c). 
In fact, given such $\,v,u,w,w'\nnh$, (\ref{cur}) implies that $\,R(w,w')v\,$ 
is a section of $\,\vt$, and so it is orthogonal to $\,u$. This gives 
(\ref{rzo}\hh-a), while (\ref{rzo}\hh-a) and (\ref{inq}\hh-ii) yield 
(\ref{rzo}\hh-b). Finally, (\ref{rzo}\hh-a) and the first Bianchi identity 
imply (\ref{rzo}\hh-c).

We will focus our discussion on the case where, in addition to (\ref{rzo}), 
the conditions
\begin{equation}\text{\rm(i)}\hskip9pt
R(\vt,\,\cdot\,,\w\hskip-2.4pt,\,\cdot\,)\,=\,0\hs,
\hskip38pt\text{\rm(ii)}\hskip9pt
R(\w\hskip-2.8pt,\w\hskip-2.8pt,\,\cdot\,,\,\cdot\,)\,=\,0
\label{rvv}
\end{equation}
hold for a \npd\ $\,\vt\,$ on a \prd.

Our interest in (\ref{rvv}) arises from the fact that a $\,2$\diml\ \npd\ 
$\,\vt\,$ with (\ref{rvv}) exists on every \cs\ \prd\ $\,(M,g)$ with $\rwo$. 
(See Lemma~\ref{precs}(ii).)
\begin{remark}\label{liebr}A \vf\ $\,w\,$ on the \ts\ $\,M\,$ of a \bd\ is 
$\,\pmb$\pj\ onto the \bmf\ $\,\bs$, where $\,\pmb:M\to\bs\,$ is the \bp, if 
and only if, for every vertical \vf\ $\,u\,$ on $\,M$, the Lie bracket 
$\,[w,u]\,$ is also vertical. This is easily verified in local coordinates for 
$\,M\,$ that make $\,\pmb\,$ appear as a Euclidean projection.
\end{remark}

\section{Pro\-ject\-a\-bil\-i\-ty of the Le\-vi-Ci\-vi\-ta 
connection}\label{plcc}
We will use the following assumption. The clause about the leaves of 
$\,\w$ follows from the other hypotheses if one replaces $\,M\,$ by a suitable 
neighborhood of any given \pt.
\begin{equation}
\begin{array}{l}
\text{\rm $\vt\,$ is an $\,\er$\diml\ \npd\ on a 
pseu\-\hbox{do\hs-}Riem\-ann\-i\-an mani-}\\
\text{\rm fold $\,\hskip1.4pt(M,g)\hs\,$ with $\,\hs\dim M=\hh n\hh$, and 
$\,\hs\w$ has contractible leaves which are the f\hs i-}\\
\text{\rm bres of a \bp\ $\hs\pmb\nnh:\nnh M\nnh\to\nnh\bs\hs$ for some 
$\er$\diml\ \bmf\ $\hh\bs$.}
\end{array}\label{cnd}
\end{equation}
Given a \npd\ $\,\vt\hs$ on a \prd\ $\,(M,g)$, we say that the \lcc\ 
$\,\nabla\hs$ of $\,g\,$ is {\it$\,\w\nh$-pro\-ject\-a\-ble}, cf.\ 
\cite{ghanam-thompson}, if, for $\,\pmb,\bs\,$ chosen, locally, as in 
(\ref{cnd}), and for any \vf s $\,w,v\,$ in $\,M\,$ such that $\,w\,$ is 
$\,\pmb\hs$\pj\ and $\,v\,$ is a section of $\,\vt\hs$ parallel along 
$\,\w\nnh$, the section $\,\nabla_{\!w}v\,$ of $\,\vt\hs$ is parallel along 
$\,\w\nh$ as well.
\begin{lemma}\label{rvveq}A \npd\/ $\,\vt\hs$ on a \prd\/ $\,(M,g)$ 
satisfies condition\/ $(\ref{rvv}${\rm\hh-i}$)$ if and only if the \lcc\ 
$\,\nabla$ is $\,\w\nh$\pj.
\end{lemma}
\begin{proof}For $\,w,v\,$ as above and any section $\,u\,$ of $\,\w\nnh$, 
(\ref{cur}) gives $\,\nabla_{\!u}\nabla_{\!w}v=R(w,u)v$, as the other two 
terms in (\ref{cur}) vanish: $\,v\,$ is parallel along $\,\w\nnh$, and so 
$\,\nabla_{\!u}v=\nabla_{[w,u]}v=0$, where $\,[w,u]\,$ is a section of 
$\,\w=\,\kerd\pmb\,$ due to $\,\pmb$-pro\-ject\-a\-bil\-i\-ty of $\,w\,$ and 
Remark~\ref{liebr}.
\end{proof}
Let us now assume (\ref{cnd}), and let the \lcc\ $\,\nabla\hskip1pt$ be 
$\,\w\nh$\pj.
\begin{enumerate}
  \def\theenumi{{\rm\roman{enumi}}}
\item For any $\,\pmb\hs$\pj\ \vf\ $\,w$, if $\,u\,$ is a section of 
$\,\w\nnh$, then so is $\nabla_{\!u}w$.
\item Sections $\,v\,$ of $\,\vt\,$ parallel along $\,\w$ are in a natural 
bijective correspondence $\,\vl$ with sections $\,\xi\,$ of $\,\tab$. It 
assigns to $\,v\,$ the $\,1$-form $\,\xi=\vl(v)\,$ on $\,\bs\,$ 
such that $\,\xi((d\pmb)w)=g(v,w):\bs\to\bbR\,$ for any $\,\pmb\hs$\pj\ \vf\ 
$\,w\,$ on $\,M$. (Here and below $\,(d\pmb)w\,$ denotes the \vf\ on $\,\bs\,$ 
onto which $\,w\,$ projects.)
\item For any $\,\phi:\bs\to\bbR\,$ treated as a function on $\,M\,$ constant 
along $\,\w\nnh$, the $\,g$-gra\-di\-ent $\,v=\nabla\nh\phi\,$ is a section of 
$\,\vt\,$ parallel along $\,\w$ and $\,\vl(v)=\xi\,$ for the section 
$\,\xi=d\phi\,$ of $\,\tab$.
\item $\nabla\,$ gives rise to a \tc\ $\,\hs\text{\rm D}\hs\,$ in the \tbd\ 
$\,\tb$, which we call the {\it$\,\w\nh$\pd\ connection\/} for $\,g\,$ and 
$\,\vt$, and which is characterized by 
$\,\text{\rm D}_{\hh(d\pmb)w}\hs\xi=\vl(\nabla_{\!w}v)\,$ whenever 
$\,\xi=\vl(v)$, for $\,\pmb\hs$\pj\ \vf s $\,w\,$ on $\,M$ and sections 
$\,v\,$ of $\,\vt\,$ parallel along $\,\w\nnh$.
\item $\vl(R(w,w\hs'\hh)\hh v)
=-\hs\xi\hs R^{\text{\rm D}}((d\pmb)w,(d\pmb)w\hs'\hh)\,$ if $\,\xi=\vl(v)$, 
with $\,\xi\hs R^{\text{\rm D}}$ as in (\ref{rid}), for sections $\,v\,$ of 
$\,\vt\,$ parallel along $\,\w\nnh$, the curvature tensors $\,R\,$ of $\,g\,$ 
and $\,R^{\text{\rm D}}$ of the $\,\w\nh$\pd\ connection 
$\,\hs\text{\rm D}\hs\,$ on $\,\bs$, and $\,\pmb\hs$\pj\ \vf s $\,w,w\hs'$ on 
$\,M$.
\item $\vt\,$ is $\,\widetilde g$\prl\ for the conformally related metric 
$\,\widetilde g=\fy^{-2}g\,$ on $\,M$, whenever $\,\fy:\bs\to(0,\infty)$, so 
that $\,\fy\,$ may be treated as a function on $\,M\,$ constant along 
$\,\w\nnh$. The \lcc\ $\,\widetilde\nabla\,$ of $\,\widetilde g\,$ then is 
$\,\w\nh$\pj, and its $\,\w\nh$\pd\ connection is 
$\,\hs\deh\hs=\hs\text{\rm D}\hs+2\hs d\phi\hs\odot\text{\rm Id}\hs\,$ with 
$\,\phi=-\hs\log\fy\,$ (notation of Lemma~\ref{smgeo}).
\end{enumerate}
In fact, (i) follows as $\,\nabla_{\!u}w=[u,w]+\nabla_{\!w}u$, while 
$\,[u,w]\,$ (or $\,\nabla_{\!w}u$) is a section of $\,\w$ by 
Remark~\ref{liebr} (or, respectively, since $\,\w$ is parallel). Next, 
$\,\xi\,$ in (ii) is well defined: two choices of $\,w\,$ having the same 
$\,(d\pmb)w\,$ differ by a section of $\,\w$ (orthogonal to $\,v$). Also, the 
function $\,g(v,w):M\to\bbR\,$ descends to $\,\bs$. Namely, 
$\,d_u[g(v,w)]=0\,$ for any section $\,u$ of $\,\w$, as 
$\,\nabla_{\!u}v=0\,$ and, by (i), $\,g(v,\nabla_{\!u}w)=0$. Injectivity of 
$\,\vl\,$ is obvious; its surjectivity easily follows since, by 
(\ref{rzo}\hh-c) and (\ref{cur}), the connections induced by $\,\nabla\,$ in 
the restrictions of $\,\vt$ to the leaves of $\,\w$ are all flat. This proves 
(ii). Next, (iii) is obvious from (ii), as $\,v=\nabla\nh\phi$, being orthogonal 
to $\,\w\nnh$, is a section of $\,\vt$, and hence so is $\,\nabla_{\!w}v\,$ 
for any \vf\ $\,w$, which gives $\,g(\nabla_{\!u}v,w)=g(\nabla_{\!w}v,u)=0\,$ 
for sections $\,u\,$ of $\,\w\nnh$. In \hbox{(iv), $\,\,\text{\rm D}$} is 
clearly well defined, and it is tor\-sion\-free: the second covariant 
derivative $\,\hs\text{\rm D}\hs d\phi\,$ of any function 
$\,\phi:\bs\to\bbR\,$ is symmetric, as (iii) yields 
$\,(\hh\text{\rm D}\hs d\phi)((d\pmb)w,(d\pmb)w\hs'\hh)
=g(\nabla_{\!w}v,w\hs'\hh)=(\nabla\nh d\phi)(w,w\hs'\hh)$ for 
$\,v=\nabla\nh\phi\,$ and $\,\pmb\hs$\pj\ \vf s $\,w,w\hs'$ on $\,M$. Finally, 
(v) follows from (iv), (\ref{cur}), (\ref{dnp}) and (\ref{rid}\hh-i), while 
(a) in Section \ref{expr}, (iii) and (i) imply (vi).

\section{Pull\-backs of covariant tensors}\label{poct}
We will say that a $\,\hs k$-times covariant \tf\ $\,\gy\,$ on a \mf\ 
$\,M\,$ {\it annihilates\/} a \dn\ $\,\vt\,$ on $\,M\,$ if $\,k\ge1\,$ and 
$\,\gy(v_1,\dots,v_k)=0\,$ whenever $\,v_1,\dots,v_k\,$ are \vf s on $\,M\,$ 
and $\,v_j$ is a section of $\,\vt\hs$ for some $\,j\in\{1,\dots,k\}$.
\begin{lemma}\label{pullb}Suppose that we have\/ $(\ref{cnd})$ and\/ 
$(\ref{rvv}${\rm\hh-i}$)$, $\,\hs\text{\rm D}\hs\,$ is the $\,\w\nh$\pd\ 
connection on $\,\bs$, cf.\ Section $\ref{plcc}$, and\/ $\,\tau\hs$ is a 
$\,k$-times covariant \tf\ on $\,\bs$.
\begin{enumerate}
  \def\theenumi{{\rm\alph{enumi}}}
\item If\/ $\,k\ge1$, the pull\-back $\,\pmb^*\tau\hs$ of\/ $\,\tau\,$ to 
$\,M\hs$ annihilates $\,\w\nnh$.
\item For $\,k\ge0$, we have 
$\,\nabla(\pmb^*\tau)=\pmb^*(\hs\text{\rm D}\hs\tau)$, with both covariant 
derivatives treated as $\,(k+1)$-times covariant \tf s.
\end{enumerate}
Also, if\/ $\,\gy\,$ is a $\,k$-times covariant \tf\ on $\,M\hs$ and\/ 
$\,k\ge1$, while\/ $\,\gy\,$ and\/ $\,\nabla\gy\,$ both annihilate $\,\w\nnh$, 
then $\,\gy=\pmb^*\tau\,$ for a unique $\,k$-times covariant \tf\/ $\,\tau\hs$ 
on $\,\bs$.
\end{lemma}
\begin{proof}Since $\,\w=\,\kerd\pmb$, (a) follows. To obtain (b), we need 
only to consider the cases $\,k=0\,$ and $\,k=1$, as the class of covariant 
tensor fields $\,\tau\,$ on $\,\bs\,$ with 
$\,\nabla(\pmb^*\tau)=\pmb^*(\hs\text{\rm D}\hs\tau)\,$ is closed under both 
addition (with any fixed $\,k$) and tensor multiplication. For $\,k=0$, (b) is 
obvious: on functions, $\,\nabla=\hs\text{\rm D}\hs=\hs d$. If $\,k=1$, so 
that $\,\tau\,$ is a section of $\,\tab$, (ii) in Section \ref{plcc} gives 
$\,\pmb^*\tau=g(v,\,\cdot\,)\,$ for a section $\,v\,$ of $\,\vt\,$ parallel 
along $\,\w\nnh$, and then 
$\,\nabla(\pmb^*\tau)=\pmb^*(\hs\text{\rm D}\hs\tau)\,$ by (iv) in Section 
\ref{plcc}. This proves (b). In the final clause, we set 
$\,\tau((d\pmb)w_1,\dots,(d\pmb)w_k)=\gy(w_1,\dots,w_k)$, for $\,\pmb\hs$\pj\ 
\vf s  $\,w_j$ on $\,M$. That $\,d_u[\hh\gy(w_1,\dots,w_k)]=0\,$ for any 
section $\,u\,$ of $\,\w$ now follows since $\,\nabla\gy\,$ annihilates $\,\w$ 
(and so $\,\nabla_{\!u}\gy=0$), while each $\,\nabla_{\!u}w_j$ is a section of 
$\,\w$ by (i) in Section \ref{plcc}.
\end{proof}

\section{Riemann extensions}\label{riex}
Let $\,M=\tab\,$ be the \ts\ of the \ctb\ of a \mf\ $\,\bs$. Our convention is 
that, as a set, 
$\,\tab=\{(\cy,\eta):\cy\in\bs\hskip5pt\mathrm{and}\hskip5pt\eta\in\tcb\}$. 
Any fixed connection $\,\hs\text{\rm D}\hs\,$ on $\,\bs\,$ gives rise to the 
\prc\ $\,h^{\text{\rm D}}$ on $\,\tab\,$ characterized by
\begin{equation}
h_x^{\text{\rm D}}(w,w)\,=\,2\hs w\vrt(d\pmb_xw\hrz)\hskip8pt\mathrm{for\ }
\,x=(\cy,\eta)\in M=\tab\,\mathrm{\ and\ }\,w\in\txm,
\label{rex}
\end{equation}
$\,w\vrt\nnh,w\hrz$ being the vertical and 
$\,\hs\text{\rm D}\hs$-hor\-i\-zon\-tal components of $\,w$, with 
$\,w\vrt\in\tacb\,$ due to the identification between the vertical space at 
$\,x\,$ and the fibre $\,\tacb$. In other words, all vertical and all 
$\,\hs\text{\rm D}\hs$-hor\-i\-zon\-tal vectors are 
$\,h^{\text{\rm D}}\nnh$-null, while 
$\,h_x^{\text{\rm D}}(\zeta,w)=\zeta(d\pmb_xw)$, for 
$\,x=(\cy,\eta)\in\tab=M$, a vertical vector $\,\zeta\in\kerd\pmb_x=\tacb$, 
and any vector $\,w\in\txm$. Here $\,\pmb:\tab\to\bs\,$ is the \bp.

Our $\,h^{\text{\rm D}}$ is one of Patterson and Walker's {\it Riemann 
extension metrics\/} \cite{patterson-walker}. In the coordinates 
$\,\y^{\hs j}\nnh,\p_j$ for $\,\tab\,$ obtained from an arbitrary local 
coordinate system $\,\y^{\hs j}$ in $\,\bs$,
\begin{equation}
h^{\text{\rm D}}=\,2\hs dp_j\,d\y^{\hs j}
-\hs2\p_j\vg_{\!kl}^j\,d\y^k\hs d\y^{\hh l}\hh,
\label{hde}
\end{equation}
where the products of differentials stand for \sy\ products, and 
$\,\vg_{\!jl}^k$ are the components of $\,\hs\text{\rm D}\hs\,$ in the 
coordinates $\,\y^{\hs j}\nnh$. Namely, due to the coordinate expression for 
the covariant derivative, the horizontal lift of a \vf\ $\,w\,$ tangent to 
$\,\bs\,$ has the components $\,\dot\y^{\hs j}=w^j$ and 
$\,\dot\p_l=\vg_{\!kl}^jw^k\p_j$ in the coordinates 
$\,\y^{\hs j}\nnh,\p_j$. Thus, (\ref{hde}) implies (\ref{rex}).

We will not use the easily-verified facts that the vertical distribution 
$\,\vt\,$ on $\,\tab\,$ is $\,h^{\text{\rm D}}\nnh$\prl, its \lcc\ is 
$\,\w\nh$\pj\ in the sense of Section \ref{plcc}, and the $\,\w\nh$\pd\ 
connection on $\,\bs\,$ coincides with the original $\,\hs\text{\rm D}\hs$.

Another interesting example \cite{walker} arises when $\,\bs\,$ is a \mf\ of 
dimension $\,\er\ge2$ with a global coordinate system $\,y^{\hs j}$ and with 
a fixed twice-co\-var\-i\-ant \sy\ \tf\ $\,\gy$, while $\,\hs\deh\hs\,$ is the 
flat \tc\ on $\,\bs\,$ with component functions in the coordinates 
$\,\y^{\hs j}$ all equal to zero. On $\,N\nh=\tab\,$ we have the \prc
\begin{equation}
\heh\,=\,2\hs d\q_j\hs d\y^{\hs j}\hs
-\,2\gy_{kl}\,d\y^k\hs d\y^{\hh l}\nh,\hskip5pt\text{\rm with $\,\gy_{kl}$ 
not depending on the coordinates $\,\q_j$,}
\label{tih}
\end{equation}
in the standard coordinate system for $\,N\nh=\tab$, this time denoted by 
$\,\y^{\hs j}\nnh,\q_j$. Then
\begin{enumerate}
  \def\theenumi{{\rm\alph{enumi}}}
\item the vertical sub\bd\ $\,\vt\,$ of $\,\tn\hs$ is a \npd\ on the \prd\ 
$\,(N,\heh)$, and it is spanned by $\,\heh$\prl\ \vf s,
\item the \lcc\ $\,\widetilde\nabla\,$ of $\,\heh\,$ is $\,\w\nh$\pj\ 
(cf.\  Section \ref{plcc}), and the corresponding $\,\w\nh$\pd\ connection is 
$\,\hs\deh\hs$,
\item $\heh\,$ is Ric\-ci-flat, while its four-times covariant 
curvature and Weyl tensors are given by 
$\,\widetilde R=\widetilde W\nnh=\pmb^*(\widetilde\op\gy)$, with 
$\,\widetilde\op\,$ as in Remark~\ref{landc} (for $\,\hs\deh\hs\,$ rather 
than $\,\hs\text{\rm D}$).
\end{enumerate}
(See also Lemma~\ref{walkr}.) As usual, $\,\pmb:N\to\bs\,$ is the \bp.

In fact, with the index ranges $\,\hs j,k,l,m\in\{1,\dots,\er\}\,$ and 
$\,\my,\ny\in\{\er+1,\dots,2\hh\er\}$, let $\,c_{\my j}$ be a nonsingular 
$\,\er\times\er\,$ matrix of constants. In the new coordinates 
$\,\y^{\hs j}\nnh,\y^\my$ for $\,N\hs$ defined by 
$\,\q_j=c_{\my j}\y^\my\nnh$, the components of $\,\heh\,$ are 
$\,\heh_{jk}=-\hs2\gy_{jk},\,\heh_{\my\ny}=0\,$ and 
$\,\heh_{\my j}=\heh_{j\my}=c_{\my j}$, so that all three-times ``covariant'' 
Christoffel symbols of $\,\heh\hs$ vanish except 
$\,\widetilde\vg_{\!jkl}=\gy_{jk,l}-\gy_{kl,j}-\gy_{jl,k}$, the comma 
denoting partial differentiation. This gives (a) and (b): the $\,\y^\my$ 
coordinate \vf s are parallel and span $\,\kerd\pmb$. Also (cf.\ (\ref{rlm})), 
$\,\widetilde R_{jklm}\nh=\widetilde\vg_{\!jlm,k}\nh
-\widetilde\vg_{\!klm,j}\nh=\gy_{mk,lj}\nh-\nh\gy_{lk,mj}\nh
-\nh\gy_{mj,lk}\nh+\nh\gy_{lj,mk}$ are the only (possibly) nonzero curvature 
components, which implies (c), Ric\-ci-flat\-ness of $\,\heh\,$ now being 
clear since $\,\heh^{\hs jk}\nh=0$.

\newpage
\part{CONFORMALLY SYMMETRIC MANIFOLDS}\label{csm}
This part deals with the central topic of the present paper. In 
Sections~\ref{bapr} --~\ref{difo} we provide successive steps leading to a 
proof, in Section \ref{pfth}, of our main classification result.

\section{Basic properties}\label{bapr}
Given a \prd\ $\,(M,g)\,$ of dimension $\,n\ge4\,$ which is \cs\ 
(that is, $\,\nabla W\nnh=0$), let $\rw\,$ be the rank of its Weyl tensor 
acting on exterior $\,2$-forms at each \pt, and let $\,\vt$ be the parallel 
\dn\ on $\,M\,$ spanned by all \vf s of the form $\,W(u,v)v'\nnh$, for 
arbitrary \vf s $\,u,v,v'$ on $\,M$.
\begin{lemma}\label{weoto}For any \cs\ \prd\ $\,(M,g)$, the following three 
conditions are equivalent\/{\rm:}
\begin{enumerate}
  \def\theenumi{{\rm\roman{enumi}}}
\item $\text{\rm rank}\,W\nnh=1$,
\item the parallel distribution\/ $\,\vt\hs$ introduced above is two\diml\ and 
null,
\item $W\nnh=\hs\ve\hskip1pt\om\nh\otimes\om\,$ for some $\,\ve=\pm\hs1\,$ and 
some parallel differential\/ $\,2$-form\/ $\,\om\ne0\,$ on\/ $\,M$, defined, 
at each \pt, only up to a sign, and having rank\/ $\,2\,$ at every \pt.
\end{enumerate}
Each of conditions {\rm(i)} -- {\rm(iii)} holds if\/ $\,(M,g)\,$ is \ecs, but 
not \rr, as defined in the Introduction.
\end{lemma}
\begin{proof}Since $\,W$ acting on $\,2$-forms is self-ad\-joint, (i) implies 
(iii): $\,\hs\text{\rm rank}\hskip3.5pt\om\hs=2\,$ as
$\,\om\hs\wedge\hs\om\hs=0\,$ by the first Bianchi identity for 
$\,W\nnh=\hs\ve\hskip1pt\om\nh\otimes\om\hs$. If (iii) holds, $\,\vt\,$ is the 
image of the \vbm\ acting by $\,u\mapsto\hs\om\hh u\,$ (see (\ref{tau}) for 
$\,\tau=\hh\om$), and so $\,\dim\vt=2$. As contractions of $\,W$ vanish, all 
vectors in $\,\vt\,$ are null. Thus, (iii) yields (ii). Finally, if we assume 
(ii) and choose, locally, a differential $\,2$-form $\,\om\hh'$ such that the 
image of the morphism $\,u\mapsto\hs\om\hh'\nh u\,$ is $\,\vt$, then 
$\,\om\hh'$ will span the image of $\,W$ acting on $\,2$-forms, proving (i). 
The final clause is a known result \cite[Theorem~9(ii)]{derdzinski-roter-80}.
\end{proof}
\begin{lemma}\label{cstwd}In any \cs\ \prd\ $\,(M,g)\,$ of dimension\/ 
$\,n\ge4\,$ such that\/ $\rwo$,
\begin{enumerate}
  \def\theenumi{{\rm\alph{enumi}}}
\item the scalar curvature\/ $\,\hs\text{\rm s}\hs\,$ is identically zero,
\item the Ricci tensor\/ $\,\ri\,$ satisfies the Codazzi equation\/ 
$\,d^{\hs\nabla}\hskip-2.3pt\ri=0$, cf.\ Remark\/ $\ref{codaz}$,
\item $R\hs=\hs W\nnh+(n-2)^{-1}\hs g\wedge\nh\ri$, with notations as at the 
end of\/ Section $\ref{expr}$,
\item for every \vf\/ $\,u$, the image $\,\ri u$, defined by\/ 
{\rm(\ref{tau})}, is a section of the \npd\/ $\,\vt\,$ appearing in 
Lemma\/ $\ref{weoto}$$(${\rm ii}$)$.
\end{enumerate}
\end{lemma}
\begin{proof}For $\,\om\,$ as in Lemma~\ref{weoto}(iii), 
$\,R_{jkm}{}^s\hs\om_{sl}=R_{jkl}{}^s\hs\om_{sm}$ by (\ref{rit}) with 
$\,\gy=\hs\om\hh$. Thus, $\,R_k^s\hs\om_{sl}=R_{jkls}\hs\om^{sj}$. Summing 
cyclically over $\,j,k,l$, we get $\,R_{jkls}\hs\om^{sj}=0\,$ (as $\,\vt\,$ is 
the image of $\,\om\hh$, and so $\,R_{klsj}\hs\om^{sj}=0\,$ by 
(\ref{rzo}\hh-b)). Hence 
$\,R_k^s\hs\om_{sl}=0$. Also, $\,W\nnh=\hs\ve\hskip1pt\om\nh\otimes\om\hs$, 
so that $\,W_{jkl}{}^s\hs\om_{sm}=\hs\om_{jk}W_l{}^s{}_{sm}=0$. As 
$\,W\nnh=R-(n-2)^{-1}\hs g\wedge\hh\sj\,$ (see the end of 
Section \ref{expr}), we thus have 
$\,0=(n-1)(n-2)\hs W_{klsj}\hs\om^{sj}=2\hs\text{\rm s}\,\om_{kl}$, and (a) 
follows. Next, since $\,W_{jkl}{}^s\hs\om_{sm}$, $\,R_k^s\hs\om_{sl}$ and 
$\,\hs\text{\rm s}\hs\,$ all vanish, symmetry of $\,R_{jkl}{}^s\hs\om_{sm}$ in 
$\,l,m$, established above, implies analogous symmetry of 
$\,R_{jl}\om_{km}-R_{kl}\om_{jm}$. Choosing a basis of a given tangent space 
in which $\,\om_{jl}=0\,$ unless $\,\{j,l\}=\{1,2\}$, then setting $\,k=1$, 
$\,m=2$, and using the latter symmetry for $\,j,l>2$, or $\,j=2\,$ and 
$\,l>2$, we obtain $\,R_{jl}=0\,$ unless $\,j,l\in\{1,2\}$. As $\,\vt\,$ is 
the image of $\,\om\hh$, this yields (d). Finally,  
$\,\hs\text{\rm div}^\nabla W\nnh=0\,$ (i.e., 
$\,W_{jkl}{}^s{}_{\nh,\hh s}=0$), since $\,\nabla W\nnh=0$. Hence, by the 
second Bianchi identity, the Codazzi equation holds for the Schouten tensor 
$\,\sj\,$ (see the end of Section \ref{expr}). As (a) gives $\,\sj=\ri$, 
(b) and (c) follow.
\end{proof}
For the case $\,\nabla R\ne0\,$ in Lemma~\ref{cstwd}, see also 
\cite[Theorem~7]{derdzinski-roter-78} and 
\cite[Theorem~7]{derdzinski-roter-80}.
\begin{lemma}\label{precs}For any \cs\ \prd\ $\,(M,g)$ with\/ $\rwo$, and for 
$\,\om\hh,\vt\,$ as in Lemma\/ $\ref{weoto}$, let\/ $\,\ri\,$ denote the \rt. 
Then
\begin{enumerate}
  \def\theenumi{{\rm\roman{enumi}}}
\item $\om\hh,\hs\ri,\hs\nabla\nnh\ri\,$ and\/ $\,W$ all annihilate 
$\,\w\nnh$, in the sense of\/ Section $\ref{poct}$,
\item $\vt\,$ satisfies the curvature conditions {\rm(\ref{rvv})}.
\end{enumerate}
\end{lemma}
\begin{proof}To obtain (i), we fix a \vf\ $\,v$. Since $\,\nabla\nnh\ri\,$ is 
totally \sy\ by Lemma~\ref{cstwd}(b), we only need to verify that the images 
of the \vbm s $\,\tm\to\tm$, obtained from $\,\om,\ri\,$ and 
$\,\nabla_{\!v}\ri\,$ by index raising (cf.\ (\ref{tau})), are all contained 
in $\,\vt$. (The case of $\,W$ in (i) is immediate from that of $\,\om\hh$, 
since $\,W\nnh=\hs\ve\hskip1pt\om\nh\otimes\om\hh$.) For $\,\om\hh$, the last 
claim is obvious: $\,\vt\,$ is the image of $\,\om\hh$. For $\,\ri$, it 
follows from Lemma~\ref{cstwd}(b). For $\,\nabla_{\!v}\ri$, it is immediate 
from the assertion about $\,\ri$, since $\,\vt\,$ is parallel. This proves 
(i). Finally, (i) applied to $\,W$ implies (ii), by Lemma~\ref{cstwd}(c), 
since $\,g(\vt,\w\nh)=0$, while $\,\ri$ annihilates $\,\w$ and $\,\vt\,$ by 
(i) and (\ref{inq}\hh-ii).
\end{proof}
\begin{theorem}\label{csvhc}Let\/ $\,(M,g)$ be a \cs\ \prd\ of dimension 
$\,n\ge4$, such that\/ $\rwo$, and let\/ $\,\vt\,$ be the $\,2$\diml\ \npd\ on 
$\,M$, described above. We may assume, locally, that the leaves of\/ 
$\,\w\nnh$ are the fibres of a \bp\ $\,\pmb\nh:M\nh\to\bs\,$ for some \su\ 
$\,\bs$, and so, by {\rm Lemma~\ref{precs}(ii)}, $\,\bs\,$ carries the 
$\,\w\nh$\pd\ \tc\ $\,\hs\text{\rm D}$ defined in Section $\ref{plcc}$. Then
\begin{enumerate}
  \def\theenumi{{\rm\roman{enumi}}}
\item the Ricci tensor $\,\ri^{\hs\text{\rm D}}$ of\/ $\,\hs\text{\rm D}\hs\,$ 
is symmetric,
\item $g\,$ has the Ricci tensor 
$\,\hs\ri\hs=(n-2)\pmb^*\hskip-2.6pt\ri^{\hs\text{\rm D}}\nnh$,
\item $d^{\hs\text{\rm D}}\nh\ri^{\hs\text{\rm D}}\hskip-2.5pt=\nh0$, that is, 
$\,\hs\text{\rm D}\hs\,$ is projectively flat, cf.\ Theorem\/ $\ref{prjfl}$,
\item $\text{\rm D}\hs\ri^{\hs\text{\rm D}}\nnh=0\,$ at\/ 
$\,\pmb(x)\in\bs\,$ if and only if\/ $\,\nabla R=0\,$ at\/ $\,x\in M$,
\item $\ri^{\hs\text{\rm D}}\nnh$ is recurrent, in the sense of 
Section $\ref{scoc}$, if and only if so is\/ $\,\ri$,
\item the $\,2$-form $\,\om\,$ of\/ Lemma\/ $\ref{weoto}$$(${\rm iii}$)$ 
equals\/ $\,\pmb^*\nnh\alpha\,$ for a $\,\hs\text{\rm D}\hs$\prl\ area form\/ 
$\,\alpha\,$ defined, locally in $\,\bs$, up to a sign.
\end{enumerate}
\end{theorem}
\begin{proof}By Lemma~\ref{precs}(i) with $\,\nabla\om=0$, the final clause of 
Lemma~\ref{pullb} applies to both $\,\gy=\om\,$ and $\,\gy=\ri$. We thus get 
(vi) (where $\,\hs\text{\rm D}\hs\alpha=0\,$ in view of Lemma~\ref{pullb}(b) 
for $\,\tau=\alpha$, as $\,\nabla\om=0$), and $\,\ri=(n-2)\hs\pmb^*\tau\,$ for 
some \sy\ covariant $\,2$-ten\-sor $\,\tau\,$ on $\,\bs$.

Lemma~\ref{cstwd}(c) gives $\,(n-2)\hh R(w,w\hs'\hh)\hh v
=g(w,v)\hs\ri\hs w\hs'\nh-\hs g(w\hs'\nnh,v)\hs\ri\hs w\,$ for any section 
$\,v$ of $\,\w$ and \vf s $\,w,w\hs'$ on $\,M$. (Notation of (\ref{tau}); the 
terms involving $\,W$ and $\,\ri\hs v$ vanish as $\,W$ and $\,\ri\,$ 
annihilate $\,\w\nnh$, by Lemma~\ref{precs}(i).) If $\,w\,$ and $\,w\hs'$ are 
$\,\pmb\hs$\pj, while $\,v\,$ is a section of $\,\vt\,$ parallel along 
$\,\w\nnh$, this equality, (v) in Section \ref{plcc} and the relation 
$\,\ri=(n-2)\hs\pmb^*\nh\tau\,$ give 
$\,R^{\hs\text{\rm D}}\nnh=\tau\wedge\hs\text{\rm Id}$ on $\,\bs\,$ (cf.\ 
(\ref{rri})); contracting, we get $\,\tau=\ri^{\hs\text{\rm D}}\nnh$, which 
proves (i) -- (ii). Now Lemmas~\ref{pullb}(b) and~\ref{cstwd}(b) yield (iii) 
-- (v).
\end{proof}

\section{Conformal changes of the metric}\label{ccot}
Given a \cs\ \prd\ $\,(M,g)\,$ with $\rwo$, let us consider the following 
conditions imposed on a function $\,\fy:M\to\bbR\hs$.
\begin{equation}
\begin{array}{rl}\arraycolsep9pt
\mathrm{a)}&(2-n)\hs\nabla\nh\df=\fy\nnh\ri\hs,\hskip5pt
\text{\rm where}\hskip5ptn=\hs\dim M\hskip4pt\text{\rm 
and}\hskip4pt\ri\hskip5.3pt\text{\rm is\ the\ Ricci\ tensor,}\\
\mathrm{b)}&\text{\rm the\ gradient\ $\,\naf\,$\ is\ a\ section\ of\ the\ 
\dn\ $\,\vt$\ (see\ Section \ref{bapr}).}
\end{array}\label{ndf}
\end{equation}
\begin{remark}\label{delfz}Condition (\ref{ndf}\hh-b) for a \npd\ $\,\vt\,$ on 
any \prd\ $\,(M,g)\,$ and a function $\,\fy:M\to\bbR\,$ implies that 
$\,\Delta\fy=g(v,v)=0\,$ for the section $\,v=\naf\,$ of $\,\vt$. In fact, the 
image of $\,\nabla v:\tm\to\tm\,$ is contained in $\,\vt$, and, since 
$\,\nabla v\,$ is self-ad\-joint, its kernel contains $\,\w$ (and $\,\vt$); we 
thus get $\,\Delta\fy=0\,$ evaluating $\,\hs\text{\rm tr}\,\nabla v\,$ in a 
basis of any tangent space $\,\txm$, a part of which spans $\,\vt_x$.
\end{remark}
\begin{lemma}\label{cfrfl}Let\/ $\,(M,g)\,$ be a \cs\ \prd\ of dimension\/ 
$\,n\ge4\,$ such that\/ $\rwo$.
\begin{enumerate}
  \def\theenumi{{\rm\roman{enumi}}}
\item Some neighborhood of any \pt\ of $\,M$ admits a function $\,\fy>0\,$ 
with $(\ref{ndf})$, and, for any such $\,\fy$, the conformally related metric 
$\,\widetilde g=\fy^{-2}g\,$ is Ric\-ci-flat.
\item If\/ $\,M\,$ is \sco, then the \vs\ of functions\/ $\,\fy:\bs\to\bbR\,$ 
with\/ $(\ref{ndf})$ is $\,3$\diml, and such a function $\,\fy\,$ is 
uniquely determined by its value and gradient at any given \pt\/ $\,x$, 
which are arbitrary elements of\/ $\,\bbR\,$ and\/ $\,\vtx$.
\end{enumerate}
\end{lemma}
\begin{proof}The formula $\,\bna(\xi,\fy)=(\nabla\hh\xi
+(n-2)^{-1}\fy\nnh\ri\hh,\hs\df-\xi)$, with $\,\ri\,$ standing for the \rt\ of 
$\,g$, defines a connection $\,\bna\,$ in the \vb\ 
$\,\tam\oplus(M\nh\times\bbR)$. (Notation as in the lines preceding 
Remark~\ref{detby}; sections of $\,\tam\oplus(M\nh\times\bbR)\,$ are pairs 
$\,(\xi,\fy)\,$ consisting of a $\,1$-form $\,\xi\,$ and a function $\,\fy$.) 
The curvature tensor $\,\br\,$ of $\,\bna$, evaluated from 
(\ref{cur}) with $\,\psi=(\xi,\fy)\,$ and (\ref{dnp}), is given by 
$\,\br(u,v)\psi=(\xi\hh'(u,v),0)$, where $\,\xi\hh'
=-\,d^{\hs\nabla}\nh\nabla\xi-(n-2)^{-1}(\xi\wedge\ri
+\fy d^{\hs\nabla}\hskip-3pt\ri)$, 
with $\,d^{\hs\nabla}\nh\nabla\xi\,$ as in (\ref{rid}) for $\,\bs=M\,$ and 
$\,\hs\text{\rm D}\hs=\nabla$. By (\ref{rid}), 
$\,d^{\hs\nabla}\nh\nabla\xi=\xi\hs R$. As 
$\,d^{\hs\nabla}\hskip-2.3pt\ri=0\,$ (see 
Lemma~\ref{cstwd}(b)), $\,\xi\hh'\nh=\hs-\hs\xi\hs R-(n-2)^{-1}\xi\wedge\ri$.

Let the sub\bd\ $\,\yz\,$ of $\,\tam\oplus(M\times\bbR)\,$ be the direct sum 
$\,\zz\oplus(M\nh\times\bbR)$, where $\,\zz\,$ is 
the sub\bd\ of $\,\tam\,$ whose sections are the $\,1$-forms $\,\xi\,$ on 
$\,M\,$ that annihilate $\,\w$ in the sense of Section \ref{poct}. The sub\bd\ 
$\,\yz\,$ is $\,\nabla$\prl, that is, invariant under $\,\nabla$\prl\ 
transports, since the \dn\ $\,\w$ is parallel. Consequently, $\,\yz\,$ is also 
$\,\bna$\prl, due to the definition of $\,\xi\wedge\ri\,$ in 
Section \ref{expr} and the fact that $\,\ri\,$ annihilates $\,\w$ (see 
Lemma~\ref{precs}(i)). Thus, $\,\bna\,$ has a restriction to a connection in 
$\,\yz$. Next, given a section $\,(\xi,\fy)\,$ of $\,\yz$, 
Lemma~\ref{cstwd}(c) allows us to replace $\,R\,$ in $\,\xi\hs R\,$ by 
$\,(n-2)^{-1}\nh g\wedge\ri$. (In fact, $\,\xi=g(v,\,\cdot\,)$, where $\,v\,$ 
is a section of $\,\vt$, so that $\,W(v,\,\cdot\,,\,\cdot\,,\,\cdot\,)=0\,$ 
and $\,\ri\hs v=0\,$ by Lemma~\ref{precs}(i) and (\ref{inq}\hh-ii).) However, 
$\,\xi\hs(g\wedge\ri)=-\hs\xi\wedge\ri$, as $\,\ri\hs v=0$. Thus, 
$\,\xi\hh'\nh=0\,$ (see the last paragraph), and so the restriction of 
$\,\bna\,$ to $\,\yz\,$ is flat.

Finally, if $\,M\,$ is \sco, the assignment $\,\fy\mapsto(\df\nh,\fy)\,$ is a 
linear isomorphism of the space $\,\mathcal{S}\hs$ of all functions 
$\,\fy:\bs\to\bbR\,$ with (\ref{ndf}) onto the space $\,\mathcal{S}\hs'$ of 
$\,\bna$\prl\ sections of $\,\yz$. Namely, injectivity of 
$\,\fy\mapsto(\df\nh,\fy)\,$ is obvious; that it maps $\,\mathcal{S}\hs$ into 
$\,\mathcal{S}\hs'\nh$, and is surjective, follows from the definition of 
$\,\bna$. Our assertion is now immediate, as Remark~\ref{delfz} and (c) in 
Section \ref{expr} imply Ric\-ci-flat\-ness of $\,\widetilde g$.
\end{proof}
\begin{proposition}\label{ecswp}If\/ $\,(M,g)\,$ is a \cs\ \prd\ of 
dimension\/ $\,n\ge4\,$ such that\/ $\rwo$, and\/ $\,\fy:M\to\bbR\,$ is a 
positive function with {\rm(\ref{ndf})}, then the conformally related metric\/ 
$\,\widetilde g=\fy^{-2}g$ on $\,M\hs$ is, locally, a Riemannian product 
of a Ric\-ci-flat metric $\,\heh\,$ with the signature\/ 
\hbox{$\,-$\hskip1pt$-$\hskip1pt$+$\hskip1pt$+$} on a four\mfd\/ $\,N$ and 
a flat metric $\,\gm\,$ on a \mf\/ $\,\mv\nnh$ of dimension\/ $\,n-4\hh$. In 
addition,
\begin{enumerate}
  \def\theenumi{{\rm\roman{enumi}}}
\item the $\,2$\diml\ \npd\/ $\,\vt\hs$ on\/ $\,M\hs$ defined in 
Section $\ref{bapr}$ is tangent to the $\,N$ factor,
\item both $\,\vt\hs$ and\/ $\,\w$ are spanned, locally, by 
$\,\widetilde g$\prl\ \vf s,
\item $\fy\,$ is constant in the direction of the $\,\mv$ factor.
\end{enumerate}
\end{proposition}
\begin{proof}From (a) in Section \ref{expr} for any \vf\ $\,u\,$ and any 
section $\,v\,$ of $\,\w\nnh$, with $\,\phi=-\hs\log\fy\,$ and 
$\,w=\nabla\nh\phi$, we get $\,\widetilde\nabla_{\!u}v\,=\,\nabla_{\!u}v\,
+\,g(u,w)\hs v\,-\,g(u,v)\hs w$, where $\,g(v,w)\hs u=0$ since 
$\,w=\nabla\nh\phi=-\fy^{-1}\naf\,$ is, by (\ref{ndf}\hh-b), a section of $\,\vt$. 
The sub\bd s $\,\w$ and $\,\vt\,$ of $\,\tm$, with $\,\vt\subset\w\nnh$, are 
thus $\,\widetilde\nabla$\prl\ (invariant under $\,\widetilde\nabla$\prl\ 
transports). Furthermore, the restriction of $\,\widetilde\nabla\,$ to the 
sub\bd\ $\,\w$ is flat: as $\,\widetilde g=\fy^{-2}g\,$ is Ric\-ci-flat 
(Lemma~\ref{cfrfl}(i)), its four-times covariant curvature tensor 
$\,\widetilde R\,$ equals its Weyl tensor $\,\widetilde W=\fy^{-2}W\nh$, and 
so it annihilates $\,\w$ by Lemma~\ref{precs}(i). This proves (ii).

Let $\,\,U\,$ be a \sco\ neighborhood of any given \pt\ of $\,M$. The 
$\,\widetilde g$\prl\ (that is, $\,\widetilde\nabla$\prl) sections of $\,\w$ 
defined on $\,\,U\,$ form a \vs\ $\,\mathcal{E}\,$ of dimension $\,n-2$, on 
which $\,\widetilde g\,$ is a degenerate symmetric bilinear form. The 
$\,2$\diml\ $\,\widetilde g$-null\-space $\,\mathcal{E}\hs'$ of 
$\,\mathcal{E}\,$ (i.e., the $\,\widetilde g$-or\-thog\-o\-nal complement of 
$\,\mathcal{E}$) consists 
of all $\,\widetilde g$\prl\ sections of $\,\vt$. Since $\,\widetilde g\,$ 
descends to a (nondegenerate) pseu\-\hbox{do\hs-}Euclid\-e\-an inner product 
on $\,\mathcal{E}/\mathcal{E}\hs'\nnh$, we may choose a $\,(n-4)$\diml\ 
vector subspace of $\,\mathcal{E}/\mathcal{E}\hs'\nnh$, on which 
$\,\widetilde g\,$ is nondegenerate, and realize it as the projection image of 
some subspace $\,\mathcal{X}\,$ of $\,\mathcal{E}\,$ with 
$\,\dim\mathcal{X}=n-4\,$ such that $\,\widetilde g\,$ is nondegenerate on 
$\,\mathcal{X}$. The \dn\ of dimension $\,n-4$ on $\,\,U\,$ spanned by the 
\vf s forming $\,\mathcal{X}\,$ is clearly $\,\widetilde g$\prl\ and 
nondegenerate. Hence, by a classical result of Thomas \cite{thomas}, it is a 
factor \dn, tangent to some factor manifold $\,(\mv\nnh,\gm)$, in a local 
Riem\-ann\-i\-an-prod\-uct decomposition of $\,(M,g)$. Obviously, $\,\gm\,$ is 
flat, as elements of $\,\mathcal{X}\,$ constitute $\,\gm$\prl\ \vf s spanning 
the \tbd\ of $\,\mv\nnh$. Denoting by $\,(N,h)\,$ the other factor manifold, we 
get (i) (since $\,\vt\,$ is orthogonal to all sections of $\,\w\nnh$, 
including elements of $\,\mathcal{X}$), as well as (iii) (which follows from 
(i) and (\ref{ndf}\hh-b)).
\end{proof}
\begin{remark}\label{ddfrh}Locally, condition (\ref{ndf}) can be conveniently 
rephrased in terms of the quotient surface $\,\bs\,$ with the $\,\w\nh$\pd\ 
connection $\,\hs\text{\rm D}\hs\,$ (cf.\ Theorem~\ref{csvhc}). Namely, 
(\ref{ndf}\hh-b) means that $\,\fy\,$ may be treated as a function 
$\,\bs\to\bbR\hs$, and then (\ref{ndf}) is equivalent to the condition 
$\,\hs\text{\rm D}\hh\df=-\fy\nnh\ri^{\hs\text{\rm D}}\nnh$. This is immediate 
from Theorem~\ref{csvhc}(ii) and Lemma~\ref{pullb}(b) applied to both 
$\,\gy=\fy\,$ and $\,\gy=\df$.

Lemma~\ref{cfrfl}(ii) could therefore be derived from an analogous statement 
about the equation $\,\hs\text{\rm D}\hh\df=-\fy\nnh\ri^{\hs\text{\rm D}}$ on 
$\,\bs$, which can be proved by the same argument, and constitutes the 
$\,2$\diml\ case of a more general result of Gardner, Kriele and Simon 
\cite{gardner-kriele-simon}.
\end{remark}

\section{Reduction of the dimension}\label{rotd}
Our next step in classifying \cs\ \mf s with $\rwo$, in dimensions $\,n\ge4$, 
consists in reducing the problem to the case $\,n=4\hh$.

For \prd s $\,(N,h),(\mv\nnh,\gm)\,$ and a $\,C^\infty$ function 
$\,\fy:N\to(0,\infty)$, the {\it warped product\/} with the {\it base\/} 
$\,(N,h)$, {\it fibre} $\,\hh(\mv\nnh,\gm)$, and {\it warping function} 
$\,\fy\,$ is the \prd\ $\,(M,g)$ given by $\,M=N\nh\times\nh\mv$ and 
$\,g\,=\,h\hs+\fy^2\gm$. (Here $\,g,\fy,h$ also stand for their own 
pull\-backs to $\,N\nh\times\mv\nh\nnh$.) See \cite[p.~237]{besse}. Note that 
$\,g\hs=\hs h\hs+\fy^2\gm\,$ is conformal\-ly related to the product metric 
$\,\fy^{-2}g\hs\,=\,\fy^{-2}h\hs+\hh\gm$.

In the next lemma, all Weyl and curvature tensors are four-times covariant 
tensors.
\begin{lemma}\label{wtwpr}Given \prd s $\,(N,\heh)\,$ and $\,(\mv\nnh,\gm)\,$ 
such that\/ $\,\heh\,$ is Ric\-ci-flat and\/ $\,\gm\,$ is flat, while 
$\,\dim N\ge4$, let $\,\fy\,$ be a positive function on $\,N\nh$, and let\/ 
$\,W^h\nnh$ denote the Weyl tensor of the conformally related metric 
$\,h=\fy^2\heh\,$ on $\,N\nh$. The Weyl tensor $\,W$ of the warped product 
metric $\,g=h\hs+\fy^2\gm\,$ on $\,N\nh\times\mv$ then equals the pull\-back 
of\/ $\,\hs W^h\nnh$ under the projection $\,N\nh\times\mv\to\mv\nh$.
\end{lemma}
This is clear from the conformal transformation rule for $\,W\nnh$. Namely, 
$\,W$ is related to the four-times covariant Weyl and curvature tensors 
$\,\widetilde W\nh,\widetilde R\,$ of the product metric 
$\,\widetilde g=\fy^{-2}g=\heh+\gm\,$ by 
$\,W\nh=\fy^2\widetilde W\nh=\fy^2\widetilde R$, where 
$\,\widetilde W\nh=\widetilde R\,$ since $\,\widetilde g\,$ is Ric\-ci-flat. 
Next, $\,\widetilde R\,$ is the pull\-back of the curvature tensor $\,R\hs'$ 
of $\,\heh\,$ under $\,N\nh\times\mv\to\mv\nnh$, while, as $\,\heh\,$ is 
Ric\-ci-flat, $\,R\hs'$ equals the Weyl tensor $\,W'\nh$ of 
$\,\heh=\fy^{-2}h$, that is, $\,\fy^{-2}W^h\nnh$.
\begin{theorem}\label{dimre}In any \cs\ \prd\/ $\,(M,g)$ of dimension\/ 
$\,n\ge4\,$ such that\/ $\rwo$, every \pt\ has a connected neighborhood 
isometric to a warped product\/ $\,(N\nh\times\mv,\,h\hs+\fy^2\gm)$, for any 
solution $\,\fy>0\,$ to {\rm(\ref{ndf})}, and
\begin{enumerate}
  \def\theenumi{{\rm\alph{enumi}}}
\item $(N,h)\,$ is a \cs\ \mf\ with $\,\dim N\nh=4\,$ and\/ $\rwho$,
\item $(\mv\nnh,\gm)\,$ is flat, $\,\dim\hs\mv\nnh=\hs n-4\ge0$, and\/ $\,\fy\hs$ 
is constant in the $\,\mv\nnh$-factor direction,
\item $\fy\,$ treated as a function on\/ $\,(N,h)\,$ satisfies 
{\rm(\ref{ndf})} with\/ $\,n\,$ replaced by $\,4\hh$.
\end{enumerate}
Conversely, every warped product\/ $\,(M,g)=(N\nh\times\mv,\,h\hs+\fy^2\gm)\,$ 
with \hbox{{\rm(a)} -- {\rm(c)}} is \cs, has\/ $\rwo$, and\/ $\,\fy\,$ 
satisfies\/ {\rm(\ref{ndf})} on $\,(M,g)$.
\end{theorem}
\begin{proof}Let $\,(M,g)\,$ satisfy the conditions $\,\nabla W\nnh=0\,$ and 
$\rwo$. By Proposition~\ref{ecswp}, every \pt\ of $\,M\,$ has a 
connected neighborhood $\,\,U\,$ such that $\,(U\nh,\fy^{-2}g)\,$ is isometric 
to a Riemannian product \mf\ $\,(N\nh\times\mv,\,\heh+\gm)$, in which the 
factor $\,(N,\heh)\,$ is $\,4$\diml\ and Ric\-ci-flat, while $\,(\mv\nnh,\gm)\,$ 
is flat; since $\,\fy$, which exists by Lemma~\ref{cfrfl}(i), may be treated 
as a function on $\,N\hs$ (Proposition~\ref{ecswp}(iii)), $\,(U\nh,g)\,$ is 
isometric to the warped product $\,(N\nh\times\mv,\,h\hs+\fy^2\gm)$, where 
$\,h=\fy^2\heh$.

As the $\,N$-factor sub\mf s are totally geodesic in $\,(M,g)\,$ and 
$\,h\,$ represents their sub\mf\ metrics, restricting to any of them the 
$\,g$\prl\ $\,2$-form $\,\om\,$ appearing in Lemma~\ref{weoto}(iii) we obtain 
an $\,h$\prl\ $\,2$-form $\,\om\hh^h$ on $\,N\nh$. At the same time, the 
relation $\,W\nnh=\hs\ve\hskip1pt\om\nh\otimes\om\hs\ne0\,$ combined with 
Lemma~\ref{wtwpr} gives 
$\,W^h\nnh=\hs\ve\hskip1pt\om\hh^h\nnh\otimes\om\hh^h\nnh\ne0$, so that 
$\,\hs\text{\rm rank}\hskip3.5pt\om\hh^h\nnh=2\,$ by the first Bianchi 
identity for $\,W^h\nnh$. Hence $\,(N,h)\,$ is \cs\ and, by Lemma~\ref{weoto}, 
$\,\hs\text{\rm rank}\hskip2.4ptW^h\hskip-2pt=\hskip-1.2pt1$. Also, 
$\,\vt\,$ (tangent to the $\,N$ factor by Proposition~\ref{ecswp}(i)), 
is the \dn\ defined as in Section~\ref{bapr} for $\,(N,h)\,$ rather than 
$\,(M,g)$. By Remark~\ref{delfz}, $\,\fy\,$ treated as a function on $\,N$ has 
vanishing $\,h$-La\-plac\-i\-an and its $\,h$-gra\-di\-ent is null, so 
that Ric\-ci-flatness of $\,\heh=\fy^{-2}h\,$ and (c) in Section \ref{expr} 
give (\ref{ndf}\hh-a) with $\,n=4$.

Conversely, let $\,(M,g)=(N\nh\times\mv,\,h\hs+\fy^2\gm)$, with (a) -- (c). As 
$\,W^h\nnh=\hs\ve\hskip1pt\om\hh^h\nnh\otimes\om\hh^h\nnh\ne0$ in $\,(N,h)$, 
where $\,\om\hh^h$ is an $\,h$\prl\ $\,2$-form (Lemma~\ref{weoto}(iii)) we get 
$\,W\nnh=\hs\ve\hskip1pt\om\nh\otimes\om\hs\ne0$, by Lemma~\ref{wtwpr}, for 
the $\,2$-form $\,\om\,$ on $\,M\,$ obtained as the pull\-back of $\,\om\hh^h$ 
under the factor projection $\,M\to N\nh$. For the \lcc\ 
$\,\widetilde\nabla\,$ of the product metric 
$\,\widetilde g=\fy^{-2}g=\heh+\gm\,$ and any \vf\ $\,u$, we have 
$\,\widetilde\nabla_{\!u}\om=
-\,2(d_u\phi)\hs\om+\om(u,\,\cdot\,)\wedge\hs d\phi$, where 
$\,\phi=-\hs\log\fy$, since the same is true for the \lcc\ of $\,\heh\,$ and 
$\,\om\hh^h$ (rather than $\,\om$) in view of (b) in Section \ref{expr} and 
the fact that $\,\om\hh^h$ is $\,h$\prl. Using (b) in Section \ref{expr} 
again, we now see that $\,\om\,$ is $\,g$\prl. Thus, by Lemma~\ref{weoto}, 
$\,(M,g)\,$ is \cs\ and $\rwo$.
\end{proof}
\begin{corollary}\label{wrppr}Theorem\/ $\ref{dimre}$ remains true also when 
the dimensions $\,4\,$ and\/ $\,n-4$ in\/ {\rm(a)} -- {\rm(c)} are replaced 
by $\,k\hs$ and $\,n-k$, for any given $\,k\in\{4,\dots,n\}$.
\end{corollary}
In fact, writing $\,\gm\,$ in Theorem~\ref{dimre} as the product 
$\,\gm\hs''\nh+\gm\hs'$ of two flat metrics, we obtain a new 
warp\-ed-prod\-uct decomposition of $\,g=h\hs+\fy^2\gm$, namely, 
$\,g=h'\nh+\fy^2\gm\hs'\nnh$, with $\,h'\nh=h\hs+\fy^2\gm\hs''$ having the 
required properties as a consequence of Theorem~\ref{dimre}.

\section{Dimension four}\label{difo}
Four\diml\ \cs\ \mf s $\,(N,h)\,$ with $\rwho$, where $\,W^h$ is 
the Weyl tensor of $\,h$, can be constructed as follows. Let there 
be given an \ea\ \pftc\ $\,\hs\text{\rm D}\hs\,$ on a surface $\,\bs\,$ 
with a $\,\hs\text{\rm D}\hs$\prl\ area element $\,\pm\hs\alpha\,$ (see 
Section \ref{eaco}), a factor $\,\ve\in\{1,-1\}$, and a twice-co\-var\-i\-ant 
\sy\ \tf\ $\,\tau$ on $\,\bs\,$ with $\,\op\tau=\ve\hs\alpha\otimes\alpha$. 
Locally, such $\,\tau\,$ always exists (Theorem~\ref{slveq}(i)). We now set 
$\,(N,h)=(\tab,\,\,h^{\text{\rm D}}\hskip-2.5pt-2\tau)$, where 
$\,h^{\text{\rm D}}$ is the metric on $\,\tab\,$ given by (\ref{rex}), and the 
symbol $\,\tau\,$ also stands for the pull\-back of $\,\hs\tau\,$ to 
$\,\tab$. It should be \pt ed out that the tensor $\,\tau$ is not really a 
parameter for the above construction: locally, up to an isometry, different 
choices of $\,\tau\,$ lead to the same metric 
$\,h^{\text{\rm D}}\hskip-2.5pt-2\tau$. See Section \ref{tlst}.
\begin{theorem}\label{clsdf}Every pseu\-\hbox{do\hs-}Riem\-ann\-i\-an 
four\mfd\/ $\,(N,h)=(\tab,\,\,h^{\text{\rm D}}\hskip-2.5pt-2\tau)$ obtained as 
above is \cs\ and\/ $\rwho$, for its Weyl tensor $\,W^h\nnh$, while 
$\,\vt\hs$ defined as in\/ Section $\ref{bapr}$ coincides with \tvn\ on 
$\,\tab$, and\/ $\,\hs\text{\rm D}$ used in the construction is the same as 
the $\,\w\nh$\pd\ connection in Theorem\/ $\ref{csvhc}$.

Conversely, every \pt\ of any \cs\ \psr\ four\mfd\/ $\,(N,h)\,$ with 
$\rwho\,$ has a connected neighborhood isometric to an open subset of a \mf\/ 
$\,(\tab,\,\,h^{\text{\rm D}}\hskip-2.5pt-2\tau)\,$ constructed as above.
\end{theorem}
We precede the proof of Theorem~\ref{clsdf} with two lemmas. Conformal 
flatness of $\,\hs\deh\,$ in Lemma~\ref{pwmcf} is a result of Walker 
\cite[p.\ 69]{walker-re}. Lemma~\ref{walkr} is due to Ruse \cite{ruse} for 
$\,\er=2\,$ (which is the only case that we need), and to Walker \cite{walker} 
for arbitrary $\,\er$.
\begin{lemma}\label{pwmcf}For a \pftc\/ $\,\hs\text{\rm D}\hs$ on a \mf\/ 
$\,\bs\,$ whose \rt\ is \sy, the metric\/ $\,h^{\text{\rm D}}$ on $\,\tab$, 
given by $(\ref{rex})$, is \cf. Furthermore, using Lemma\/ $\ref{smgeo}$ we 
can find, locally in $\,\bs$, a function $\,\fy>0\,$ such that, for 
$\,\xi=-\hs d\hh\log\fy$, the connection $\,\hs\deh\hs
=\hs\text{\rm D}\hs+2\hs\xi\odot\nh\text{\rm Id}\hs\,$ is flat. Choosing any 
local coordinates $\,\y^{\hs j}$ for $\,\bs\,$ in which the components of\/ 
$\,\hs\deh\,$ all vanish and letting $\,\y^{\hs j}\nnh,\p_j$ be the 
corresponding  coordinates in $\,\tab$, we then have the 
sym\-met\-ric-prod\-uct relations
\begin{equation}
\aligned
&\mathrm{a)}\hskip5pt
\vg_{\!kl}^j\,d\y^k\hs d\y^{\hh l}\hs=\,2\fy^{-1}\df\,d\y^{\hs j}\nh,
\hskip33pt
\mathrm{b)}\hskip5pt
\fy^{-2}h^{\text{\rm D}}\hs=\,2\,d(\fy^{-2}\p_j)\,d\y^{\hs j}\nh,\\
&\mathrm{c)}\hskip5pt
\fy^{-1}R_{kl}\,d\y^k\hs d\y^{\hh l}\hs=\,d\hs[\hh\partial_j(\fy^{-1})]\,d\y^{\hs j}\nh,
\hskip10pt
\mathrm{where}\hskip6pt\partial_j=\hs\partial/\partial\y^{\hs j}\nh,
\hskip6pt
\label{phg}
\endaligned
\end{equation}
$\vg_{\!kl}^j$ and $\,R_{kl}$ being the components of 
$\,\hs\text{\rm D}\hs\,$ and its \rt\ $\,\ri^{\hs\text{\rm D}}\nnh$. 
\end{lemma}
\begin{proof}Since 
$\,\fy\vg_{\!kl}^j=\fy_{,k}\delta_l^j\nh+\fy_{,l}\delta_k^j$ (cf.\ 
Lemma~\ref{smgeo}), (a) follows. Now (\ref{hde}) gives (b). Conformal flatness 
of $\,\hs h^{\text{\rm D}}\nh$ is in turn obvious from (b): the metric 
$\,\fy^{-2}h^{\text{\rm D}}\nnh$, having constant component functions in the 
new coordinates $\,\y^{\hs j},\fy^{-2}\p_j$, is flat. Finally, by 
Remark~\ref{ddfrd}, 
$\,\hs\text{\rm D}\hh\df=-\fy\nnh\ri^{\hs\text{\rm D}}\nnh$, that is, 
$\,\nh\fy\nh R_{kl}=-\fy_{,kl}$, while 
$\,\nh\fy_{,kl}=\partial_l\partial_k\fy-\vg_{\!kl}^j\,\partial_j\fy$. Thus, 
(a) implies (c).
\end{proof}
\begin{lemma}\label{walkr}Every \prc\ $\,\heh\,$ on a \mf\/ $\,N$ of dimension 
$\,2\er\hs$ with an $\,\er$\diml\ $\,\heh$-null \dn\/ $\,\vt\hs$ spanned by 
$\,\heh$\prl\ \vf s has, locally, the form $(\ref{tih})$ in some coordinates 
$\,\y^{\hs j}\nnh,\q_j$ such that the $\,\q_j$ coordinate \vf s span $\,\vt$. 
In addition, the \lcc\ $\,\widetilde\nabla$ of\/ $\,\heh\,$ is $\,\vt$\pj, 
cf.\ Section $\ref{plcc}$, its $\,\vt$\pd\ connection\/ $\,\hs\deh\hs$ on a 
local leaf space $\,\bs\,$ for $\,\vt\,$ is flat, and the components of\/ 
$\,\hs\deh\hs\,$ are all zero in the local coordinates for $\,\bs\hs$ provided 
by the functions $\,\y^{\hs j}\nnh$.
\end{lemma}
\begin{proof}Letting $\,\,U\,$ stand for increasingly small neighborhoods of 
any given \pt\ $\,x\in N\nh$, we may choose $\,\er\hs$ functions 
$\,\y^{\hs j}$ on $\,\,U\,$ whose $\,\heh$-gra\-di\-ents 
$\,v_j=\widetilde\nabla\y^{\hs j}$ are linearly independent $\,\heh$\prl\ 
sections of $\,\vt$. The Ricci identity (\ref{rid}\hh-i) for 
$\,\hs\text{\rm D}\hs=\widetilde\nabla$, its curvature tensor 
$\,\widetilde R$, and $\,\xi=\heh(v_j,\,\cdot\,)\,$ shows that 
$\,\widetilde R\,$ annihilates $\,\vt\,$ (cf.\ Section \ref{poct}), and so 
$\,\vt$-pro\-ject\-a\-bil\-i\-ty of $\,\widetilde\nabla\,$ follows from 
Lemma~\ref{rvveq}. Being constant along $\,\w\nnh=\vt$, our $\,\y^{\hs j}$ 
descend to $\,\bs$, where they form a local coordinate system with 
$\,\hs\deh\hs$\prl\ differentials $\,d\y^{\hs j}$ (by (iii) -- (iv) in 
Section \ref{plcc}); thus, the component functions of the $\,\vt$\pd\ 
connection $\,\hs\deh\,$ all vanish. Any fixed sub\mf\ $\,\bs\hs'$ of 
$\,\,U\,$ such that $\,x\in\bs\hs'$ and the projection $\,\pmb:U\to\bs\,$ sends 
$\,\bs\hs'$ dif\-feo\-morphic\-ally onto $\,\bs\,$ gives rise to $\,\er\hs$ 
functions $\,\q_j$ on $\,\,U\,$ with $\,\q_j=0\,$ on $\,\hs U\nh\cap\bs\hs'$ 
and $\,(d\q_j)(v_k)=\delta_{jk}$. (As the \vf s $\,v_j$ commute, they form 
coordinate vector fields on each leaf $\,Y\hs$ of $\,\vt$, for the coordinates 
$\,\q_j$ on $\,Y\nh$.) Thus, $\,\tau(v_k,\,\cdot\,)=0\,$ for 
$\,\tau=\heh-2\hs d\q_j\hs d\y^{\hs j}$ and each $\,k$, and so 
$\,\tau=-\hs2\gy_{kl}\,d\y^k\hs d\y^{\hh l}$ for some $\,\gy_{kl}$. In the 
coordinates $\,\y^{\hs j}\nnh,\y^\my$ for $\,N\hs$ obtained from 
$\,\y^{\hs j}\nnh,\q_j$ as at the end of Section \ref{riex}, the relations 
$\,\widetilde\nabla v_j=0\,$ now give 
$\,\partial_\my\gy_{jk}=\widetilde\vg_{\!jk\my}=0$, as 
$\,\heh_{j\my}=c_{\my j}$ are constants and $\,\heh_{jk}=-\hs2\gy_{jk}$.
\end{proof}
By (\ref{hde}), 
$\,h^{\text{\rm D}}\hskip-2.5pt-2\tau\,=\,2\hs dp_j\,d\y^{\hs j}
-\hs2(\hs\p_j\vg_{\!kl}^j+\hs\tau_{kl})\,d\y^k\hs d\y^{\hh l}$ in the 
coordinates $\,\y^{\hs j}\nnh,\p_j$ for $\,\tab$ corresponding to local 
coordinates $\,\y^{\hs j}$ in $\,\bs$. This expression can be further 
simplified if $\,\y^{\hs j}\nnh$, rather than being arbitrary, are chosen, 
along with a function $\,\fy>0$, as in Lemma~\ref{pwmcf}. In the new 
coordinates $\,\y^{\hs j}\nnh,\fy^{-2}\p_j$ for $\,\tab$, we then get, from 
(\ref{phg}\hh-b),
\begin{equation}
\heh\,\,=\,\,2\hs d(\fy^{-2}\p_j)\,d\y^{\hs j}\hs
-\,2\fy^{-2}\tau_{jk}\,d\y^{\hs j}\hs d\y^k\hh,\hskip7pt\text{\rm where
$\,\heh=\fy^{-2}(\hh h^{\text{\rm D}}\hskip-2.5pt-2\tau)$.}
\label{tge}
\end{equation}
\begin{proof}[Proof of Theorem\/ $\ref{clsdf}$]Our $\,\heh\,$ in 
(\ref{tge}) is a special case of (\ref{tih}). Hence, by (c) in 
Section \ref{riex}, $\,\heh\,$ has the Weyl tensor 
$\,\widetilde W\nnh=\pmb^*(\widetilde\op(\fy^{-2}\tau))$, for 
$\,\widetilde\op\,$ associated with the connection $\,\hs\deh$ \hbox{of 
Lemma~\ref{pwmcf},} so that $\,\hs\deh\hs=\hs\text{\rm D}\hs
+2\hs\xi\odot\nh\text{\rm Id}\hs\,$ for 
$\,\xi=-\hs d\hh\log\fy$. Lemma~\ref{trfru}(a) now gives 
$\,\widetilde{\op}(\fy^{-2}\tau)=\fy^{-2}\op\tau$, and so the Weyl tensor of 
$\,h=h^{\text{\rm D}}\hskip-2.5pt-2\tau=\fy^2\heh\,$ is 
$\,W^h\nnh=\fy^2\widetilde W=\pmb^*(\op\tau)
=\hs\ve\hskip1pt\om\hh^h\nnh\otimes\om\hh^h\nnh$, where 
$\,\om\hh^h\nnh=\pmb^*\nnh\alpha$. Also, \tvn\ $\,\vt\,$ on $\,\tab\,$ is 
$\,\heh$-null and $\,\heh$\prl, while the \lcc\ $\,\widetilde\nabla\,$ of 
$\,\heh\,$ is $\,\w\nh$\pj\ and $\,\hs\deh\hs\,$ is its $\,\w\nh$\pd\ 
connection. (See (a), (b) in Section \ref{riex}.) Hence, by (vi) in Section 
\ref{plcc}, the same is true for $\,\vt\,$ if one replaces 
$\,\heh,\hs\deh\hs\,$ by $\,h,\hs\text{\rm D}\hs$. As $\,\alpha\,$ is 
$\,\hs\text{\rm D}\hs$\prl, Lemmas~\ref{rvveq} and~\ref{pullb}(b) for 
$\,\tau=\alpha\,$ now imply that $\,\om\hh^h$ is $\,h$\prl. The first part of 
Theorem \ref{clsdf} thus follows from Lemma~\ref{weoto}. (The image of 
$\,\om\hh^h$ is $\,\vt$, as $\,\om\hh^h$ annihilates $\,\vt=\w\nnh$.)

For the second part of Theorem \ref{clsdf}, let $\,(N,h)\,$ be \cs, with 
$\,\dim N\nh=4\,$ and $\rwho$. By Lemma~\ref{cfrfl}(i), every $\,x\in N$ has a 
neighborhood $\,\,U\,$ such that both (\ref{ndf}) and (\ref{cnd}) hold, for 
$\,\er=2$, some function $\,\fy:U\to(0,\infty)$, a suitable quotient \su\ 
$\,\bs$, and the \dn\ $\,\vt\,$ defined in Section \ref{bapr}, 
with $\,\,U\,$ instead of $\,M$ in (\ref{cnd}). According to 
Remark~\ref{ddfrh}, $\,\fy\,$ descends to a function $\,\bs\to\bbR\hs$, and 
$\,\hs\text{\rm D}\hh\df=-\fy\nnh\ri^{\hs\text{\rm D}}$ for 
$\,\fy:\bs\to\bbR\,$ and the $\,\w\nh$\pd\ connection $\,\hs\text{\rm D}\hs\,$ 
on $\,\bs$, which in turn means that the connection $\,\hs\deh\hs
=\hs\text{\rm D}\hs+2\hs\xi\odot\nh\text{\rm Id}\hs\,$ on $\,\bs$, with 
$\,\xi=-\hs d\hh\log\fy$, is flat (Remark~\ref{ddfrd}).

Proposition~\ref{ecswp}(ii) (with 
$\,n=4,\hs\w\nh=\vt,\hs g=h,\hs\widetilde g=\heh$) and Lemma~\ref{walkr} imply 
that, locally, $\,\heh\,$ has the form (\ref{tih}), for some 
twice-co\-var\-i\-ant \sy\ \tf\ $\,\gy\,$ on $\,\bs$, and local coordinates 
$\,\y^{\hs j}$ in $\,\bs$, in which the components of $\,\hs\deh\hs\,$ all 
vanish. By (c) in Section \ref{riex}, 
$\,\widetilde W\nnh=\pmb^*(\widetilde\op\gy)$. Therefore, 
$\,\ve\hskip1pt\om\hh^h\nnh\otimes\om\hh^h\nh=W^h\nnh=\fy^2\widetilde W
=\pmb^*(\op\tau)\,$ for $\,\tau=\fy^2\gy\,$ and $\,\op$ corresponding to 
$\,\hs\text{\rm D}\hs\,$ as in (\ref{opf}), with the successive equalities due 
to Lemma~\ref{weoto}(iii), the relation $\,h=\fy^2\heh$, and 
Lemma~\ref{trfru}(a). However, by Theorem~\ref{csvhc}(vi), 
$\,\om\hh^h\nnh=\pmb^*\nnh\alpha$, where $\,\pm\hs\alpha\,$ is a 
$\,\hs\text{\rm D}\hs$\prl\ area element on $\,\bs$. Hence 
$\,\op\tau=\ve\hs\alpha\otimes\alpha$. By (\ref{tge}), $\,h\,$ and 
$\,h^{\text{\rm D}}\hskip-2.5pt-2\tau$ are isometric, as they have the same 
form in suitable coordinates.
\end{proof}
All four\diml\ \cs\ \mf s with $\rwo\,$ belong to a class of conformally 
recurrent four\mfd s for which Olszak \cite{olszak} provided a local 
classification. Olszak's result states that the \mf s in question have, 
locally, the form $\,(\tab,h^{\text{\rm D}}\hskip-2.5pt-2\tau)\,$ for some 
\tc\ $\,\hs\text{\rm D}\,$ on a \su\ $\,\bs\,$ and some twice-co\-var\-i\-ant 
\sy\ \tf\ $\,\tau\,$ on $\,\bs$. To derive Theorem~\ref{clsdf} from it, one 
needs to determine what requiring $\,h^{\text{\rm D}}\hskip-2.5pt-2\tau\,$ to 
be \cs\ means for $\,\hs\text{\rm D}\hs\,$ and $\,\tau$. As our discussion 
shows, the answer is: $\,\op\tau=\ve\hs\alpha\otimes\alpha\,$ for some 
$\,\hs\text{\rm D}\hs$\prl\ area element $\,\pm\hs\alpha$.

\section{The local structure theorem}\label{tlst}
Although a local classification of \cs\ \mf s with $\rwo$ could now be easily 
derived from Theorems~\ref{dimre} and~\ref{clsdf}, we state it differently, 
for reasons explained in Remark~\ref{notwp} below. Namely, let the following 
objects be given:
\begin{enumerate}
  \def\theenumi{{\rm\roman{enumi}}}
\item an integer $\,n\ge4\hh$,
\item a surface $\,\bs\,$ with a \pftc\ $\,\hs\text{\rm D}\hs$,
\item a $\,\hs\text{\rm D}\hs$\prl\ area element $\,\pm\hs\alpha\,$ on 
$\,\bs\,$ (see Section \ref{eaco}),
\item a sign factor $\,\ve=\pm\hs1$,
\item a \rvs\ $\,\mv\hs$ of dimension $\,n-4\hh$,
\item a pseu\-\hbox{do\hs-}Euclid\-e\-an inner product $\,\lr\,$ on $\,\mv\nnh$.
\end{enumerate}
We are also assuming the existence of a twice-co\-var\-i\-ant \sy\ tensor 
$\,\tau\,$ on $\,\bs\,$ with 
$\,\op\tau
=\ve\hs\alpha\otimes\alpha$. Locally, such $\,\tau\,$ always exists, by 
Theorem~\ref{slveq}(i).

The data (i) -- (vi) give rise to the $\,n$\diml\ \prd
\begin{equation}\label{hgt}
(M,g)\,\,=\,\,(\tab\,\times\,\mv,\,\,h^{\text{\rm D}}\hskip-2.5pt-2\tau+\gm
-\theta\ri^{\hs\text{\rm D}})\,,
\end{equation}
where $\,h^{\text{\rm D}}$ is the metric (\ref{rex}) on $\,\tab\,$ and 
$\,\gm\,$ is the constant \prc\ on $\,\mv$ corresponding to the inner product 
$\,\lr$, while $\,\theta:\mv\nh\to\bbR\,$ is given by $\,\theta(v)=\lg v,v\rg$. 
As before, the symbols for functions or covariant tensor fields on $\,\bs$, or 
on the factor manifolds $\,\tab\,$ and $\,\mv\nnh$, are also used to represent 
their pull\-backs to $\,\tab\,\times\,\mv\nnh$.

For coordinate descriptions of the metric $\,g\,$ in (\ref{hgt}), see formulae 
(\ref{ged}) in Section \ref{pfth}.

We treat $\,n,\bs,\hs\text{\rm D}\hh,\alpha,\ve,\mv$ and $\,\lr\,$ in (i) -- 
(vi) as parameters for our construction, while $\,\tau\,$ is merely an object 
assumed to exist, even though the metric $\,g\,$ in (\ref{hgt}) clearly 
depends on $\,\tau$. The reason is that, if the data (i) -- (vi) are fixed, 
the metrics corresponding to two choices of $\,\tau\,$ are, locally, isometric 
to each other (Remark~\ref{plbck} below).

We can now state our local classification result; for a proof, see 
Section \ref{pfth}.
\begin{theorem}\label{class}The \prd\/ {\rm(\ref{hgt})} obtained as above from 
data {\rm(i)} -- {\rm(vi)} with the stated properties is \cs\ and has\/ 
$\rwo$. Also,
\begin{enumerate}
  \def\theenumi{{\rm\alph{enumi}}}
\item the \mf\/ {\rm(\ref{hgt})} is \ls, or \rr, if and only if so is\/ 
$\,\text{\rm D}\hh$,
\item the \dn\ $\,\vt\hs$ defined in\/ Section $\ref{bapr}$ is tangent to the 
$\,\tab\,$ factor, and coincides with \tvn\ of\/ $\,\tab$, 
\item the $\,\w\nnh$\pd\ connection for $\,\vt$, cf.\ Theorem\/ $\ref{csvhc}$, 
is the original\/ $\,\hs\text{\rm D}\hh$.
\end{enumerate}
Conversely, in any \cs\ \prd\ with\/ $\rwo$, every 
\pt\ has a connected neighborhood isometric to an open subset of a manifold\/ 
{\rm(\ref{hgt})} constructed above from some data {\rm(i)} -- {\rm(vi)}.
\end{theorem}
\begin{remark}\label{notwp}A general local form of a \cs\ metric with $\rwo\,$ 
is the warped product $\,h\hs+\fy^2\gm\,$ with 
$\,h=h^{\text{\rm D}}\hskip-2.5pt-2\tau\,$ and 
$\,\hs\text{\rm D}\hh\df=-\fy\nnh\ri^{\hs\text{\rm D}}$ on $\,\bs\,$ (for 
details, see Theorems~\ref{dimre}, \ref{clsdf} and Remark~\ref{ddfrh}). As 
shown in Section \ref{pfth}, this is equivalent to (\ref{hgt}). We chose to 
state the classification theorem in terms of (\ref{hgt}) to avoid using, in 
addition to $\,\tau$, yet another object (the function $\,\fy$), which is not 
a genuine parameter of the construction, as the local isometry type of the 
resulting metric does not depend on it.
\end{remark}

\section{Proof of Theorem~\ref{class}}\label{pfth}
The metric $\,g\,$ defined by (\ref{hgt}) has the lo\-cal-co\-or\-di\-nate 
expressions
\begin{equation}
\begin{array}{rl}
\mathrm{i)}\hskip10pt&
g\,=\,2\hs dp_j\,d\y^{\hs j}-\hs(2\hh\p_j\vg_{\!kl}^j\nh+2\tau_{kl}\nh+
\nnh\gm_{ab}\x^a\x^bR_{kl})\,d\y^k\hs d\y^{\hh l}+\hs\gm_{ab}\,d\x^a\hs d\x^b,
\\
\mathrm{ii)}\hskip10pt&
\fy^{-2}g\,=\,2\hs d\q_j\hs d\y^{\hs j}\hs
-\,2\fy^{-2}\tau_{jk}\,d\y^{\hs j}\hs d\y^k\hs+\,\gm_{ab}\,du^a\hs du^b.
\end{array}\label{ged}
\end{equation}
In (\ref{ged}\hh-i) $\,\y^{\hs j}\nnh,\p_j,\x^a$ are product coordinates for 
$\,\tab\,\times\,\mv$ formed by the coordinates $\,\y^{\hs j}\nnh,\p_j$ in 
$\,\tab\,$ obtained as usual from any given local coordinates $\,\y^{\hs j}$ 
in $\,\bs$, and linear coordinates $\,\x^a$ in $\,\mv$ corresponding to a basis 
$\,e_a$, with $\,\gm_{ab}=\lg e_a,e_b\rg$, while $\,\vg_{\!kl}^j$ and 
$\,R_{kl}$ stand for the components of $\,\hs\text{\rm D}\hs\,$ and its \rt\ 
$\,\ri^{\hs\text{\rm D}}\nnh$. Thus, (\ref{ged}\hh-i) is obvious from 
(\ref{hde}). As for \hbox{(\ref{ged}\hh-ii),} we start from the coordinates 
$\,\y^{\hs j}\nnh,\p_j$ for $\,\tab\,$ based on coordinates $\,\y^{\hs j}$ in 
$\,\bs\,$ that, instead of being arbitrary, are chosen (along with the 
function $\,\fy$) as in Lemma~\ref{pwmcf}. We then replace the resulting 
product coordinates $\,\y^{\hs j}\nnh,\p_j,\x^a$ in $\,\tab\,\times\,\mv$ (with 
$\,\x^a$ as before) by the new coordinates $\,\y^{\hs j}\nnh,\q_j,u^a$ such 
that $\,2\q_j=2\fy^{-2}\p_j+\gm_{ab}\x^a\x^b\fy^{-3}\hs\partial_j\fy\,$ and 
$\,u^a=\fy^{-1}\x^a\nh$. By (\ref{tge}), 
$\,\fy^{-2}(\hh h^{\text{\rm D}}\hskip-2.5pt-2\tau)
=2\hs d(\fy^{-2}\p_j)\,d\y^{\hs j}\hs
-\,2\hs\fy^{-2}\tau_{jk}\,d\y^{\hs j}\hs d\y^k\nnh$, while 
$\,\fy^{-2}(\gm-\theta\ri^{\hs\text{\rm D}})
=\fy^{-2}(\gm_{ab}\,d\x^a\hs d\x^b\nh
-\gm_{ab}\x^a\x^bR_{kl}\,d\y^k\hs d\y^{\hh l})$, and so (\ref{phg}\hh-c) gives 
$\,\fy^{-2}(\gm-\theta\ri^{\hs\text{\rm D}})
=\fy^{-2}\gm_{ab}\,d\x^a\hs d\x^b\nh
-\gm_{ab}\x^a\x^b\fy^{-1}d\hs[\hh\partial_j(\fy^{-1})]\,d\y^{\hs j}\nh$, 
proving (\ref{ged}\hh-ii).

In view of (\ref{ged}\hh-ii), $\,g\,$ is, locally, a warped product 
$\,h\hs+\fy^2\gm\,$ with the factor metrics 
$\,h=2\fy^2\hs d\q_j\hs d\y^{\hs j}\nh-2\tau_{jk}\,d\y^{\hs j}\hs d\y^k$ (in 
dimension $\,4$, coordinates: $\,\y^{\hs j}\nnh,\q_j$) and 
$\,\gm=\gm_{ab}\,du^a\hs du^b$ (dimension $\,n-4$, coordinates: $\,u^a$). 
Since $\,\gm_{ab}=\lg e_a,e_b\rg$ are constants, $\,\gm\,$ is flat. On the 
other hand, $\,h\,$ is isometric to the metric 
$\,h^{\text{\rm D}}\hskip-2.5pt-2\tau\,$ defined as in Section \ref{difo} 
using the data (i) - (iv) of Section \ref{tlst}. Namely, by (\ref{tge}), 
$\,h^{\text{\rm D}}\hskip-2.5pt-2\tau\,$ equals 
$\,2\fy^2\hs d\bq_j\hs d\y^{\hs j}\nh-2\tau_{jk}\,d\y^{\hs j}\hs d\y^k$ in 
coordinates $\,\y^{\hs j}\nnh,\bq_j$ formed by our $\,\y^{\hs j}$ and some 
$\,\bq_j$, while $\,\fy\,$ and $\,\tau_{kl}$ depend only on 
$\,\y^{\hs j}\nnh$. By Theorem \ref{clsdf}, $\,(N,h)\,$ is \cs\ and $\rwho$, 
so that Theorem~\ref{dimre} yields the same for (\ref{hgt}), while 
Theorems~\ref{clsdf} and~\ref{csvhc}(iv),\hskip2.2pt(v) imply (a) -- (c).

Conversely, let $\,(M,g)\,$ be \cs, with $\rwo$. Thus, locally, 
$\,g=h\hs+\fy^2\gm\,$ with $\,h,\gm,\fy\,$ as in Theorem~\ref{dimre}, and, by 
Theorem~\ref{clsdf}, $\,h=h^{\text{\rm D}}\hskip-2.5pt-2\tau$, for suitable 
$\,\hs\text{\rm D}\hs\,$ and $\,\tau$. Now (\ref{tge}) gives 
$\,g=2\fy^2d(\fy^{-2}\p_j)\,d\y^{\hs j}\nh
-2\hs\tau_{jk}\,d\y^{\hs j}\hs d\y^k\nh+\fy^2\gm_{ab}\,du^a\hs du^b$ in some 
coordinates $\,\y^{\hs j}\nnh,\fy^{-2}\p_j,u^a\nnh$, with constants 
$\,\gm_{ab}$ such that $\,\gm=\gm_{ab}\,du^a\hs du^b\nnh$. (Note that 
$\,\hs\text{\rm D}\hh\df=-\fy\nnh\ri^{\hs\text{\rm D}}\nnh$, by 
Remark~\ref{ddfrh}; hence, according to Remark~\ref{ddfrd}, $\,\fy\,$ is 
chosen as in Lemma~\ref{pwmcf}.) This gives \hbox{(\ref{ged}\hh-ii)} for 
$\,\q_j=\fy^{-2}\p_j$, and Theorem~\ref{class} follows.
\begin{remark}\label{plbck}For $\,(M,g)\,$ in (\ref{hgt}) and 
$\,(M,\hs g\hs'\hh)=(\tab\,\times\,\mv,\,\,h^{\text{\rm D}}\hskip-2.5pt
-2\tau\hs'\nh+\gm-\theta\ri^{\hs\text{\rm D}})$ constructed as in 
Section \ref{tlst} from the same data \hbox{(i) -- (vi),} but with (possibly) 
different $\,2$-ten\-sors $\,\tau\,$ and $\,\tau\hs'\nnh$, an isometry 
$\,\js\,$ of $\,(M,g\hs'\hh)\,$ onto $\,(M,g)\,$ can be defined by 
$\,\js(\cy,\eta,\x)=(\cy,\eta+\xi(\cy)-(d\lg L,\nh H\x\rg)_\cy\hs,H\x+L(\cy))\,$ 
(notation of Section \ref{riex}), for any function $\,L:\bs\to\mv$ with 
$\,\hs\text{\rm D}\hh dL=-\hh L\nh\ri^{\hs\text{\rm D}}\nnh$, any 
$\,\lr$-pre\-serv\-ing linear isomorphism $\,H:\mv\nh\to\mv\nnh$, and any 
$\,1$-form $\,\xi\,$ on $\,\bs\,$ such that 
$\,2(\tau-\tau\hs'\hh)=\lo(\xi+d\lg L,L\rg/4)\,$ (cf.\ (\ref{lxi})); 
$\,\lg L,\nh H\x\rg\,$ and $\,\lg L,L\rg\,$ are functions $\,\bs\to\bbR\hs$. 
In fact, $\,\js^*\nh g=g\hs'\nnh$, since $\,\js^*\nh g\,$ is obtained by 
replacing $\,\p_j$ and $\,\x^a$ in \hbox{(\ref{ged}\hh-i)} with 
$\,\p_j+\hh\xi_j-\gm_{ab}H^b_c\x^c\partial_j\hh L^b$ and $\,H^a_bx^b+L^a\nnh$. 
Cf.\ \cite[\S8]{patterson-walker}.

Let $\,\bs\,$ now be \sco. First, $\,(M,g\hs'\hh)\,$ is isometric to 
$\,(M,g)$, as one sees choosing $\,\js\,$ with 
$\,L=0,\hs H=\hs\text{\rm Id}\,$ and $\,\xi\,$ such that 
$\,2(\tau-\tau\hs'\hh)=\lo\xi\,$ (which exists by Theorem~\ref{slveq}(ii)). 
Secondly, for every $\,\cy\in\bs$, the set 
$\,\{\cy\}\times\hs\tcb\hh\times\hs\mv\,$ is contained in an orbit of the 
isometry group of $\,(M,g)$. In fact, $\,L\,$ with 
$\,\hs\text{\rm D}\hh dL=-\hh L\nh\ri^{\hs\text{\rm D}}$ and $\,\xi\hh'$ with 
$\,\lo\xi\hh'\nh=0\,$ realize all values $\,L(\cy)\in\mv$ and 
$\,\xi\hh'\hskip-3.4pt_\cy\in\tacb$, cf.\ Remarks~\ref{ddfrh}, \ref{detby} and 
Lemma~\ref{immer}.
\end{remark}

\section{Compactness the quotient surface}\label{cotq}
We do not know whether there exist any {\it compact\/} \ecs\ \mf s. However, 
some such \mf s do have a com\-pact\-ness-type property, namely, the leaf 
space of $\,\vt\hs$ is a globally well-defined (Haus\-dorff) closed \su:
\begin{example}\label{ctype}Given a closed surface $\,\bs$, an integer 
$\,n\ge4$, and a metric signature 
\hbox{$-$\hskip1pt$-$\hskip2pt$\dots$\hskip0pt$+$\hskip1pt$+$} with $\,n\,$ 
signs, containing at least two minuses and at least two pluses, there exists 
an \ecs\ \mf\/ $\,(M,g)\,$ of dimension $\,n\,$ such that
\begin{enumerate}
  \def\theenumi{{\rm\roman{enumi}}}
\item $g\,$ is not \rr, and has the prescribed signature\/ 
\hbox{$-$\hskip1pt$-$\hskip2pt$\dots$\hskip0pt$+$\hskip1pt$+\hs$,}
\item the leaves of the \dn\ $\,\vt\hs$ defined in Section $\ref{bapr}$ are the 
fibres of a \bd\ with the \ts\ $\,M\,$ and base $\,\bs$.
\end{enumerate}
\end{example}
In fact, we may choose $\,(M,g)\,$ to be the manifold (\ref{hgt}) obtained 
from the construction in Section \ref{tlst}, using $\,\lr\,$ of the 
appropriate signature, and $\,\hs\text{\rm D}\hs$, along with 
$\,\pm\hs\alpha\,$ and $\,\tau\,$ on our \su\ $\,\bs$, that satisfy the 
required assumptions: they exist by Theorem~\ref{eqcls}, while (i) follows 
from Theorem~\ref{class}(a), since $\,\hs\text{\rm D}\hs\,$ is not \rr.

\section{Further comments}\label{comm}
Let us consider triples $\,(M,g,x)\,$ formed by a \mf\ $\,M$, a \pt\ 
$\,x\in M$, and a \cs\ metric $\,g\,$ on $\,M\,$ having a fixed signature 
\hbox{$-$\hskip1pt$-$\hskip2pt$\dots$\hskip0pt$+$\hskip1pt$+$} with two or 
more minuses and two or more pluses, and satisfying the condition $\rwo$. 
(Note that the signature determines $\,\dim M$.) We call two such triples 
$\,(M,g,x)\,$ and $\,(M'\nnh,g\hh'\nnh,x\hh')$ {\it equivalent\/} if some 
isometry of a neighborhood of $\,x\,$ in $\,M\,$ onto a neighborhood of 
$\,x\hh'$ in $\,M'$ sends $\,x\,$ to $\,x\hh'\nnh$, and we refer to the set of 
equivalence classes of this relation as the {\it local moduli space\/} of \cs\ 
metrics of the given signature with $\rwo$.

We similarly define the local moduli space of \ea, \pftc s on \su s, to 
be the set of equivalence classes of quadruples 
$\,(\bs,\hs\text{\rm D}\hs,\pm\hs\alpha,\cy)$ formed by a \su\ $\,\bs\,$ with 
such a connection $\,\hs{\rm D}\hs\,$ (see Section \ref{eaco}), a 
$\,\hs\text{\rm D}\hs$\prl\ area element $\,\pm\hs\alpha\,$ on $\,\bs$, and 
a \pt\ $\,\cy\in\bs$. The equivalence relation is defined similarly, except 
that isometries are replaced by (local) unimodular affine diffeomorphisms.

By Theorem~\ref{class}, given 
\hbox{$-$\hskip1pt$-$\hskip2pt$\dots$\hskip0pt$+$\hskip1pt$+\hs$}, the former 
moduli space is in a natural one-to-one correspondence with the latter: the 
correspondence sends the equivalence class of $\,(M,g,x)\,$ to that of 
$\,(\bs,\hs\text{\rm D}\hs,\pm\hs\alpha,\cy)$ obtained in Theorem~\ref{csvhc}, 
with $\,\cy=\pmb(x)$.

Our next comment concerns Ric\-ci-\-re\-cur\-rence. In Theorem~\ref{class} one 
cannot simply replace the condition $\rwo\,$ by `not being \rr' (both for 
$\,(M,g)\,$ and for $\,\hs\text{\rm D}\hs\,$ in (iii) of Section \ref{tlst}), 
and still obtain a classification result with an analogous final clause. 
Namely, there is no principle of unique continuation, either for \cs\ metrics 
$\,g\,$ with $\rwo$, or for \pftc s $\,\hs\text{\rm D}$ \hbox{on \su s.} 
(Thus, neither of the two can in general be made real-analytic by a suitable 
choice of local coordinates.) In fact, both $\,g,\hs\text{\rm D}\hs\,$ can be 
\rr\ on some nonempty open set, without being so everywhere. Examples are 
immediate from Theorem~\ref{rcrec}: it suffices to deform a cylinder surface 
so as to make it non-cylindrical just in a small subset. The construction of 
Section \ref{tlst}, applied to the corresponding cen\-tro\-af\-fine connection 
$\,\hs\text{\rm D}\hs$, then yields a metric $\,g\,$ with the stated property.

Finally, any \cs\ \prd\ $\,(M,g)\,$ of dimension $\,n\ge4\,$ with $\rwo$, and 
any pseu\-\hbox{do\hs-}Euclid\-e\-an inner product $\,\lr'$ on a $\,k$\diml\ 
real \vs\ $\,\mv'\nnh$, treated as a \trinv\ metric $\,\gm\hs'$ on $\,\mv'\nnh$, 
give rise to the \cs\ metric 
$\,g\hh'\nh=g-(n-2)^{-1}\theta\hs'\hskip-1.6pt\ri+\gm\hs'$ with $\rwo\,$ on 
$\,M\times\mv\nnh$, where $\,\ri\,$ is the Ricci tensor of $\,g\,$ and 
$\,\theta\hs':\mv\to\bbR$ is defined by $\,\theta\hs'\nh(v)=\lg v,v\rg'\nnh$. 
(Notation of (\ref{hgt}).) In fact, $\,g\hh'$ can also be constructed as in 
Section \ref{tlst}, since, by Theorem~\ref{csvhc}(ii), we may replace 
$\,\ri^{\hs\text{\rm D}}$ in (\ref{hgt}) with $\,(n-2)^{-1}\ri$.

\end{document}